\documentclass[12pt,leqno]{report}    %   LEFT NUMBERRING OF EQUATIONS
\usepackage{amsmath}
\usepackage{latexsym}
\usepackage{amssymb}
\usepackage{amsmath,amssymb,amsfonts}
\usepackage{graphicx}
\usepackage{psfrag}
\usepackage{mathabx}
\usepackage{bbm}
\usepackage[small,nohug,heads=vee]{diagrams}
\diagramstyle[labelstyle=\textstyle]

\def\ca{{\mbox{a}}}
\def\sqb{{\mbox{b}}}
\def\sqp{{\mbox{p}}}

\def\sqa{{\mbox{a}}}
\def\sqq{{\mbox{q}}}
\def\sqm{{\mbox{m}}}

\numberwithin{equation}{section}

\newcommand{\leftmapsto}{ \leftarrow\hspace{-2pt}\mapstochar}
\newcommand{\cX}{ {\cal X} }
\newcommand{\cD}{ {\cal D} }
\newcommand{\cL}{ {\cal L} }
\newcommand{\cS}{ {\cal S} }
\newcommand{\cG}{ {\cal G} }
\newcommand{\cR}{ {\cal R} }
\newcommand{\cC}{ {\cal C} }
\newcommand{\cO}{ {\cal O} }
\newcommand{\cY}{ {\cal Y} }
\newcommand{\cZ}{ {\cal Z} }
\newcommand{\cN}{ {\cal N} }
\newcommand{\cK}{ {\cal K} }
\newcommand{\cW}{ {\cal W} }

\newcommand{\Fbb}{\mathbb{F}_{\bullet}}
\newcommand{\be}{ \begin{equation} }
\newcommand{\ee}{ \end{equation} }

\begin{document}
\pagestyle{plain}

\addcontentsline{toc}{chapter}{Contents}
\addcontentsline{toc}{section}{Contents}
\addcontentsline{toc}{subsection}{Contents}

\title{Non-Additive Prolegomena\\
{\small (to any future Arithmetic that will be able to present itself as a Geometry)}}

\author{Shai M. J. Haran\\
Technion- Israel Institute of Technology\\
haran@tx.technion.ac.il}

\maketitle

\chapter*{Introduction} 
Andr\'{e} Weil (in a letter to his sister Simone Weil \cite{W40}) compares the situation of the one-dimensional objects of mathematics to that of the Rosetta stone. The Rosetta stone contains three languages talking about one reality -
two languages (Greek and Demotic) were well understood, and one (Hieroglyphic) was a complete mystery.
The situation in mathematics is somewhat different in that we have one language - the language of addition and multiplication, 
i.e. commutative rings - and with this language we are talking about three similar but different realities, two Geometric realities (curves over finite fields and curves over $\mathbb C$, i.e. compact Riemann surfaces) 
are well understood, and one Arithmetic reality (number fields) is a complete mystery.
Comparing Arithmetic and Geometry we find that we have the wrong language, 
and especially addition is behind the three basic problems of Arithmetic:
{ \bf the real prime}, {\bf the arithmetical plane}, and {\bf the absolute point}.

\noindent {\bf The real prime:}
In the geometry of curves over a field $k$ we realize that if we want to have theorems we need to pass from affine to projective geometry, and in particular we have to add the point at infinity $\infty$ to the affine line 
${\mathbb A}^1_k$ and look at the 
projective line
${\mathbb P}^1_k= {\mathbb A}^1_k \cup \{ \infty \}$.
In our language this translates into the fact that the fraction 
field of the polynomial ring $k[X]$,
the field of rational functions $k(X)$, 
embedds into the field of Laurent series
$K_{\alpha}= k((x-\alpha))$ for each $\alpha \in k$, 
as well as into the field $K_{\infty}=k((\frac{1}{x}))$.
These fields have one dimensional local subrings, 
the subrings of power series 
${\cal O}_{\alpha}= k[[x-\alpha]]\subseteq K_{\alpha}$
and ${\cal O}_{\infty}= k [[\frac{1}{x}]]\subseteq K_{\infty}$.
Similarly in arithmetic we have to add the real and 
complex primes of a number field to the finite primes, 
and in particular add the real prime $\eta$ to the finite primes 
$spec({\mathbb Z})= \{2, 3, 5, 7, 11 \ldots \}$.
The fraction field of the ring of integers $\mathbb Z$,
the field of rational numbers $\mathbb Q$,
embedds into the field of $p$-adic numbers ${\mathbb Q}_p$
for each finite prime $p$, as well as into the field of real numbers 
${\mathbb Q}_{\eta}= {\mathbb R}$.
For the finite primes $p$, 
the fields ${\mathbb Q}_p$ 
contains the one dimensional local subring of $p$-adic 
integers ${\mathbb Z}_p \subseteq {\mathbb Q}_p$.
But for the real prime we find that ${\mathbb Z}_{\eta}$, 
the interval $[-1, 1]$,
is not closed under addition!

\noindent {\bf The arithmetical plane:}
Having two geometrical objects we can take their product, 
and a product of two curves gives us a surface.
In particular, the product of the affine line with itself 
gives us the affine plane,
${\mathbb A}_k^1 \times {\mathbb A}_k^1 \cong {\mathbb A}_k^2$.
This translates in our language to the fact that 
$k[x] \otimes_k k[x] \cong k[x_1, x_2]$.
But in arithmetic we find that
${\mathbb Z} \otimes {\mathbb Z}  = {\mathbb Z} $, 
the arithmetical plane reduces to its diagonal!

\noindent {\bf The absolute point:}
The category of $k$-algebras has $k$ as an initial object, 
hence in geometry over $k$ we have a point 
$spec(k)$ as a final object.
Having addition in our language forces the integers ${\mathbb Z}$ 
to be the initial object of the category of commutative rings,
hence $spec({\mathbb Z})$ is our final geometric object,
and we are missing the absolute point $spec({\mathbb F}_1)$,
where the "field with one element"
${\mathbb F}_1$ is the common field of all finite fields 
${\mathbb F}_p , p=2,3,5,7 \ldots$

\vspace{10pt}

It was Kurokawa, Ochiai, and Wakayama \cite{KOW} who were the first to suggest abandoning addition,
and work instead with the language of
multiplicative monoids. This idea was further described in Deitmar 
\cite{De}, but note that the spectra of monoids always looks like the spectra of a local ring:
the non-invertible elements are the unique maximal ideal.
For Kurokawa there is also a ``zeta world'' of analytic functions that encode geometry,
where the field with one element ${\mathbb F}_1$,
is encoded by the identity function of ${\mathbb C}$,
see Manin \cite{M}.
Soul\'{e} in \cite{S} tries to capture ${\mathbb F}_1$ 
by defining ``${\mathbb F}_1$-varieties'' as a subcollection of 
${\mathbb Z}$-varieties.
(It was C. Soul\'{e} who awakened me from my dogmatic slumber).

In \cite{H07} we gave our first non-additive language for geometry 
based on the language of ${\mathbb F}$-rings
and ring-categories. A ring-category is a category $A$ 
with two symmetric monoidal structures,
``direct-sum'' $\oplus: A \times A \rightarrow A$, and 
``tensor-product'' $\otimes:A \times A \rightarrow A$,
with $\otimes$ naturally distributive over $\oplus$.
The category of ring-categories has ${\mathbb F}_1$ as its initial object.
As a category ${\mathbb F}_1$ has for objects the 
finite sets and for morphisms the partial bijections;
the operation $\oplus$ corresponds to disjoint sum; 
the operation $\otimes$ corresponds to direct-product.
An ${\mathbb F}$-ring is a ring-category $A$ such that the canonical map 
${\mathbb F}_1 \rightarrow A$ is the identity on objects.
The motivation for this language came from the hint of the real prime $\eta$:
while ${\mathbb Z}_{\eta} = [-1 , 1]$ is not closed under addition, 
it is closed under contraction, and we have by the fundamental 
Cauchy-Schwartz inequality that 
$(x, y) = x_1 \cdot y_1 + \cdots + x_n \cdot y_n \, \,  \mbox{is in} \, \, {\mathbb Z}_{\eta}= [-1, 1]$
if the vectors $x$, $y$ are in $({\mathbb Z}_{\eta})_n$, 
that is $|x|_{\eta} \leq 1$ and $|y|_{\eta} \leq 1$, 
with the $l_2$-norm 
$|x|_{\eta}  = (\sum\limits_i |x_i|^2)^{\frac{1}{2}}$.

In \cite{Du} Durov gives a more algebraic language of  generalized rings based on monads of sets.
But this forces him to replace the $l_2$-norm by the $l_1$-norm,
and thus he gets the wrong results at the real prime $\eta$.
Thus in Durov's language $GL_n({\mathbb Z}_{\eta})$ 
is the finite group of symmetries of the $l_1$-polytope
$\{  x \in {\mathbb R}^n , \sum\limits_i |x_i| \leq 1\}$, 
while it should be the orthogonal group $O_n$
(and the unitary group $U_n$ for a complex prime)
 which is the symmetry group of the $l_2$-ball
\[\{  x \in {\mathbb R}^n , \sum\limits_i |x_i|^2 \leq 1\} = ({\mathbb Z}_{\eta})_n\]
Indeed, MacDonald \cite{Mac} gives a $q$-analog interpolation 
between the zonal-spherical functions on 
$GL_n({\mathbb Q}_p) / GL_n({\mathbb Z}_p) $ and 
the zonal-spherical functions on $GL_n({\mathbb R}) / O_n$ and 
$GL_n({\mathbb C}) / U_n$.
Similarly,
there is a $q$-analog interpolation where the representation theory
of the quantum group $GL_n$ gives in the ``p-adic limit'' 
the representation theory of $GL_n({\mathbb Z}_p)$,
and in the ``real or complex limits'' the representation theory of
$O_n$ and $U_n$,
see \cite{H08}, $[0]$.
But Durov succeeded in proving the negative result that the tensor 
product of $\mathbb Z$ with itself, 
both in his category of generalized rings, and in our category of 
$\mathbb F$-rings, reduces again to $\mathbb Z$, 
and thus the ``arithmetical-plane'' still reduces to its diagonal.

\vspace{10pt}

Here we give a new language for geometry based on a new concept of a generalized ring.
A generalized ring $A$ is given by a sequence of sets with embeddings 
$A_0 = \{0\} \subseteq A_1 \subseteq A_2 \subseteq \cdots \subseteq A_n \subseteq \cdots$ 
and projecions 
$\pi_n : A_{n+1} \rightarrow A_n$, with $\pi_n|_{A_n} = id_{A_n}$.
The set $A_n$ carries an action by the symmetric group $S_n$, 
and the embeddings / projections are covariant.
There are two covariant operations, the operation of multiplication
\[\begin{array}{lcrl}
A_k \times (A_{n_1} \times \cdots \times A_{n_k})& \longrightarrow& A_n & , n=n_1+\cdots n_k\\
a\quad , \quad b= (b_1, \ldots , b_k)& \mapsto & a \circ b&
\end{array}\]
and the operation of contraction
\[\begin{array}{lcr}
A_n \times (A_{n_1} \times \cdots \times A_{n_k})& \longrightarrow& A_k\\
a\quad , \quad b & \mapsto & (a,b)
\end{array}
\]

These operations are required to satisfy the following axioms:
\begin{description}
\item[associativity:] $(a \circ b ) \circ c = a \circ ( b \circ c)$;
\item[adjunction:] $((a,b),c) = (a, c\circ b) \quad , \quad (a,(b,c)) = (a \circ c , b)$;
\item[linearity:] $a \circ (b,c) = (a \circ b , c) \quad , \quad (a,b) \circ c = ( a \circ \tilde{c} , \tilde{b})$,\\
where $\tilde{c}$ is obtained from $c$ by diagonal embedding, \\
and similarly $\tilde{b}$ from $b$; 
\item[unit:] have $1 \in A_1$ satisfying\\
$ a = 1 \circ a = a \circ 1 = (a,1)$
\end{description}
It follows from the axioms that $A_1$ is a commutative monoid with unit; 
$(A_1)^n$ acts on $A_n$;
and $a^t=(1,a)$ is an involution of $A_1$.
We say that $A$ is 
\underline{self-adjoint} if we have $a^t = a $ for all $a \in A_1$.
We say that $A$ is \underline{commutative} if we have 
$a \circ \tilde{b} = b \circ \tilde{a}$ for $a \in A_n$, $b \in A_m$.
We develop a geometry for self-adjoint generalized rings.
We do not need the commutativity axioms, but the self-adjointness is 
indispensible.
On the other hand, all our examples are commutative, and the self-adjointness axiom complicates the description of ``formulas'' in our language (i.e. of the free object).

\vspace{10pt}

In Chapter \S 1 we give the precise definition of a generalized ring,
in a more functorial (but equivalent) form.

\vspace{10pt}

In Chapter \S 2 we give our basic examples of generalized-rings:
\begin{itemize}
\item[] the initial object - the field with one element $\mathbb F$;
\item[] the generalized-ring ${\cG}(A)$ attached to a commutative (semi-)ring $A$;
\item[] the local generalized-ring ${\cal O}_{\eta}$ associated to a real or complex prime $\eta$ and its associated residue field ${\mathbb F}_{\eta}$;
\item[] the generalized-ring ${\mathbb F}[M]$ associated to a monoid $M$;
\item[] the free generalized-rings $\Delta^w$ (resp. $\Delta$ ) 
representing the functors $A \mapsto A_w$ 
(resp. $A \mapsto \lim\limits_{\stackrel{\longleftarrow}{\mathbb F}}  A_w$).
\end{itemize}

In Chapter \S 3 we consider equivalence-ideals, ideals, homogeneous-ideals,
and the correspondences between them.
There is also the usefull notion of $h$-ideals, which for self-adjoint generalized-rings coincide with homogeneous-ideals.

\vspace{10pt}

In Chapter \S 4 we consider the topological space $spec(A)$,
consisiting of the primes of the generalized-ring $A$, with its Zariski-topology,
and its dense subset $Espec(A)$ consisiting of the stable primes.
We show the contravariant functors $A \mapsto spec(A)$,
and $A \mapsto Espec(A)$, 
take generalized-rings and homomorphisms into compact, sober (Zariski),
topological spaces and continuous maps.

\vspace{10pt}

In Chapter \S 5 we consider localizations of generalized-rings,
 and obtain the sheaf ${\cO}_A$ of generalized-rings over 
$spec(A)$ for a self-adjoint generalized-ring $A$.

\vspace{10pt}

In Chapter \S 6 we consider the category 
${\cal LGRS}$ of locally-generalized-ringed-spaces,
its objects are pairs $(X, \cO_X)$ 
of a topological space $X$ and a sheaf $\cO_X$ over $X$ of (self-adjoint) generalized-rings with local fibers.
We show that the functors $A \mapsto (specA, \cO_A)$ and 
$(X, \cO_X) \mapsto \cO_X(X)$
are adjoint.
The category of Grothendieck-generalized-schemes is the full subcategory 
${\cal GGS} \subseteq {\cal LGRS}$ consisting of 
$(X, \cO_X) \in {\cal LGRS}$
such that there exists an open cover $X = \bigcup\limits_i U_i$ with 
$(U_i, \cO_X|_{U_i}) \cong spec \cO_X (U_i)$.
The category of generalized-schemes ${\cal GS}$ 
is the category of pro-objects of ${\cal GGS}$.
We show that in ${\cal GS}$ there is a compactification 
$\overline{spec({\mathbb Z})}$ of $spec({\mathbb Z})$,
and similarly for any number-field $K$ there is a 
compactification $\overline{spec(\cO_K)}$ of the finite-primes 
$spec(\cO_K)$.
These compactifications follow very closely the ones given in \cite{H07}, (and reproduced in \cite{Du}).

\vspace{10pt}

In Chapter \S 7 we consider the tensor-product of generalized-rings, and the
resulting fiber-products in the categories ${\cal GS} \supseteq {\cal GGS}$ of\\ (Grothendieck) generalized-schemes.
In particular, we have in ${\cal GS}$ the arithmetical plane 
$\overline{spec(\mathbb Z)} \times_{{\mathbb F} [\pm 1]} \overline{spec(\mathbb Z)}$,
and its open and dense affine generalized scheme given by the spectrum of 
${\cal G}(\mathbb Z) \otimes_{{\mathbb F} [\pm 1]} {\cal G}(\mathbb Z)$.
We give an explicit description of  
 the generalized ring 
${\cal G}(\mathbb N) \otimes_{\mathbb F} {\cal G}(\mathbb N)$,
and of ${\cal G}(\mathbb N) \bigotimes\limits_{\mathbb F}
        {\cal G}(\mathbb N) \bigotimes\limits_{\mathbb F}
         \mathbb F [\pm 1]$
which map surjectively onto
$\cG(\mathbb Z ) \bigotimes\limits_{\mathbb F[\pm 1]}
\cG(\mathbb Z )$.

\vspace{10pt}

In Chapter \S 8 we sketch the theory of divisors. We define the sheaf of meromorphic functions in (\ref{sec8.1}), and the abelian group of (Cartier) divisors in (\ref{sec8.2}).
In (\ref{sec8.3}) we mention the associated invertible 
$\cO_X$-module,
(although we do not enter into the linear
and homological algebra of modules 
in this paper).
In (\ref{sec8.4}) we define effectivenness of a divisor.
In (\ref{sec8.5}) we define the ordered abelian group
$Div(\cX)$ for 
$\cX = \{ \cX_N \}$
a generalized scheme.
The point is that while we have some philosophies about the distribution of the primes within discrete data
(e.g. conjugacy classes in Galois groups),
we are completely ignorant about the distributuion of the primes (or prime powers) within the continuum ${\mathbb R}^+$.
We feel we need to recreate the continuum
${\mathbb R}^+$ in an arithmetical way
(as $Pic(\overline{spec {\mathbb Z}})$),
cf. (\ref{sec8.6}).
In (\ref{sec8.7}) we
enter the twilight-zone and conjecture the existence of
intersection numbers for divisors and of
 ''Frobenus divisors'' on
$\overline{spec {\mathbb Z}} \prod\limits_{{\mathbb F} [\pm 1]} 
 \overline{spec {\mathbb Z}}$.
This is our liebster Jugendtraum of \cite{H89}-
to give a Scheme for the proof of the Riemann Hypothesis.

\vspace{10pt}

The reader should be warned that while admitedly, 
we have used the Riemann 
Hypothesis as a compass in the developement of this language,
there are many open problems in arithmetics which
are more ''elementary'' than the Riemann Hypothesis
(e.g. an analog of Hurwitz genus formula
for the map
$\underline{f} : \overline{spec {\mathbb Z}} \rightarrow  
{\mathbb P}^1$ 
associated with
$f \in {\mathbb Q}^*$,
cf. (6.4.13),
will give the ABC conjecture),
and there are still many more elementary problems in our theory.
We need to redo 
''non-additive commutative algebra'',
and as a first step,
we need to understand the combinatorics of the
''non-additive polynomials'',
which are given by finite rooted trees with certain labelings and boundary
identifications taken up to commutativity and self-adjunction.

\setcounter{equation}{0}
\setcounter{chapter}{-1}
\chapter{Notations}
For a category $ C $ we write $ X \in C $ for ``$ X $ is an object of
$ C $'', and we let $ C(X,Y) $ denote the set of maps in $ C $ from $
X $ to $ Y $. We denote by $Set_0$ the category with objects sets
$X$ with a distinguished element $O_X\in X$, and with maps
preserving the distinguished elements
\begin{equation}    %%%%    0.1     %%%%
 Set_0(X,Y) = \Big\{f\in Set(X,Y) \;,\; f(O_X)=O_Y\Big\}.    
\end{equation}
The category $Set_0$ has direct and inverse limits. The set
$[o]=\big\{o\big\}$ is the initial and final object of $Set_0$. For
$f\in Set_0(X,Y)$ we have
\begin{equation}    %%%%    0.2     %%%%
   ker f=f^{-1}(O_Y) \quad , \quad coker f=Y/f(X) .
\end{equation}
There is a cannonical map
\begin{equation}    %%%%    0.3     %%%%
   coker \ ker f= X/f^{-1}(O_Y) \quad \longrightarrow \quad ker \ coker  f=f(X)  .
\end{equation}
We denote by $\mathbb{F}_0$ the subcategory of $Set_0$ with objects
the \underline{finite} sets, and with maps
\begin{equation}\begin{array}{cl}
   \mathbb{F} _0(X,Y) & = \big\{ f \in Set_0(X,Y),\ coker \ ker f
                    \xrightarrow{\sim} ker \ coker  f \\
& \hspace{5cm}     \mbox{ is  an  isomorphism} \big\} \\
                  & = \big\{ f \in Set_0(X,Y),\; f|_{X \setminus
                  f^{-1}(O_Y)}\, \mbox{ is  an  injection} \big\}.
\end{array}
\end{equation}
We let $ Set_{\bullet} $ denote the category with objects sets and with
partially defined maps
\begin{equation}    %%%%    0.5     %%%%
   Set_{\bullet}(X,Y)= \underset{X' \subseteq X}{\perp\mkern-7mu\perp} Set(X',Y). 
\end{equation}
Thus to $ f \in Set_{\bullet}(X,Y) $ there is associates its domain $ D(f)
\subseteq X $, and $ f \in Set(D(f),Y) $.\\
We have an isomorphism of categories
\begin{equation}    \label{eq:0.6}   %%%%    0.6     %%%%
 Set_0 \textstyle \stackrel{\thicksim }{ \longleftrightarrow } Set_{\bullet}
\end{equation}
given by
\begin{equation}    %%%%    0.7     %%%%
\begin{array}{ll}
X & \longmapsto X_+ = X \smallsetminus \big\{ O_X \big\} \\
f & \longmapsto f_+ \;,\; D(f_+) = X \smallsetminus f^{-1}( O_Y );
\end{array}
\end{equation}
and inversly
\begin{equation}    %%%%    0.8     %%%%
\big\{ O_X \big\}  \coprod X  = X_0 \leftmapsto X 
\end{equation}

\[\left.\begin{array}{cr}
x \in D(f): & f(x)\\
x \not \in D(f): \; & O_Y
\end{array}\right\}
 = f_0 (x) \leftmapsto f
\]
We let $ \mathbb{F}_{\bullet} $ denote the subcategory of $ Set_{\bullet} $
corresponding to $ \mathbb{F}_0 $ under the isomorphism
(\ref{eq:0.6}), it has objects the finite sets, and maps are the
partial bijections
\begin{equation}    %%%%    0.9     %%%%
   \mathbb{F}_{\bullet}(X,Y) = \big\{ f:D(f) \xrightarrow{\thicksim} f(X)
\big.  \ bijections,\  \big. D(f) \subseteq X,\ f(X) \subseteq Y
\big\} . 
\end{equation}
To avoid problems with set theory we shall work with a countable
model of $ \mathbb{F}_{\bullet} $ that contains $ [0] = \varnothing $ (the
empty set, the initial and final object), $ [1] = \big\{ 1 \big\}, \
\ldots, \ [n] = \big\{ 1,\ldots,n \big\} \ ,\ldots $ and is closed
under the operations of pull-back and push-out. In particular $
\mathbb{F}_{\bullet} $ is closed under disjoint union, the categorical sum,
which we denote by $ X \oplus Y $ ; it is associative, commutative,
and has $ [0] $ as unit. Moreover, $ \mathbb{F}_{\bullet} $ is closed under
the ``tensor-product'' operation
$$ X \otimes Y = \big\{ (x,y) , x \in X , y \in Y \big\}
  = coker \big\{ X \perp\mkern-7mu\perp  Y \rightarrow X \top\mkern-7mu\top Y \big\} ; $$
it is associative, commutative, and has $ [1] $ as unit, and is
distributive over $\oplus$. With these two operations, $\oplus$ and
$\otimes$, $ \mathbb{F}_{\bullet} $ is a ``ring category'' in the sense of
\cite{H07}. In fact, $ \mathbb{F}_{\bullet} $ is the initial object of the
category \emph{ RingCat} of ring categories. The $ \mathbb{F} $-Rings (on
which the geometry of \cite{H07} is based upon) are the ring categories A
$\in$ \emph{RingCat} such that the cannonical map $ \mathbb{F}_{\bullet}
\rightarrow A $ is bijection on objects. Inspired by \cite{S}, 
we refer to $ \mathbb{F}_{\bullet} $ as the ``field with
one element''. The category $ \mathbb{F}_{\bullet} $ has involution
\begin{equation}    %%%%    0.10    %%%%
\label{eq:0.10}
\begin{array}{c}
\mathbb{F}_{\bullet} (X,Y) \xrightarrow{\thicksim} \mathbb{F}_{\bullet} (Y,X) \\
\left( f:D(f) \stackrel{\thicksim}{\rightarrow} f(X) \right)  \longmapsto
f^t: f(X) \xrightarrow{\thicksim} D(f).
\end{array}
\end{equation}
We usually let $ X,Y,Z,W $ denote objects of $ \mathbb{F}_{\bullet} $,
without explicitly saying so, and when we consider ``$ Set_{\bullet}(X,Y) $''
it is usually implicitly assumed that $ X,Y \in \mathbb{F}_{\bullet} $.

\setcounter{section}{0}
\chapter{The Definition of Generalized Rings}
A \emph{generalized ring} $A$ consists of the following
structure:

\setcounter{section}{1}
\setcounter{equation}{0}

\vspace{10pt}

\noindent (1.1) A \emph{functor} $ A \in (Set_0)^{\mathbb{F}_{\bullet}}$. 
Thus for $ X
\in \mathbb{F}_{\bullet}$, we have $ A_X \in Set_0 $, and for $ f \in
\Fbb(X,Y)$, we have $ f_A \in Set_0(A_X,A_Y)$, such that
\begin{equation} 
 (g \circ f)_A = g_A \circ f_A \quad , \quad (id_X)_A = id_{A_X}
\end{equation}
In particular, for an injection $ f:X \hookrightarrow Y $ in
$\mathbb{F}_{\bullet}$, 
we have an injection $f_A:A_X \hookrightarrow A_Y$
in $ Set_0$, and a surjection
$ f^t_A:A_Y \twoheadrightarrow A_X $ in $ Set_0$, 
$ f^t_A \circ f_A = id_{A_X}$, 
and we usually identify $A_X$ as a subset of $A_Y$,
and write $a|_X$ for $f^t_A(a) \;,\; a \in A_Y$.
\begin{equation}\hspace{-2in}
\mbox{We assume:} \, A_{[0]} = \{ 0 \}
\end{equation}
For $f \in Set_{\bullet}(X,Y)$ , with $X,Y \in \Fbb$, we define
\begin{equation} 
A_f = \prod\limits_{y \in Y} A_{f^{-1} (y)}
    =  \left\{ b=( b^{(y)})_{y \in f(X)} \;,\;\; b^{(y)}
       \in A_{f^{-1} (y)} \subseteq A_X \right\} 
\end{equation}

\setcounter{section}{2}
\setcounter{equation}{0}
We have an operation of \emph{multiplication} for $f \in
Set_{\bullet}(X,Y)$,
\begin{equation}
\circ : A_Y \times A_f \rightarrow A_X \quad , \quad 
        a,b \mapsto a \circ b = a \circ_f b
\end{equation}
and we have an operation of \emph{contraction}
\begin{equation}
 (\, , \, ) : A_X \times A_f \rightarrow A_Y \quad , \quad 
       a,b \mapsto (a , b) = (a , b)_f. 
\end{equation}
We assume these operations respect the zero elements
\begin{equation} 
  0_Y \circ b = 0_X = a \circ 0_f, \,
     \mbox{here} \, 
     0_f = ( 0_{f^{-1}(y)} ) \in A_f  
\end{equation}
\begin{equation}
 (0_X,b) = 0_Y = (a,0_f).  
\end{equation}
Thus for the cannonical map $  c_X \in Set_{\bullet}(X,[1])$, 
the unique map with $  D(c_X) = X $, 
we have $ A_{c_X} \cong A_X$, 
and we get the operations
\begin{equation}
  \label{eq1.2.5}
\circ : A_{[1]} \times A_X \rightarrow A_X   , \quad \mbox{and}
\end{equation}
\begin{equation}
\label{eq1.2.6}
 (\, , \, ) : A_X \times A_X \rightarrow A_{[1]}.  
\end{equation}
For the identity map $  id_X \in Set_{\bullet}(X,X)$, 
we have\\ $  A_{id_X} = \prod\limits_{x \in X} A_{ \{ x\}} \cong (A_{[1]})^X$, 
and we get the operations
\begin{equation}
\label{eq1.2.7}
 \circ  \quad \mbox{and} \quad (\, , \, ):   A_X \times ( A_{[1]} )^X
\rightrightarrows A_X.  
\end{equation}
Before we can write the axioms satisfied by $ \circ $ and $(\, , \, )$,
we have to extend the definition of these operations. 
For $  f \in Set_{\bullet}(X,Y)$, \\
$  g \in Set_{\bullet}(Y,Z)$, and 
$ z \in Z$, we obtain by
restricting $ f $ the map
\begin{equation}
 f|_z \in Set_{\bullet} \left( ( g \circ f )^{-1} (z) \;,\; g^{-1} (z) \right)
\end{equation}
For $  y \in g^{-1} (z) $ we have identification of fibers
$  ( f|_z )^{-1} (y) = f^{-1} (y)$, 
and we have a projection
\begin{equation}
  A_f \twoheadrightarrow A_{f|_z} \quad,\quad
b = (b^{(y)}) \mapsto b|_z = (b^{(y)})_{y \in g^{-1} (z)}
\end{equation}
The extended operations are given by
\begin{equation}
  \circ : A_g \times A_f \rightarrow A_{g \circ f}  
\end{equation}
\[ a \; , \; b \mapsto a \circ _f b  \quad \mbox{with} \quad
    ( a \circ _f b )^{(z)} = a^{(z)} \circ_{f|_z } (b|_z ) 
\]
and by
\begin{equation}
 A_{g \circ f} \times A_f \rightarrow   A_g 
\end{equation}
\[  c \; , \; b \mapsto (c \;,\; b)_f  \quad \mbox{with} \quad
  (c \;,\; b)_f ^{(z)} = \left( c^{(z)} \;,\; (b|_z) \right)_{f|_z}. 
\]
We can now write the axioms.

\setcounter{section}{3}
\setcounter{equation}{0}

\vspace{10pt}
\noindent{\bf (1.3) Associativity:}
  \label{sec1.3}
  For $  W \xleftarrow{h} Z \xleftarrow{g} Y \xleftarrow{f} X  $ in
$ Set_{\bullet}$, 
and for $d \in A_h$, $c \in A_g$, $b \in A_f$, we have in
$ A_{h \circ g\circ f}$:
\begin{equation}
    d \circ (c \circ b) = (d \circ c ) \circ b  
\end{equation}

\setcounter{section}{4}
\setcounter{equation}{0}

\vspace{10pt}
\noindent{\bf (1.4) Left-Adjunction:}
  \label{sec1.4}
For  $  W \xleftarrow{h} Z \xleftarrow{g} Y \xleftarrow{f} X  $ in
$ Set_{\bullet}$,  and for $d \in A_{h \circ g\circ f}$, $ a \in A_g$,
$c \in A_f$, we have in $ A_h$:
\begin{equation}
 ( d , a \circ c) = \left( (d , c ) , a \right)  
  \label{eq1.4.1}
\end{equation}

\setcounter{section}{5}
\setcounter{equation}{0}

\vspace{10pt}
\noindent{\bf (1.5) Right-Adjunction}
\label{sec1.5}
 For $  W \xleftarrow{h} Z
\xleftarrow{g} Y \xleftarrow{f} X$ in $ Set_{\bullet}$,  
and for $d \in A_{h \circ g}$, $ a \in A_{g\circ f}$, $c \in A_f$, we have
in $ A_h$:
\begin{equation}
\label{eq1.5.1}
  ( d \circ c , a) = \Big( d , (a,c) \Big) 
\end{equation}

\setcounter{section}{6}
\setcounter{equation}{0}

\vspace{10pt}
\noindent{\bf (1.6) Left-Linear}
  \label{sec1.6}
 For $  W \xleftarrow{h} Z \xleftarrow{g} Y \xleftarrow{f} X$ 
     in $ Set_{\bullet}$,  
and for $d \in A_h$, $ a \in A_{ g\circ f}$, $c \in A_f$, 
we have in $ A_{h \circ g}$:
\begin{equation}
\label{eq1.6.1}
   ( d \circ a , c) = d \circ (a,c)  
\end{equation}

\setcounter{section}{7}
\setcounter{equation}{0}

\vspace{10pt}
\noindent{\bf (1.7) Right-Linear}
\label{sec1.7}
 For $  Z \xrightarrow{g} Y \xleftarrow{f} X  $
 in $ Set_{\bullet}$, we form the cartesian square,
%\begin{minipage}{2.5in}
\begin{equation}
\label{eq1.7.1}
%"picture"    :(     %%%%%
%
%$$      Z  \underset{Y}{\pi}  X    $$
%$$    \tilde{f} \swarrow \qquad \qquad \searrow \tilde{g}    $$
%$$    Z  \qquad \qquad \qquad \qquad X    $$
%$$    g \searrow  \qquad \qquad  \swarrow f   $$
%$$    Y     $$
%
  Z  \prod\limits_{Y}  X = \left\{ (z,x) \in D(g) \times D(f) ,
           g(z) = f(x)  \right\} 
\end{equation}
%\end{minipage} \hspace{10pt} \begin{minipage}{2in}
\begin{center}
\begin{diagram}
 & &Z \prod\limits_Y X & &\\
 &\ldTo^{\tilde{f}}& &\rdTo^{\tilde{g}}& \\
Z& & & &X\\
 &\rdTo_g&&\ldTo_f&\\
 & & Y& & &
\end{diagram}
\end{center}
%\end{minipage}

with $  \tilde{f}$, $\tilde{g}$ the natural projections. Note that
we have an identification of fibers 
\[ \tilde{f}^{-1}(z) \cong f^{-1} ( g (z) )\,  \mbox{for} \, z \in Z
\]
 and we get a map
\begin{equation} 
\label{eq1.7.2}
 A_f \rightarrow A_{\tilde{f}} ,\quad c \mapsto \tilde{c} \in A_{\tilde{f}} 
\quad \mbox{with} \quad \tilde{c}^{ (z) } = c^{ ( g (z))}. 
\end{equation}
Similarly, $   \tilde{g}^{-1} (x) \cong g^{-1} ( f (x) )$  for   
$x\in X$, and we get a map 
\[
  A_g \rightarrow A_{\tilde{g}} \quad,\quad 
  a \mapsto \tilde{a} \in A_{\tilde{g}} \quad \mbox{with}  
  \quad \tilde{a}^{ (x) } = a^{ ( f (x) ) } 
\]
The axiom of Right-Linearity states that for 
$f, g, \tilde{f}, \tilde{g}$ as above, and $  W \xleftarrow{h} Y$ in
$ Set_{\bullet}$, and for $   d \in A_{h \circ f}$, $a \in A_g$,
$ c \in A_f$,
we have in $  A_{h \circ f}$:
\begin{equation}
  ( d , c ) \circ a = (d \circ \tilde{a} , \tilde{c} )  
\end{equation}
{\bf Remark:} In the above axioms,  (1.3) through (1.7),
we have a sequence of identities, one for each fiber over a given
element $  w \in W$. We can assume without loss of generality that
$  W = [1]$, $ h$ is the cannonical map, and we get identities in
$  A_X $ (resp.  $   A_Z$, $A_Z$, $A_Y$, $A_Z$ ) for $   d \in A_Z  $
(resp.  $   A_X$, $A_Y$, $A_Z$, $A_X   $ ) in  (1.3) (resp.
(1.4) , (1.5) , (1.6) , (1.7)).

\setcounter{section}{8}
\setcounter{equation}{0}

\vspace{10pt}
\noindent{\bf (1.8) Unit axioms}
  \label{sec1.8}
  We have a distinguished element
$   1 = 1_A \in A_{[1]}$. 

Hence for each singleton $  \{ x \} \in \Fbb$,
using the unique isomorphism $   [1] \xrightarrow{ \sim} \{ x \}$ 
we get $    1_x \in A_{ \{ x \}}$; 
and for $ f \in \Fbb (X,Y)    $ 
we have $   1_f = ( 1_{ f^t(y) } )_{ y \in f(X) } \in A_f$.
 
We have for all $       f \in Set_{\bullet}(X,Y)$,
$a \in A_f$:
\begin{equation}
  \label{eq1.8.1}
\begin{array}{c}
a \circ 1_{id_X}  = a , \\
 1_{id_Y} \circ a \; = a ,  \\
 (a \,,\, 1_{id_X}) = a
\end{array}
\end{equation}
{\bf Remark:} In the above axiom (\ref{eq1.8.1}), we can assume without
loss of generality that $   Y = [1]$, $f=c_X$ the cannonical map,
and (\ref{eq1.8.1}) read
\begin{equation}
    a \circ 1_{id_X} = 1 \circ a = (a , 1_{id_X}) = a 
\ \mbox{for} \ a \in A_X  
\end{equation}

\vspace{10pt}

Note that $     (1,a)   $ makes sense only for $     a \in A_{[1]}$,
we define
\begin{equation}
      a^t = (1,a) \quad, \; a \in A_{[1]}.     
\end{equation}
It is an involution of $    A_{[1]}$:
\begin{equation}
   (a^t)^t = \left( 1, (1,a) \right) = ( 1 \circ a ,1 ) = (a,1)=a
\end{equation}
 It preserves the operation of multiplication, and the
special elements $0$, $1$:
\begin{equation}
   (a \circ b)^t = (1, a \circ b) = \Big( (1,b),a \Big) =
 \Big( 1 \circ (1,b),a \Big) =
\end{equation}
\[= (1,a) \circ (1,b) = a^t \circ b^t 
\]
\[
   0^t = (1,0)=0 \qquad , \quad 1^t = (1,1)=1  \]
 For $  f \in \Fbb(X,Y)$, and for
 \[   a = ( a^{ (y) } ) \in A_f = \prod\limits_{y \in Y} A_{ f^{-1} (y) }
 = \prod\limits_{x \in D(f)} A_{ \{ x \} } = ( A_{[1]} )^{D(f)},
\]
 we put
\begin{equation} 
  \label{eq1.8.6}
     a^t = ( a^{t(x)} ) \in A_{f^t} = \prod_{x \in X}
  A_{ \big\{ f(x) \big\} } = ( A_{[1]} )^{f(X)}   , 
\end{equation}
\[
   a^{t(x)} = ( a^{ ( f(x) ) } )^t = ( 1,a^{( f(x) )} ) 
\]
  and we have
\begin{equation}
  (1_{ f^t \circ f } , a) = a^t    
\end{equation}
\begin{equation}
\label{eq1.8.8}
\mbox{ {\bf Definition: } We shall say that } \,  A
\mbox{ is \emph{self-adjoint} if}
\ee
 \[    a^t=a  
\mbox{ for all} \, a \in A_{[1]}.\] 

\vspace{10pt}

{\bf Remark:} In $   A_{[1]}  $ we have $    0 \circ 1 = 1 \circ 0 = 0
\circ 0 = 0$, and $    1 \circ 1 = 1$, hence for $   f \in
\Fbb (X,Y)$,  $g \in \Fbb (Y,Z)$, we have 
\begin{equation}
   1_g \circ 1_f = 1_{g \circ f}  
\end{equation}
 Also we have 
(c.f. (0.10), (\ref{eq1.8.6})),
\[  \left( 1_f \right)^t = 1_{f^t}.  \]
We obtain for all $  a \in A_X$, $f \in  \Fbb (X,Y)$, 
\begin{equation}
  a \circ 1_{f^t} = (  a \circ 1_{f^t}, 1_{id_Y} ) = ( a,
(1_{id_Y},1_{f^t} ) ) = (a, 1_f). 
\end{equation}
This gives a structure of a functor $\Fbb \rightarrow Set_0$ on 
$X \mapsto A_X$, 
which we require is the given stucture (1.1):
Our unit axiom (\ref{eq1.8.1}) is completed with the requirement that
for all $a \in A_X$, $f \in \Fbb (X,Y)$,
\begin{equation}
\label{eq1.8.11}
   a \circ 1_{f^t} = (  a,1_f ) = f_A(a).    
\end{equation}  
\setcounter{section}{9}
\setcounter{equation}{0}
\noindent{\bf (1.9) Remark:} Given a commutative diagram in $   Set_{\bullet}$, 
\noindent (1.9.1)
\setcounter{equation}{1}
\begin{center} 
 \begin{diagram}
X&\rTo^{f}&Y\\
\dTo^{\varphi}& &\dTo_{\psi}\\
X' & \rTo^{f'}&Y'
\end{diagram}
\end{center}
with 
 $    \varphi \in \Fbb (X,X') , 
       \psi \in \Fbb (Y,Y')$,
we get a map
\begin{equation}
          ( \varphi , \psi )_A: A_f \longrightarrow
 A_{f'} 
\quad,
\end{equation}
\[   ( \varphi , \psi )_A(b) = 1_{\psi} \, \circ \, b \, \circ \, 1_{\varphi^t} \;\, for \;\, b \in A_f.  
\]
   The operations of multiplication and contraction are functorial
   in the following sense. 
   For $  a \in A_{Y'}$, we have in $  A_{X'}$,
\begin{equation}
    a \circ (\varphi,\psi)_A (b) = a \circ 1_\psi \circ b \circ
1_{\varphi^t} = \varphi_A \left(  \psi_A^t(a) \circ b  \right) 
\end{equation}
For $  a \in A_{X'}$, we have in $    A_{Y'}$, 
\begin{equation}
   \Big(  a, (\varphi,\psi)_A(b) \Big) = \Big(  a, 1_\psi \circ b
\circ 1_{\varphi^t} \Big) = \Big(  a \circ 1_\varphi , 1_\psi \circ
b \Big) = 
\end{equation}
\[=\Big(   (a \circ 1_\varphi , b), 1_\psi  \Big) = \psi_A
\Big(   (\varphi_A^t(a) ,b)   \Big)
  \]
  
\setcounter{section}{10}
\setcounter{equation}{0}
{\bf (1.10) Remark:} If follows from (1.3), (1.8),
  that $    A_{[1]} = A_{id_{[1]}}   $  is an associative monoid
  with unit 1. It also follows that it is commutative: 
\begin{equation}
     a \circ b = (a \circ b,1) = (1 \circ a, (1,b)) =
(1,(1,b))\circ a = b \circ a   
\end{equation}
For $  X \in \Fbb$, the commutative monoid $  A_{id_X} =
(A_{[1]})^X $ acts on the right on $  A_X$, cf. (\ref{eq1.2.7}).
Note that for $  a=(a^{(x)}) \in A_{id_X}$, and $  b \in A_X$,
\begin{equation}
  b \circ a = (b \circ a,1_{id_X})  = (b,(1_{id_X},a)) = (b,a^t) 
\end{equation}
i.e. the two actions of (\ref{eq1.2.7}) are related by the
involotion.
 The monoid $  A_{[1]} $ also acts on the left on $ A_X$, 
cf. (\ref{eq1.2.5}),  and this is nothing but the diagonal
right action: For $  a \in A_{[1]}$, $b \in A_X$, putting 
$\tilde{a} \in A_{id_X}$, $\tilde{a}^{(x)} = a $ for all $  x \in X$,
\begin{equation}
\label{eq1.10.3}
  a \circ b = (1, a^t) \circ b = (1 \circ b, \tilde{a}^t) =
(b,\tilde{a}^t ) = b \circ \tilde{a}. 
\end{equation}
More generally, for $  f \in Set_{\bullet}(X,Y)$, $b \in A_f$, 
$a= (a^{(y)})
\in A_{id_Y} = (A_{[1]})^Y$, we let
\begin{equation}
  \tilde{a} \in A_{id_X} = (A_{[1]})^X \;,\;
\tilde{a}^{(x)} =  a^{( f(x) )} \, \mbox{for} \, x \in X,  
\end{equation}
and we have in $  A_f$, 
\begin{equation}
\label{eq1.10.5}
  a \circ b = b \circ \tilde{a}.   
\end{equation}
Indeed checking (\ref{eq1.10.5}) at a given fiber
$   A_{f^{ -1 } (y) }$, $y \in Y$,
reduces to (\ref{eq1.10.3}). 
The action of $   (A_{[1]})^X    $ on $  A_X  $ is self-adjoint
with respect to the pairing (1.2.6): for
$  a \in (A_{[1]})^X$,  $b,d \in A_X$, we have 
\begin{equation}
   ( b \circ a, d ) = ( b, (d,a) ) = ( b, d \circ a^t )  
\end{equation}

\setcounter{section}{11}
\setcounter{equation}{0}
{\bf (1.11) Remark:} There is enough commutativity in our axioms to 
make the theory work.
However, all our examples (and calculations) will satisfy the extra axiom of\\
\underline{\bf Commutativity}: For
$    Z \xrightarrow{g} Y \xleftarrow{f} X      $ in $   Set_{\bullet}$,
and for $   a \in A_g$, $c \in A_f$, we have in
$    A_{ g \circ \tilde{f} } = A_{ f \circ \tilde{g}}$, 
\begin{equation}
\label{eq1.11.1}
  a \circ \tilde{c} = c \circ \tilde{a}
\end{equation}
using the notations of (\ref{eq1.7.1}) , (\ref{eq1.7.2}).

\vspace{10pt}

Eventually, all our generalized rings will be self-adjoint and commutative,
but we choose to emphasize the role of self-adjunction and of commutativity in the development of the theory.
The reader may assume for simplicity that all generalized rings are self-adjoint and commuative throughout the paper.

\vspace{10pt}

\setcounter{section}{12}
\setcounter{equation}{0}
{\bf (1.12) Remark:} For $  f,g \in  Set_{\bullet}(X,Y)$ we write 
$   f \leq g  $ for
\begin{equation}
 D(f) \subseteq D(g)   \, \mbox{and} \,    g|_{D(f) } = f    
\end{equation}
or equivalenty, if for all $  y \in Y$, $f^{-1}(y) \subseteq  g^{-1}(y)  $.
For $   f \in Set_{\bullet}(X,Y)$,  $h \in Set_{\bullet}(X,Z)$, 
$g \in Set_{\bullet}(Y,Z)$
such that $   h \leq g \circ f$, we have\\
$  ( A_h,A_f )  \subseteq  A_g  $ via
\begin{equation}
  A_h  \times  A_f  \hookrightarrow  A_{g \circ f}  \times  A_f  \xrightarrow{ (,) }  A_g     
\end{equation}
Given such $ f $ and $ h$, there exists $ g $ such that
$    h \leq g \circ f$, if and only if  
\begin{equation}
\label{eq1.12.3}
   D(h) \subseteq D(f) \, \mbox{ and for } \,
   x_1,x_2 \in D(h) : 
\end{equation}
\[
h(x_1) \neq h(x_2) \Rightarrow  f(x_1) \neq f(x_2)  
\]
or equivalenty,
\begin{equation}
\label{eq1.12.4}
   D(h) \subseteq D(f) , \, \mbox{ and for } \, 
   z_1 \neq  z_2 , z_i \in Z :
\end{equation}
\[
f( h^{-1} (z_1) ) \cap  f( h^{-1} (z_2) )  =  \emptyset
\] 
In this case we have $  h/f \in Set_{\bullet}(Y,Z)   $  with
$  D(h/f) = f( D(h))$, and
\begin{equation}
  h/f( f(x) ) = h(x)  
\end{equation}
This function $  h/f  $ is the minimal solution:
for all $   g \in Set_{\bullet}(Y,Z)  $ such that $   h \leq g \circ f$,
we have $   h/f   \leq    g$. 
Thus when (\ref{eq1.12.3}) or (\ref{eq1.12.4})
are satisfied we have a well-defined contraction
\begin{equation}
\label{eq1.12.6}
  ( A_h,A_f ) \subseteq  A_{h/f} 
\end{equation}
Similarly, given $ f $ and $ h$, there is a $ g $ such that
$  h \geq g \circ f  $ and $  D(g) = f(X)$, if and only if 
\begin{equation}
   D(h) \supseteq D(f), \, \mbox{   and for } \,
   z_1 \neq  z_2 , z_i \in Z :
\end{equation}
\[
 f( h^{-1} (z_1) ) \cap  f( h^{-1} (z_2) )  = \emptyset
\]
In this case $ g = f \setminus h  $ is unique, $  D(f \setminus h) =
f(X)$, $f \setminus h( f(x) ) = h(x)   $ and we have a well-defined
contraction
\begin{equation}
\label{eq1.12.8}
A_h \times A_f \twoheadrightarrow A_{(f \setminus h) \circ f} \times A_f
 \xrightarrow{( \, ,\, )} A_{f \setminus h} , \quad
 (A_h , A_f ) \subseteq A_{h / f}
\end{equation}
Putting (\ref{eq1.12.6}) and (\ref{eq1.12.8}) together, 
we have for $f \in Set_{\bullet} (X, Y)$, $h \in Set_{\bullet} (X, Z)$,
such that
\begin{equation}
\mbox{for } \, z_1 \neq z_2 , z_i \in Z , \quad
f(h^{-1} (z_1) ) \cap f(h^{-1} (z_2) ) = \emptyset 
\end{equation}
a well defined contraction
\begin{equation}
(A_h , A_f ) \subseteq A_{h/f} , \, \mbox{with } \,
D(h/f) = f(D(h)) , 
\end{equation}
\[ \mbox{and } \,
h/f (f(x)) = h(x)
\]

\vspace{10pt}

Note that in the axioms of right-adjunction (\ref{eq1.5.1}),
and left-linearity (\ref{eq1.6.1}), 
it is possible for the left hand side of the equality to be defined,
while the right hand side is not defined ( the contraction $(a,c)$ is not defined, but $(d \circ c , a)$  is defined).
Thus some care is needed in passing from the left to the right of 
(\ref{eq1.5.1}), (\ref{eq1.6.1}).
If the right side of (\ref{eq1.5.1}) or  (\ref{eq1.6.1})
is defined, so is the left side, and we have an equality.

In particular, the axioms allow to transform any formula in the operations of multiplication and contraction, to an equivalent formula with only
\underline{one} contraction.
Thus expressions of the form
\[
 (a,b) = (a_1 \circ a_2 \circ \cdots \circ a_n , b_1 \circ \cdots \circ b_m)\]
are closed under multiplication and contraction, and we have the formulas:
\begin{equation}
\label{eq1.12.11}
(a,b) \circ (c,d) = (a \circ \tilde{c} , d \circ \tilde{b} )
\end{equation}

\begin{equation}
\label{eq1.12.12}
\left( (a,b) , (c,d) \right) = (a \circ \tilde{d} , 
c \circ \tilde{b}).
\end{equation}

\setcounter{section}{13}
\setcounter{equation}{0}
{\bf (1.13) Definition:} A homomorphism
 of generalized rings
$\varphi : A \rightarrow A'$
is a natural transformation of functors 
(so for $X \in \Fbb$,
we have $\varphi_X \in Set_0 (A_X , A'_X)$,
and for $f \in \Fbb (X,Y)$,
we have $\varphi_Y \circ f_A = f_{A'} \circ \varphi_X$),
such that $\varphi$ preserves multiplication, contraction, and the unit:
\begin{equation}
\label{eq1.13.1}
\varphi(a \circ b) = \varphi (a) \circ \varphi (b) ,
\end{equation}
\begin{equation}
\label{eq1.13.2}
\varphi((a , b)) = (\varphi (a),  \varphi (b)) ,
\end{equation}
\begin{equation}
\label{eq1.13.3}
\varphi (1_A) = 1_{A'}
\end{equation}
Thus we have a category of generalised rings and homomorphisms 
which we denote by $\cG \cR$.
We let $\cG\cR_{\cC} \subseteq \cG\cR$  
denote the full subcategory consisiting of commutative generalized rings,
 i.e. those satisfying the extra commutativity axiom (\ref{eq1.11.1}).
We let $\cG\cR^+$ denote the full subcategory consisiting of 
self-adjoint generalized rings, 
i.e. those satisfying (\ref{eq1.8.8}).
We denote by ${\cal GR_C^+ = GR_C \cap GR^+} $
the commutative self-adjoint generalized rings.

We remark that for $A, A' \in \cG\cR$,
a collection of maps 
$\varphi = \{\varphi_X \in Set_0 (A_X , A'_X )\}$
satisfying 
(\ref{eq1.13.1}), (\ref{eq1.13.2}), (\ref{eq1.13.3}), 
is a homomorphism; 
i.e. it is automatically a natural transformation of functors by (\ref{eq1.8.11}).

\chapter{Examples of Generalised Rings}
\setcounter{equation}{0}
\section{The  field with one element $\mathbb{F}$ }

We write $\mathbb{F}$ for the functor 
$(0.8):\mathbb{F}_{\bullet} \stackrel{\thicksim}{\rightarrow} \mathbb{F}_0 \subseteq Set_0$.
Thus
\begin{equation}
\mathbb{F}_X=X_0=X \coprod  \{0_X \} , \, \mbox{and} \,
\mathbb{F}_f = \prod_y\left(f^{-1}(y) \right)_0 \, \mbox{for} \,
f \in Set_{\bullet}(X,Y).
\end{equation}

The multiplication is given by
\begin{equation}
\begin{array}{ll}
\mathbb{F}_Y \times \mathbb{F}_f & \rightarrow  \mathbb{F}_X\\
y_0 , (x_f^{(y)}) & \mapsto y_0 \circ x_f = x_f^{(y_0)} \in 
          \left(f^{-1}(y_0)\right)_0\subseteq X_0 
\end{array}
\end{equation}
The contraction is given by
\[
\begin{array}{ll}
\mathbb{F}_X \times \mathbb{F}_f & \rightarrow  \mathbb{F}_Y\\
x_0 , (x_f^{(y)}) & \mapsto (x_0 , x_f) = 
   \left\{\begin{array}{ll}
y & \, \mbox{if} \, x_0 = x_f^{(y)}\\
0 & \mbox{otherwise}.
\end{array}\right.
\end{array}
\]
It is easy to check that $\mathbb{F}$ is a (commutative, self-adjoint) generalised ring.
For $A \in {\cal GR}$, and for $x \in X \in \mathbb{F}_{\bullet}$,
put 
\begin{equation}
\varphi_X(x) = (1_x , 1_{j_x}) = 1_x \circ 1_{j_x^t} \in A_X , \,
\mbox{with} \, j_x \in \mathbb{F}_{\bullet}(\{x\}, X) \,
\mbox{the inclusion.}
\end{equation}
We have $\varphi_X \in Set_0 (\mathbb{F}_X , A_X)$, 
and the collection of $\varphi_X$ define a homomorphism 
$\varphi \in {\cal GR}(\mathbb{F}, A)$.
It is easy to check this is the only possible homomorphism,
and $\mathbb{F}$ is the initial object of ${\cal GR}$.

\section{Commutative Rings}
\label{sec2.2}
%reset equation counter
\setcounter{equation}{0}

For a commutative ring $A$, let ${\cal G}(A)_X = A \cdot X = A^X$
be the free $A$-module with basis $X$.
It forms a functor ${\cal G}(A): \mathbb{F}_{\bullet} \rightarrow A \cdot mod \subseteq Set_0$.
We define the multiplication for $f \in Set_{\bullet}(X, Y)$
by 
\begin{equation}
{\cal G}(A)_Y \times {\cal G}(A)_f = A^Y \times \prod_{y \in Y} A^{f^{-1}(y)}
        \rightarrow A^X = {\cal G}(A)_X 
\end{equation}
\[ a = (a_y) , b = (b_x^{(y)})_{x \in f^{-1}(y)} \mapsto (a \circ b )_x 
       = a_{f(x)} \cdot b_x^{f(x)}
\]
We define the contraction by
\begin{equation}
{\cal G}(A)_X \times {\cal G}(A)_f = A^X \times \prod_{y \in Y} A^{f^{-1}(y)}
        \rightarrow A^Y = {\cal G}(A)_Y 
\end{equation}
\[
a = (a_x) , b = (b_x^{(y)})_{x \in f^{-1}(y)} \mapsto (a , b )_y =
        \sum_{x \in f^{-1}(y)} a_x \cdot b_x^{(y)}
\]
It is straightforward to check that ${\cal G}(A)$ is a 
commutative self-adjoint generalized ring.
A homomorphism of (commutative) rings $\varphi \in Ring (A,B)$ 
gives a homomorphism
${\cal G}(\varphi) \in {\cal GR}({\cal G}(A), {\cal G}(B))$, 
thus we have a functor
\begin{equation}
\label{eq2.2.3}
{\cal G} : Ring \rightarrow {\cal GR}_C^+
\end{equation}
It is fully-faithful: if $\varphi \in {\cal GR}({\cal G}(A), {\cal G}(B))$, 
and $a = (a_x)\in {\cal G}(A)_X$,
then $\varphi_X(a)_x =\left(\varphi_{[1]}(a_x)\right)$ by functoriality 
over $\mathbb{F}_{\bullet}$, so $\varphi$ is determined by 
$\varphi_{[1]}: A \rightarrow B$;
the map $\varphi_{[1]}$ is multiplicative (and preserves $1$),
but it is also additive
\[
\varphi_{[1]} (a_1 + a_2) = \varphi_{[1]} \left(\left( (a_1, a_2) , (1,1)\right)\right) =
((\varphi_{[1]}(a_1) , \varphi_{[1]}(a_2)) , (1,1)) =\]
\[=\varphi_{[1]}(a_1) + \varphi_{[1]}(a_2)
\]
Thus $\varphi_{[1]} \in Ring(A, B)$ , and 
$\varphi = {\cal G }(\varphi_{[1]})$; and we have
\begin{equation}
Ring(A,B) = {\cal GR}({\cal G}(A) ,{\cal G}(B) )
\end{equation}

Note that for every $X \in \mathbb{F}_{\bullet}$, 
we have a distinguished element
\begin{equation}
\mathbbm{1}_X  \in {\cal G}(A)_X , \quad
\left(\mathbbm{1}_X\right)_x=1 \, \mbox{for all} \, x \in X
\end{equation}
hence for $f \in Set_{\bullet}(X,Y)$, the element 
$\mathbbm{1}_f = (\mathbbm{1}_{f^{-1}(y)}) \in {\cal G}(A)_f$.
These elements satisfy
\begin{equation}
\label{eq2.2.6}
\mathbbm{1}_g \circ \mathbbm{1}_f = \mathbbm{1}_{g \circ f} ; \quad
\mathbbm{1}_{[1]} =1 ;
\end{equation}
\[
\mathbbm{1}_X \circ 1_{f^t} = \mathbbm{1}_{f(X)} \in {\cal G}(A)_{f(X)}
\subseteq  {\cal G}(A)_Y \, \mbox{for} \, f \in Set_0(X,Y)
\]
Note that any element 
$a  = (a_x) \in  {\cal G}(A)_X$, gives an element of the monoid
\begin{equation}
\langle a \rangle = (a_x) \in  {\cal G}(A)_{id_X}= A^X
\end{equation}
and the vector $\mathbbm{1}_X$ is cyclic in the sense that
\begin{equation}
\label{eq2.2.8}
a = \mathbbm{1}_X \circ \langle a \rangle
\end{equation}

The definition of ${\cal G}(A)$ does not use subtraction and 
thus works for semi-rings $A$ (having two associative and commutative 
operations, addition $x+y$ , and multiplication $x \cdot y$,
with units $0$ and $1$, and multiplication  distributive over addition).
We have the generalized rings ${\cal G} (\mathbb{N})$ and 
${\cal G}([0, \infty ))$,
as well as the ``tropical'' examples ${\cal G} (\mathbb{N}_t)$ and 
${\cal G}([0, \infty )_t)$
where we replace addition by the operation of taking the maximum 
$\max \{x,y\}$.

\vspace{10pt}

The generalized ring ${\cal G} (\mathbb{N})$ contains the cyclic vectors 
$\mathbbm{1}_X \in {\cal G} (\mathbb{N})_X = \mathbb{N}^X$,
satisfying (\ref{eq2.2.6}),
 and these vectors generate ${\cal G} (\mathbb{N})$.
For $n = (n_x) \in {\cal G} (\mathbb{N})_X$,
we have the set over $X$.
\begin{equation}
\pi_n : X_n \rightarrow X
\end{equation}
\[X_n = \left\{ (x,j) \in X \times \mathbb{N}_+ , j \leq n_x \right\} \, , \, \pi_n(x,j) = x
\]
and we have
\begin{equation}
n=\left( \mathbbm{1}_{X_n} ,\mathbbm{1}_{\pi_n} \right).
\end{equation}
We can view the elements of ${\cal G}(\mathbb{N})_X$ as 
isomorphism classes of sets over $X$.

Let $A\in {\cal GR}$, and let $a =(a_X)$,
with $a_X \in A_X$, $a_{\{x\}}= 1_x$,
and with $a_g \circ a_f = a_{g \circ f}$ (where $a_f = (a_{f^{-1}(y)})$ 
for $f \in Set_{\bullet}(X,Y)$).
Define a map $\varphi = \varphi^{(a)}:{\cal G}(\mathbb{N}) \rightarrow A $ by
\begin{equation}
\varphi_X : \mathbb{N}^X \rightarrow A_X  , \quad \varphi_X(n) = (a_{X_n}, a_{\pi_n}).
\end{equation}
For $f \in Set_{\bullet}(X,Y) $, and for $m = (m^{(y)}) \in 
{\cal G}(\mathbb{N})_f$, we have similarly 
$m^{(y)}= \left( \mathbbm{1}_{X_m^{(y)}} ,\mathbbm{1}_{\pi_m^{(y)}} \right)$
with $\pi_m^{(y)}: X_m^{(y)} \rightarrow f^{-1}(y)$ the restriction of 
$\pi_m : X_m \rightarrow D(f) \subseteq X$.
We have by formula (\ref{eq1.12.12}),
\begin{equation}
\label{eq2.2.12}
\left((a_{X_n}, a_{\pi_n}), (a_{X_m}^{(y)} , a_{\pi_m}^{(y)}) \right) =
\left( a_{X_n} \circ \widetilde{a}_{\pi_m}, a_{X_m} \circ \widetilde{a}_{\pi_n}\right)=
\end{equation}
\[=
\left( a_{X_n \prod\limits_X X_m} , a_{\pi_Y}\right) \in A_Y
\]

%\begin{figure}[h]
%\begin{center}
%\psfrag{v1}{$[1]$}
%\psfrag{v2}{$Y$}
%\psfrag{v3}{$X_m$}
%\psfrag{v4}{$X_n \Pi_{X} X_m$}
%\psfrag{v5}{$X_n$}
%\psfrag{v6}{$X$}
%\includegraphics[scale=1]{penta.eps}
%\caption{\,}
%\end{center}
%\end{figure}
with $\pi_Y:X_n \prod\limits_X X_m \rightarrow X_m \rightarrow Y$.
Note that for $y \in Y$ we have,
%\end{minipage}

\begin{center}
\begin{diagram}
 & &X_n \prod\limits_X X_m & & \\
 &\ldTo^{\widetilde{a}_{\pi_m}}& &\rdTo^{\widetilde{a}_{\pi_n}} &\\
X_n& & & & X_m\\
\dTo^{a_{X_n}}&\rdTo^{a_{\pi_n}}&&\ldTo^{a_{\pi_m}}&\dTo_{a_{X_m}}\\
 & &X & & \\
 & \ldTo& & \rdTo_f &\\
[1]&&\lTo& & Y
\end{diagram}
\end{center}

\begin{equation}
\sharp \pi^{-1}_Y (y) = \sum_{x \in f^{-1}(y)}n_x \cdot m_x = (n,m)_y
\end{equation}

Similarly, for $n = (n_y) \in {\cal G}(\mathbb{N})_Y = \mathbb{N}^Y $,
we have by formula (1.12.11),
\begin{equation}
\label{eq2.2.14}
(a_{Y_n}, a_{\pi_n}) \circ (a_{X_m^{(y)}} , a_{\pi_m^{(y)}}) =
(a_{Y_n} \circ \widetilde{a}_{X_m}, a_{\pi_m} \circ \widetilde{a}_{\pi_n}) =
(a_{Y_n \prod_Y X_m} , a_{\pi_X}) \in A_X
\end{equation}

with $\pi_X: Y_n \prod\limits_Y X_m \rightarrow X_m \rightarrow X$.
Note that for $x \in X$, we have

\begin{center}
\begin{diagram}
 & & Y_n \prod\limits_Y X_m & &\\
 &\ldTo^{\widetilde{a}_{X_m}}& &\rdTo^{\widetilde{a}_{\pi_n}}& \\
Y_n& & & & X_m\\
\dTo^{a_{Y_n}}&\rdTo^{a_{\pi_n}}& & \ldTo^{a_{X_m}}&\dTo_{a_{\pi_m}}\\
 & & Y &\lTo_f & X\\
&\ldTo & & &\\
[1]& & & & 
\end{diagram}
\end{center}

\begin{equation}
\label{eq2.2.15}
\sharp \pi^{-1}_{X}(x) = n_{f(x)}\cdot m_x = (n \circ m)_{(x)}
\end{equation}

Thus $\varphi$ is a homomorphism of generalized rings, and we have a natural bijection
\begin{equation}
\label{eq2.2.16}
{\cal GR}({\cal G}(\mathbb{N}), A) =
\left\{a = (a_X) \in \prod_{X \in \mathbb{F}}A_X , a_{\{x\}} = 1_x , a_g \circ a_f = a_{g \circ f}
\right\}
\end{equation}

\section{The Real Prime}
\label{sec2.3}
\setcounter{equation}{0}

Let $|\, | :K \rightarrow [0, \infty)$ be a non-archimedian absolute value on a field $K$.
Let \begin{equation}
{\cal O}_X = \left\{ a = (a_x) \in {\cal G}(K)_X \, ,\, 
\sum\limits_{x \in X}|a_x|^2 \leq 1 \right\}
\end{equation}
Note that $X \mapsto {\cal O}_X$ is a subfunctor of 
${\cal G}(K): \mathbb{F}_{\bullet} \rightarrow Set_0$.
But it is also a 
\emph{sub-generalized-ring}, in the sense that it is closed under 
the operations of multiplication and contraction, and $1 \in {\cal O}_{[1]}$.
The proof for multiplication is straightforward:

For $a =(a_y) \in {\cal O}_Y$, 
    $b =(b_x^{(y)}) \in {\cal O}_f$, we have 
\begin{equation}
\sum\limits_{x \in X} |(a \circ b)_x|^2 =
\sum\limits_{x \in X} |a_{f(x)}|^2 \cdot |b_x^{(f(x))}|^2 =
\end{equation}
\[
=\sum\limits_{y \in Y} |a_y|^2 \sum\limits_{x \in f^{-1}(y)}|b_x^{(y)}|^2 \leq
\sum\limits_{y \in Y}  |a_y|^2 \leq 1
\]
and so $a \circ b \in {\cal O}_X$. 

The proof for contraction is just the Cauchy-Schwartz inequality:
For $a =(a_x) \in {\cal O}_X$, 
    $b =(b_x^{(y)}) \in {\cal O}_f$ , we have
\begin{equation}\begin{array}{ll}
\sum\limits_{y \in Y} |(a , b)_y|^2 &=
\sum\limits_{y \in Y}\left|\sum\limits_{x \in f^{-1}(y)} a_x \cdot b_x^{(y)}
                     \right|^2 \leq\\

 &\leq \sum\limits_{y \in Y} \left(\sum\limits_{x \in f^{-1}(y)}|a_x|^2\right) \cdot
                     \left(\sum\limits_{x \in f^{-1}(y)}|b_x^{(y)}|^2\right)\\
&\leq \sum\limits_{y \in Y}\sum\limits_{x \in f^{-1}(y)}|a_x|^2 \leq
 \sum\limits_{x \in X}|a_x|^2 \leq 1
\end{array}\end{equation}
and so $(a,b) \in {\cal O}_Y$.

Thus ${\cal O}$ is a generalized ring. Note that
\begin{equation}
m_X = \left\{ a = (a_x) \in {\cal O}_X , \sum\limits_{x \in X}|a_x|^2  <1 \right\}
\end{equation}
forms a subfunctor of ${\cal O}$, and it is (the unique maximal)
\emph{ideal} of ${\cal O}$,
in the sense that we have
\begin{equation}
\label{eq2.3.5}
{\cal O} \circ m \,  , \,  m \circ {\cal O} \, , \,   (m,{\cal O})\, ,\, ({\cal O},m) \subseteq m
\end{equation}
By collapsing $m_X$ to zero we obtain the quotient (in $Set_0$):
\begin{equation}
k_X = {\cal O}_X / m_X =
\left\{  a = (a_x) \in {\cal O}_X , \sum\limits_{x \in X}|a_x|^2  =1\right\} \amalg 
\{ 0_X \}
\end{equation}
There is a canonical projection map $\pi_X:{\cal O}_X -\hspace{-5pt}-\hspace{-10pt}\gg k_X$,
with $\pi_X(m_X) = 0_X$.
Now (\ref{eq2.3.5}) imply that there is a (unique) structure of 
 a generalized ring on $k$,
such that $\pi$ is a homomorphism,
$\pi \in {\cal GR}({\cal O}, k)$. It is given by
\begin{equation}
a \circ b = \left\{\begin{array}{cl}a \circ b & \, \mbox{if} \ 
                                             ||a \circ b||=1\\
0 & \, \mbox{if}\ ||a \circ b||< 1
 \end{array}\right.
\quad,
(a,b)=\left\{\begin{array}{cl} (a,  b) & \, \mbox{if} \ 
                                             ||(a, b)||=1\\
0 & \, \mbox{if}\ ||(a, b)||< 1
 \end{array}\right.
\end{equation}
with the $l_2$ -norm
\begin{equation}
||(a_x) || = \left( \sum\limits_{x \in X} |a_x|^2 \right)^{\frac{1}{2}}
\end{equation}

\section{Monoids}
\label{sec2.4}
\setcounter{equation}{0}

We say that $M$ is a 
\emph{monoid} if it has an associative and commutative operation $\cdot$, 
an involution, and there are (neccessarily unique) elelments $0, 1 \in M$
such that 
\begin{equation}
m \cdot 1 = m \, , \, m \cdot 0 = 0 \, , \, \mbox{for all} \, m \in M
\end{equation}
the involution preserves the operation and special elelments
\[(m_1 \cdot m_2)^t = m_1^t \cdot m_2 ^t \, ,\, 1^t = 1 \, ,\, 0^t = 0 \,,
\, \mbox{and} \, (m^t)^t = m
\]
A map of monoids $\varphi : M \rightarrow M'$ is required to 
preserve the operation 
$\varphi (m_1 \cdot m_2) = \varphi(m_1) \cdot \varphi(m_2)$,
as well as the involution and special elements:
$\varphi (m^t) = \varphi (m)^t \, , \, \varphi (1) = 1 \, , \, \varphi (0 ) =0$.
Thus we have the category of monoids which we denote by $Mon.$
For $A \in {\cal GR}$, we have 
$A_{[1]} \in Mon$,
cf. (1.10), giving rise to a functor
\begin{equation}
{\cal GR} \rightarrow Mon \, ,\, A \mapsto A_{[1]}
\end{equation}
This functor has a left adjoint. For $M \in Mon$, let
\begin{equation}
\mathbb{F} [M]_X = \coprod\limits_{x \in X} M = 
\left( X \times (M \setminus \{ 0 \})\right) \amalg \{ 0 \} =
\left[X \times M \right] / [x,0] \sim 0
\end{equation}
It forms a functor $X \mapsto \mathbb{F} [M]_X : \mathbb{F}_{\bullet}
\rightarrow Set_0$.
We define the operation of multiplication by
\begin{equation}
\circ:\mathbb{F}[M]_Y \times \mathbb{F}[M]_f \rightarrow \mathbb{F}[M]_X
\end{equation}
\[[y_0 , m_0 ] \, ,\, [x^{(y)} , m^{(y)} ] \mapsto
[x^{(y_0)} , m_0 \cdot m^{(y_0)}]
\]
We define the operation of contraction by
\begin{equation}
( \,\, , \,\, ) :\mathbb{F}[M]_X \times \mathbb{F}[M]_f \rightarrow \mathbb{F}[M]_Y
\end{equation}
\[ [x_0 , m_0 ] \, ,\, [x^{(y)} , m^{(y)} ] \mapsto
[ y , m_0 ( m^{(y)})^t] \, \mbox{if } \, x_0 =x^{(y)}, \,
\mbox{otherwise} \, 0.
\]
It is straightforward to check that with these operations 
$\mathbb{F}[M]$ forms a (commutative, but non-self adjoint) generalized ring,
and for $A \in {\cal GR}$, we have adjunction
\begin{equation}
{\cal GR} ( \mathbb{F}[M] , A ) = Mon(M, A_{[1]})
\end{equation}
\[\varphi \mapsto \varphi_{[1]}
\]
\[\widetilde{\psi}_X ([x, m] ) = ( \psi(m) , 1_{j_x}) \leftmapsto \psi
\]
\[\mbox{with} \, j_x \in \mathbb{F}([1], X ) , j_x (1) =x
\]
We let $Mon^+ \subseteq Mon$
denote the full-subcategory of monoids with trivial involution,
we have similarly adjoint functors
\begin{equation}
\begin{array}{ccc}
Mon^+ & \stackrel{\textstyle \curvearrowright}{\curvearrowbotleft}& {\cal GR}^+\\
M & \mapsto & \mathbb{F}[M]\\
A_{[1]} & \leftmapsto & A
\end{array}
\end{equation}

\begin{equation}
\label{eq2.4.8}
{\cal GR}^+ (\mathbb{F}[M] , A) = Mon^+(M , A_{[1]}) 
\end{equation}

\section{Commutative trees and $\Delta$}
\setcounter{equation}{0}

By a \emph{tree} $F$
we shall mean a finite set with a distinguished element $0_F \in F$,
the \emph{root},
and a map ${\cal S} = {\cal S}_F :F \setminus  \{ 0_F \} \rightarrow F $,
such that for all $a \in F$ there exists $n$ with 
$\cS^n (a) = 0_F$;
we write $n = ht(a)$, and put $ht(0_F) =0$.
For $a \in F$, we put $\nu (a) = \sharp \cS^{-1}(a)$. 
The boundary of $F$ is the set
\begin{equation}
\partial F = \{ a \in F \, ,\, \nu (a) =0 \}
\end{equation}
The unit tree is the tree with just a root,
$\{ 0 \}$, and
$\partial \{ 0 \} = \{ 0 \} $.
The zero tree is the empty set, $\emptyset$, and $\partial \emptyset = \emptyset$.
Given a subset $B \subseteq \partial F$, we have the reduced tree 
$F|_B$, 
obtained from $F$ by omitting all the elements of $\partial F \setminus B$,
and then omitting all elements $a \in F $ such that all the elements of 
$\cS^{-1}(a)$ have been omitted and so on;
we have $\partial (F|_B)=B$.

\newpage
\begin{equation} 
\label{eq2.5.2}
\end{equation}
\begin{center}
\includegraphics[scale=1]{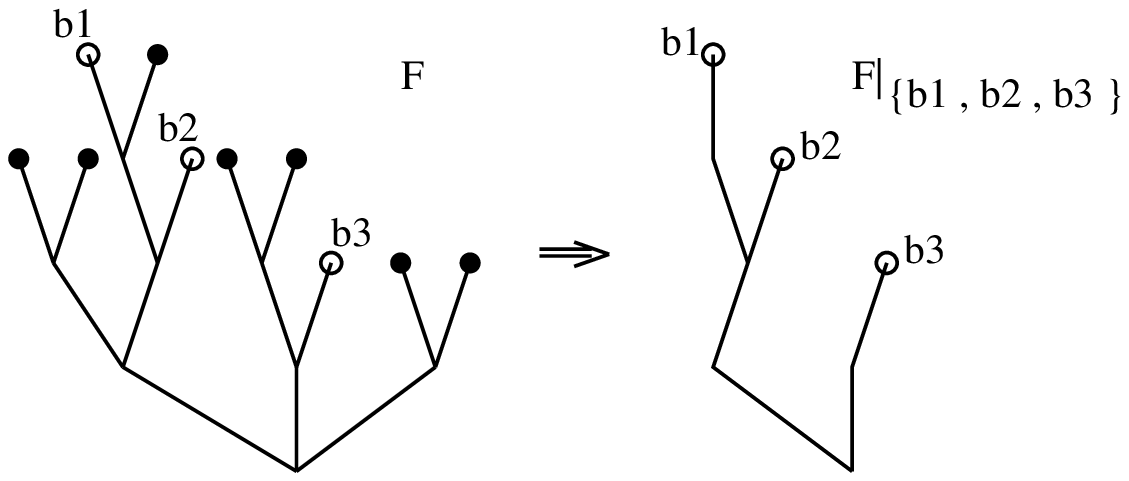}
\end{center}
Given for each $b \in \partial F $, a tree $G_b$,
let $B =\{b \in \partial F \, , \, G_b \neq \emptyset \}$, and form the tree
\begin{equation}
\label{eq2.5.3}
F \ltimes  G := F|_B \amalg \coprod\limits_{b \in B} \left( G_b \setminus \{ 0_{G_b}\} \right)
\end{equation}
with $\cS (a) = b $ if $a \in G_b$ and $\cS_{G_b}(a) = 0_{G_b}$, 
and otherwise $\cS$ being the restriction of the given $\cS_F$ and $\cS_{G_b}$.

\vspace{10pt}

Let for $X \in \mathbb{F}_{\bullet}$, 
\[\Delta_X = \{ F = (F_1 ; \{ \bar{F}_x \} ; \sigma )\}/\hspace{-3pt}\approx \]
consist of the data of a tree $F_1$,
a tree $\bar{F}_x$ for each $x \in X$, and a bijection 
$\sigma : \partial F_1 \stackrel{\thicksim}{\rightarrow} \coprod\limits_{x \in X} \partial \bar{F}_x$.
We view such data $F$ as being the same as 
$F' = (F_1' ; \{\bar{F}_x'\} ; \sigma'   )$
if there are isomorphisms of trees 
$\tau_1 : F_1 \stackrel{\thicksim}{\rightarrow} F'_1 $,
$\tau_x : \bar{F}_x \stackrel{\thicksim}{\rightarrow} \bar{F}'_x $,
with 
$\sigma' \circ \tau_1 (b) = \tau_x \circ \sigma (b)$
for $b \in \partial F_1$, $\sigma (b) \in \partial \bar{F}_x $.

Note that for such data $F=(F_1 ; \{\bar{F}_x\} ; \sigma)$,
we have an associated map
\begin{equation}
\underline{\sigma} : \partial F_1 \rightarrow X
\end{equation}
\[\underline{\sigma}(b) = x \, \mbox{iff} \, \sigma(b) \in \partial \bar{F}_x 
\]
For $f \in Set_{\bullet} (X,Y)$ , we have 
$\Delta_f = \prod\limits_{y \in Y} \Delta_{f^{-1}(y)}$,
its elements are isomorphism classes of data 
$F=\left( \{ F_y \}_{y \in f(X)} ; \{ \bar{F}_x \}_{x \in D(f)} ; \sigma_y \right)$,
with bijections 
$\sigma_y : \partial F_y \stackrel{\sim}{\rightarrow} \coprod\limits_{x \in f^{-1}(y)}
  \partial \bar{F}_{x}$,
for all $y \in Y$.

We define the operations of multiplication by
\begin{equation}
\label{eq2.5.7}
\circ : \Delta_Y \times \Delta_f \rightarrow \Delta_X
\end{equation}
\[G \circ F = (G_1 ; \bar{G}_y ; \tau) \circ (F_y ; \bar{F}_x ; \sigma_y) = 
   (G_1 \ltimes F_{\underline{\tau}(\, \, )} ;
\bar{F}_x \ltimes \bar{G}_{f(x)} ;
 \tau \circ \sigma)
\]
with
\begin{equation}
 \tau \circ \sigma : \partial (G_1 \ltimes F_{\underline{\tau}(\, \,)}) =
\hspace{2.5in}
\end{equation}
\[=
\partial G_1 \times \partial F_{\underline{\tau}(\, \,)} 
\stackrel{\stackrel{\textstyle \tau}{\textstyle \sim}}{\rightarrow}
\coprod\limits_{y \in Y} \partial \bar{G}_y \times \partial F_y
\stackrel{\textstyle \amalg \sigma_y}{\stackrel{\textstyle \sim}{\textstyle \rightarrow}}
\coprod\limits_{y \in Y} \coprod\limits_{x \in f^{-1}(y)} \partial \bar{G}_y  \times
\partial \bar{F}_x =\]
\[\hspace{1.8in}=
\coprod\limits_{x \in X} \partial \bar{F}_x  \times
\partial \bar{G}_{f(x)} =
 \coprod\limits_{x \in X} \partial ( \bar{F}_x \ltimes \bar{G}_{f(x)} )
\]
We define the operation of contraction by
\begin{equation}
\label{eq2.5.9}
( \, \, , \,\, ) : \Delta_X \times \Delta_f \rightarrow \Delta_Y
\end{equation}
\[(G , F ) = \left( (G_1 ; \bar{G}_x ; \tau),(F_y ; \bar{F}_x ; \sigma_y)\right) = 
   \left( G_1 \ltimes \bar{F}_{\underline{\tau}(\, \, )}; F_y \ltimes \bar{G}_{\underline{\sigma_y}(\, )} ; (\tau ,\sigma) \right)
\]
with
\begin{equation}
(\tau , \sigma): \partial (G_1 \ltimes \bar{F}_{\underline{\tau}(\,)})= 
\partial G_1 \times \partial \bar{F}_{\underline{\tau}(\,)}
\stackrel{\tau}{\stackrel{\sim}{\rightarrow}}
\coprod\limits_{x \in X} \partial \bar{G}_x \times \partial \bar{F}_x =
\end{equation}
\[=
\coprod\limits_{y \in Y} \coprod\limits_{x \in f^{-1}(y)} \partial \bar{F}_x  \times
\partial \bar{G}_x 
\stackrel{\amalg \sigma_y^{-1}}{\stackrel{\sim}{\rightarrow}}
\coprod\limits_{y \in Y}  \partial F_y \times \partial \bar{G}_{\underline{\sigma_y}(\,)}=
\coprod\limits_{y \in Y}  \partial (F_y \ltimes \bar{G}_{\underline{\sigma_y}(\,)})
\]

It is straightforward to check that $\Delta$ with these operations satisfy all our axioms
of a generalized ring except the axioms of right-linearity
 (1.7.3)
(also the axioms of commutativity (1.11.1)
and self-adjunction (1.8.8) fails ).

\vspace{10pt}

For a tree $F$, a subset $G \subseteq F$ will be called a (full) \emph{subtree}
with root $b \in G$, if 
$(G, \cS_F|_{G \setminus \{ b\}})$ is a tree, 
i.e. $\cS_F (G \setminus \{ b \} ) \subseteq G$, 
and if moreover for all $a \in G$,
\begin{equation}
\mbox{either} \, \cS_F^{-1}(a) \subseteq G 
\, \mbox{or} \,
\cS_F^{-1}(a)  \cap G = \emptyset
\end{equation}
Given such a subtree $G \subseteq F$, 
and given for each $a \in \partial G$, a (full) subtree
$H_a \subseteq F$ with root $a$, 
and given for $a \in \partial G$, an isomorphism of trees 
$\sigma_a : H_a \stackrel{\sim}{\rightarrow} H$, 
(with a fixed tree $H$),
we form the \emph{commuted}-tree $C^{\sigma}_{G,H}F$:
\begin{equation}
\label{eq2.5.12}
C^{\sigma}_{G,H}F = \left( F \setminus ( G \cup \coprod\limits_{a \in \partial G} H_a ) \right)
\coprod (H \setminus \partial H ) \coprod ( \partial H \times G)
\end{equation}
with $\cS = \cS_{C^{\sigma}_{G,H}F}$ given by
\begin{equation}
\label{eq2.5.11}
\cS (x) = \cS_F(x) \, \mbox{for} \, x \in F \setminus 
          (G \cup \coprod\limits_{a \in \partial G } H_a) , \, \mbox{and} \,
          \cS_F(x) \not\in \coprod\limits_{a \in \partial G} \partial H_a ,
\end{equation}
\[\begin{array}{ll}
 = (y ,a) & \, \mbox{for} \, x \in F \setminus 
          (G \cup \coprod\limits_{a \in \partial G } H_a) , \, \mbox{and} \,
          \cS_F(x) =y \in  \partial H_a , a \in \partial G , \\
= \cS_H(x) & \, \mbox{for} \, x \in H \setminus \partial  H , \, \mbox{and} \,
            x \neq 0_H , \\
= \cS_F(b) & \, \mbox{for} \, x = 0_H , \\
= (y , \cS_G(z)) & \, \mbox{for} \, x = (y , z) \in \partial H \times G , 
                      \, \mbox{and} \, z \neq 0_G ,\\
=\cS_H(y) & \, \mbox{for} \, x = (y , 0_G) \in \partial H \times G 
\end{array}
\]
\begin{center}
\includegraphics[scale=1]{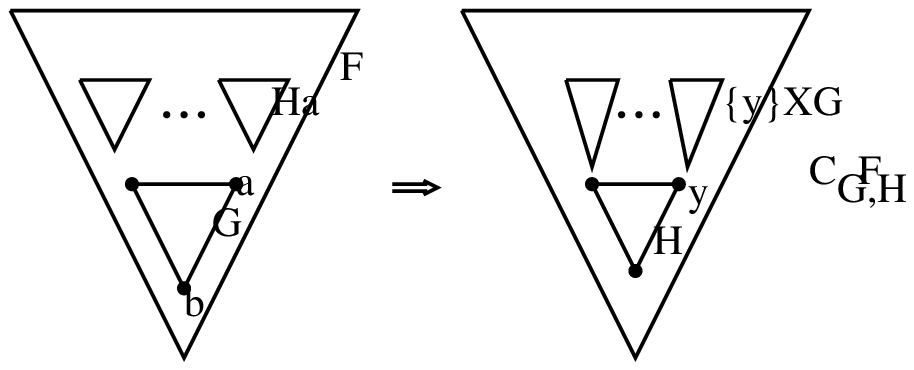}
\end{center}

We let $\approx$ denote the equivalence relation on trees generated by the above operation,
and we let $[F]$ denote the equivalence class of the tree $F$.
Thus $[F] =[F']$ if and only if there is a sequence of trees 
$F= F_0 , F_1 , \ldots  , F_N = F'$,
and there are subtrees $G_j \subseteq F_j$, and 
$H_{j,a} \subseteq F_j$
(for $a \in \partial G_j$), and isomorphisms
$\sigma_{j,a} : H_{j,a} \stackrel{\sim}{\rightarrow} H_j$,
such that $F_{j+1} = C_{G_j , H_j}^{\sigma_j} F_j$.
Note that there is a canonical identification 
$\partial F = \partial (C_{G,H}^{\sigma}F)$,
and hence if $[F] =[F']$ then $\partial F = \partial F'$.

\vspace{10pt}

We can now return to our non-example $\Delta$, 
and redefine $\Delta_X$ as the collection of isomorphism classes of data 
$\left([F_1] ; \{ [\bar{F}_x] \} ; \sigma \right)$,
where the trees $F_1$, and $\bar{F}_x$,
are taken up to equivalence.
It is easy to check that the operations of multiplication (\ref{eq2.5.7}),
and contraction (\ref{eq2.5.9}),
are well defined, and satisfy all our axioms 
(including commutativity (1.11.1)
but not self-adjunction (1.8.8)),
and we have a (commutative but non self adjoint)
generalized ring $\Delta \in {\cal GR}$.

Note that for $f \in \Fbb (X,Y)$ we have 
$1_f = (\{ 0_y \}_{y \in f(x)} ; 
        \{ 0_x \}_{x \in D(f)} ; \sigma_y (0_y) = 0_{f^{-1}(y)})$,
and the structure of $\Delta$ as a functor is given by
\begin{equation}
\left( [F_1] ; \{ [ \bar{F}_x ] \}_{x \in X} ; \sigma \right) \circ 1_{f^t} =
\left( ( [F_1 ] ; \{ [ \bar{F}_x ] \} ; \sigma ) , 1_f \right) =
\end{equation}
\[=
\left( [F_1|_{\sigma^{-1}(\coprod\limits_{x \in D(f)} \partial F_x)} ] ;
       \{ [ \bar{F}_{f^{-1}(y)} ] \}_{y \in f(X)} ; \sigma \right) 
\]  

Note that $\Delta$ is non self-adjoint and the involution is non-trivial on 
$\Delta_{[1]}$, we have:
\[([F_1] ; [\bar{F}_1] ; \sigma )^t = ([\bar{F}_1] ; [F_1] ; \sigma^{-1} )\]

Note that for $X \in \mathbb{F}_{\bullet}$, we have the element 
$\delta_X = ([\{ 0 \} \amalg X ] ; \{ [ 0_x ] \}_{x \in X} ; \sigma ) \in \Delta_X$,
where $\{ 0 \} \amalg X$ is a tree via $\cS (X) =0$,
and $\sigma : X \stackrel{\sim}{\rightarrow} \coprod\limits_{x \in X} \{ 0_x \}$ 
is the obvious map.
These elements generate $\Delta$, and they satisfy 
$\delta_X \circ 1_{f^t} = \delta_{f(X)} \in \Delta_{f(X)} \subseteq \Delta_Y$.

For any commutative generalized ring $A$, we have a natural bijection,
\begin{equation}
\label{eq2.5.15}
\cG \cR ( \Delta , A ) = \{ a = (a_X) \in \prod\limits_{X \in \mathbb{F}} A_X ,
 a_X \circ 1_{f^t} = a_{f(X)},
\ \mbox{all} \ f \in \Fbb (X,Y) \}
\end{equation}
\[\begin{array}{ccc} \varphi &\mapsto& ( \varphi ( \delta_X ))\\
                     \varphi^{(a)} &\leftmapsto & a
\end{array}
\]
To define the homomorphism $\varphi^{(a)}$ it is useful to have the following convention.
For a tree $F$, with height $h = ht(F) = \max \{ ht(b) , b \in F \}$, 
form the tree
\begin{equation}
F^T = F \amalg \coprod\limits_{b \in \partial F} \{ ( b, j) , ht (b)  < j \leq h \} 
\end{equation}
with
\[\begin{array}{l}
\cS (x) = \cS_F(x) \, \mbox{for} \, x \in F ,\\
\cS (b,j) = ( b , j-1) \, \mbox{for} \, b \in \partial F , ht(b) +1 < j \leq h , \\
\cS (b , ht(b) +1) = b
\end{array}
\]
The tree $F^T$ is of constant height.
\[ht(b) = h \quad \mbox{for all } \, b \in \partial F^T ,
\]
and it can be viewed as a sequence of maps
\begin{equation}
\label{eq2.5.17}
F^T_h \stackrel{S_h}{\rightarrow} F^T_{n-1}\rightarrow \cdots \rightarrow
 F^T_1 \stackrel{S_1}{\rightarrow} F^T_0 = \{ 0_F \}
\end{equation}
\[F^T_j = \{ b \in F^T , ht(b) = j \} , \cS_j = \cS |_{F^T_j}
\]
with fibers over $y \in F^T_j$,
\begin{equation}
S_{j+1}^{-1}(y_0) = \left\{\begin{array}{ll}
\cS_F^{-1}(b) & y = b \in F \setminus \partial F , \\
\{ (b , j+1) \} & y= (b,j) 
\end{array}\right.
\end{equation}
With our element $a = (a_X)$, we get elements 
$a_{\cS_{j+1}}^{(y)} \in A_{\cS_{j+1}} = \prod\limits_{y \in F^T_j} 
A_{\cS^{-1}_{j+1}(y)}$, by
\begin{equation}
a^{(y)}_{\cS_{j+1}} = \left\{\begin{array}{ll}
a_{\cS^{-1}_F (b)} & y = b \in F \setminus \partial F ,\\
1 & y = (b,j) ,
\end{array}\right.
\end{equation}
hence a well-defined element
\begin{equation}
a_{[F]}= a_{\cS_1} \circ a_{\cS_2} \circ \cdots \circ a_{\cS_h} \in A_{\partial F}
\end{equation}

Note that since $A$ is assumed to be commutative, equivalent trees give the same element of $A$.
Given the \emph{forest of trees} over $X$,
$\{ [F_x ] , x \in X\}$,
we have similarly a sequence of maps
\begin{equation}
F^T_h \stackrel{\cS_h}{\rightarrow}
F^T_{h-1} \rightarrow
\cdots \rightarrow F^T_1 \stackrel{\cS_1}{\rightarrow}
X \, \mbox{with} \,
h = \max \{ ht (F_x ) , x \in X\}
\end{equation}
With our element $a = (a_X)$ we have similarly a well defined element
\begin{equation}
a_{\{[F_x]\}} = a_{\cS_1} \circ  \cdots \circ a_{\cS_h} \in
A_{\pi} = \prod_{x \in X} A_{\partial F_x} 
\end{equation}
with
\[
\pi = \cS_1 \circ \cdots \circ \cS_h : F^T_h =
\coprod\limits_{x \in X} \partial F_x \rightarrow X
\]
The homomorphism $\varphi^{(a)} : \Delta \rightarrow A $ is given by
\begin{equation}
\label{eq2.5.23}
 \varphi^{(a)}
\left( [F_1] ; \{[\bar{F}_x]\} ; \sigma \right) =
\left( a_{[F_1]} \circ 1_{\sigma^t}  , a_{\{[\bar{F}_x]\}} \right) \in A_X
\end{equation}
That $\varphi$ carries the operations (\ref{eq2.5.7}), (\ref{eq2.5.9})
into the operation of multiplication and contraction in $A$,
follows from the formulas
(1.12.11) , (1.12.12).

Note that the inclusion of the set
 (\ref{eq2.2.16})
in the set (\ref{eq2.5.15}),
correspond to the cannonical homomorphism 
$\pi  \in \cG\cR ( \Delta , \cG ( \mathbb{N} ))$,
\begin{equation}
\pi_X : \Delta_X \rightarrow \cG ( \mathbb{N} )_X = \mathbb{N}^X
\end{equation}
\[\pi_X ( [F_1] ; \{[\bar{F}_x]\} ; \sigma )_{x_0} = \sharp \partial \bar{F}_{x_0}
\]
Finally we remark that the equivalnece relation on trees generated by the commutativity operation (\ref{eq2.5.12}),
is generated by the simpler transposition operation:
given a tree $F$, and $b \in F$,
and for each $a \in \cS_F^{-1}(b)$ a bijection 
$\sigma_a :\cS_F^{-1}(a) \stackrel{\sim}{\rightarrow} [n] $,
we can form the
\emph{transposed tree}
\begin{equation}
\label{eq2.5.25}
C_b^{\sigma} F = ( F \setminus \cS_F^{-1}(b) ) \amalg [n]
\end{equation}
with
\[\begin{array}{lll}
\cS_{c_b^{\sigma} F } (x)& = \cS_F (x) & \mbox{for}\,  \, 
       x \in F\setminus (\cS_F^{-1}(b) \amalg \cS_F^{-2}(b)),\\
    &= \sigma_{\cS_F (x)}(x) \in [n], & \mbox{for} \, \, x \in \cS_F^{-2}(b),\\
    &=b                , & \mbox{for} \, \, x \in [n]
\end{array}
\]
The transposition (\ref{eq2.5.25}) is a special case of the commutation 
(\ref{eq2.5.12}),
but each commutation is obtained by a sequence of transpositions.
It follows that every tree is equivalent to a unique reduced tree,
where $F$ is \emph{reduced} if for all $b \in F$,
either $\nu |_{S^{-1}(b)}$ is non-constant (i.e. there are 
$a_1 , a_ 2 \in \cS^{-1}(b)$ with 
$\nu(a_1 ) \neq \nu(a_2)$),
or it is constant $u_b = \nu(a)$, any $a \in \cS^{-1}(b) $, 
but then $u_b \geq \nu(b)$.
The reduced
trees are characterized, within their equivalence classes, as those trees having the least number of elements.

\section{Free generalized ring $\Delta^W$}
Fix $W \in \mathbb{F}_{\bullet}$.
We say that a tree $F$ is $W$-\emph{labeled}, 
if we are given for all $b \in F$ an injection
\begin{equation}
\mu_b = \mu_b^F : \cS^{-1}_F (b) \hookrightarrow W
\end{equation}
We view $\mu$ as a map $\mu : F \setminus \{ 0_F \} \rightarrow W$,
$\mu (b) = \mu_{\cS (b)} (b)$.
We can give the commutation operation for $W$-labeled trees ananlogous to 
(\ref{eq2.5.12})
(the isomrphisms $\sigma_a$ there,
have to respect the labeling),
we give only the simpler transpositions.

Given a $W$-labeled tree $F$, and $b \in F$ such that for all 
$a \in  \cS^{-1}_F (b) $ we have  
$\mu ( \cS^{-1} (a) ) \equiv W_0 \subseteq W$
we can form the transposed tree
\begin{equation}
\label{eq2.6.2}
C_b F = ( F \setminus \cS^{-1}_F (b) ) \coprod W_0
\end{equation}
with
\[\begin{array}{lll}
\cS_{C_bF}(x) & = \cS_F (x) & \mbox{for}\,  \, x \in 
               F \setminus ( \cS^{-1}_F (b) \coprod \cS^{-2}_F (b) ),\\
             & = \mu (x) \in W_0 & \mbox{for}\,  \, x \in \cS^{-2}_F (b), \\
            & =b & \mbox{for}  \, \, x \in W_0
\end{array}
\]
Note that $C_bF$ is again $W$-labeled by
\begin{equation}
\mu_b^{C_bF} = id_{W_0} : \cS^{-1}_{C_bF} (b) = W_0 \hookrightarrow W
\end{equation}
\[\mu_{w_{\small 0}}^{C_bF} : \cS^{-1}_{C_bF} (w_0) \stackrel{\cS_F}{\stackrel{\rightarrow}{\sim}} \cS^{-1}_F (b) \stackrel{\mu_b^F}{\hookrightarrow} W , \, \mbox{for} \, w_0 \in W_0
\]
\[\mu_x^{C_bF} = \mu_x^F \, \, \mbox{for}\, x \in F \setminus (\{b \} \coprod \cS^{-1}_b (b) )
\]
Note that we have a natural identification 
$\partial (C_bF) = \partial F$.
We write now $[F]_W$ for the equivalence class of the $W$-labeled tree $F$;
we have $[F]_W = [F']_W$ if and only if there is a sequence of $W$-labeled trees
$F = F_0 , F_1 , \ldots , F_N = F'$,
and elements $b_j \in F_j$, such that
for all $a \in \cS^{-1}_{F_j}(b_j)$, 
$\mu^{F_j}(\cS^{-1}(a) ) \equiv W_j$, and $F_{j+1} = C_{b_j}F_j$.

Note that the operations of multiplication (\ref{eq2.5.7}), and contraction
(\ref{eq2.5.9}),
when appplied to trees equipped with $W$-labeling,
yield trees with natural $W$-labeling.
Moreover, these operations descend to operations on equivalence classes of $W$-labeled trees.
We get in this way a commutative generalized ring $\Delta^W$,
\begin{equation}
\label{eq2.6.4}
\Delta_X^W = \left\{\begin{array}{r} ( [F_1]_W ; \{[\bar{F}_x]_W \}_{x \in X} ; \sigma ) 
             /\hspace{-3pt}\approx \, , \,  F_1 , \bar{F}_x \, \,  W\mbox{-labeled trees} \, , \\
              \sigma : \partial F_1 \stackrel{\sim}{\rightarrow}
                \coprod \partial \bar{F}_x \, \mbox{bijection} 
             \end{array}\right\}
\end{equation}
The element $\delta^W = \left( [\{  0 \} \amalg W] ; \{ 0_w \}_{w \in W} ;
            \sigma  \right) \in \Delta_W^W$
generates $\Delta^W$.

For  any commutative generalized ring $A$, we have a bijection
\begin{equation}
\label{eq2.6.5}
\cG\cR (\Delta^W , A) = A_W
\end{equation}
\[\begin{array}{c}\varphi \mapsto \varphi (\delta^W)\\
             \varphi^{(a)} \leftmapsto a
\end{array}
\]
The homomorphism $\varphi^{(a)}$ is given by (\ref{eq2.5.23}), where for a $W$-labeled tree $F$, with associated
sequence of maps (\ref{eq2.5.17}), the element 
$a_{[F]_W} = a_{\cS_1} \circ \cdots \circ a_{\cS_h} \in A_{\partial F}$
is given by
\begin{equation}
a_{\cS_{j+1}}^{(y)} = \left\{\begin{array}{l}
a \circ 1_{\mu_b^F} \in A_{\cS^{-1}_F(y)} , \, \mbox{for} \,
      y \in F \setminus \partial F , ( \mbox{here} \,\mu_y^F \in \mathbb{F}_{\bullet} (\cS^{-1}_F(y) , W))\\
1 \quad \mbox{for} \, y = (b, j)
\end{array}\right.
\end{equation}
and similarly for the forest of $W$-labeled trees $[F_x]_W$,
we get
\[a_{\{[F_x]_W \}} \in A_{\Pi} , \quad \Pi : \coprod\limits_{x \in X} \partial F_x \rightarrow X\]

Fixing a finite family $W_1 , \ldots , W_N \in \mathbb{F}_{\bullet}$,
we say that a tree $F$ is  \emph{$\{W_i\}$- labeled}, 
if for all $b \in F \setminus \partial F$ we are given $j(b) \in [N]$, 
and an injection $\mu_b : \cS^{-1}_F(b) \hookrightarrow W_{j(b)}$.
Given such a tree $F$, and $b \in F$ such that for all 
$a \in \cS^{-1}_F(b)$ we have $j(a) \equiv j$,
and $\mu_a ( \cS^{-1}_F (a) ) \equiv W_0 \subseteq W_j$, 
we can form the transposed tree (\ref{eq2.6.2}).
This  operation generates an equivalence relation, 
and we write $[F]_{W_i}$ for the equivalence class of the $\{ W_i\}$-labeled tree $F$.
Repeating the construction (\ref{eq2.6.4}), with $\{W_i \}$-labeled trees, 
we get a commutative generalized ring $\Delta^{W_1 \cdots W_N}$.
The elements $\delta^{W_i} = \left( 
[\{0\} \amalg W_i] ; \{ 0_w \}_{w \in W_i} ; \sigma
\right) \in \Delta_{W_i}^{W_1 \cdots W_N}$
generate $\Delta^{W_1 \cdots W_N}$, 
and we have for any commutative generalized ring $A$,
\begin{equation}
\label{eq2.6.7}
\cG \cR (\Delta^{W_1 \cdots W_N} , A) = A_{W_1} \times \cdots \times A_{W_N}
\end{equation}
\[\varphi \mapsto \left(\varphi ( \delta^{W_1}) , \ldots , \varphi (\delta^{W_N})\right)
\]
(i.e. $\Delta^{W_1 \cdots W_N} = \Delta^{W_1} \bigotimes\limits_{\mathbb{F}} 
                             \cdots 
                        \bigotimes\limits_{\mathbb{F}} \Delta^{W_N}$).

\vspace{10pt}

Given an injection $j: W \hookrightarrow W'$, 
every $W$-labeled tree $F$ is naturally $W'$-labeled,
and we have an injective homomorphism
\begin{equation}
\Delta^j \in \cG \cR (\Delta^W , \Delta^{W'} ) ,
\Delta^j_X ([F_1]_W ; [\bar{F}_x]_W ; \sigma ) = 
([F_1]_{W'} ; [\bar{F}_x]_{W'} ; \sigma )
\end{equation}
It is dual via (\ref{eq2.6.5}) to the map 
$A_{W'} \twoheadrightarrow A_W$, $a \mapsto a \circ 1_j$.

Converslely, if $F$ is a $W'$-labeled tree, and
\begin{equation}
B_F = \{ b \in \partial F , \mu(\cS^n (b) ) \in W \mbox{for} \, 
         n=0, \ldots , ht(b)-1 \}
\end{equation}
than the reduced tree $F|_{B_F}$, c.f. (\ref{eq2.5.2}), is $W$-labeled.
We have a surjective homomorphism
\begin{equation}
\Delta^{j^t} \in \cG \cR (\Delta^{W'} , \Delta^W ) ,\ee
\[
\Delta^{j^t}_X ([F_1]_{W'} ; [\bar{F}_x]_{W'} ; \sigma ) = 
([F_1|_B]_{W} ; [\bar{F}_x|_{\sigma(B) \cap \partial \bar{F}_x}] ; 
           \sigma|_B )
\]
\[ \mbox{with} \,  B=B_{F_1} \cap \sigma^{-1} (\coprod\limits_{x \in X} B_{\bar{F}_x})
\]
It is dual to the map $A_W \hookrightarrow A_{W'}$, 
$a \mapsto a \circ 1_{j^t}$, and we have 
$\Delta^{j^t} \circ \Delta^j = id_{\Delta^W}$.

\vspace{10pt}

Given a map $f \in Set_{\bullet}(Z,W)$, we have the element
\begin{equation}
\begin{array}{r}
\delta^f = ( [ \{0\} \amalg f(Z) \amalg D(f) ] ;
\{ [0_z]\}_{z \in D(f)} ; \sigma = id_{D(f)} ) \in
\Delta_Z^{W, f^{-1}(w)} =\\
 = (\Delta^W \otimes \bigotimes_{w \in W} 
\Delta^{f^{-1}(w)})_Z
\end{array}
\end{equation}
where the tree $F = \{0\} \amalg f(Z) \amalg D(f)$,
has $\cS_F|_{D(f)}=f$, 
$\cS_F(f(Z)) \equiv 0$, and is labled by 
$\mu_0 : f(Z) \hookrightarrow W$,
and $\mu_w = id_{f^{-1}(w)}$ for $w \in f(Z)$.
This element gives a homomorphism of generalized rings, co-multiplication,
\begin{equation}
\Delta_0^f \in \cG \cR ( \Delta^Z , \Delta^{W,f^{-1}(w)})
\end{equation}
which is dual to multiplication 
$A_W \times \prod\limits_{w \in W} A_{f^{-1}(w)} \rightarrow A_Z$.
On the other hand we have the element
\begin{equation}
\begin{array}{r}
\varepsilon^f = ( [ \{0\} \amalg D(f) ] ;
\{ [ \{0\} \amalg f^{-1}(w) ] \}_{w \in f(Z)};
\sigma = id_{D(f)}) \in 
\Delta_W^{Z, f^{-1}(w)}=\\
=
( \Delta^Z \otimes \bigotimes_{w \in W} \Delta^{f^{-1}(w)} )_W
\end{array}
\end{equation}
giving rise to a homomorphism of generalized rings, co-contraction,
\begin{equation}
\Delta^f_{(\, , \, )} \in \cG \cR (\Delta^W , \Delta^{Z,f^{-1}(w)})
\end{equation}
which is dual to contraction 
$A_Z \times \prod\limits_{w \in W} A_{f^{-1}(w)} \rightarrow A_W$.

\vspace{10pt}

The functor $\Delta : \mathbb{F}_{\bullet} \rightarrow  \cG \cR_C  $, 
with its structure (of co-multiplication, co-contraction,
and co-unit) is thus a co-generalized-ring-object in the tensor category 
$( \cG \cR_C , \bigotimes\limits_{\mathbb{F}})$,
i.e. the dual of our axioms are satisfied.
(Just as the polynomial ring $\mathbb{Z} [X]$, with co-multiplication 
$\Delta_{\bullet} (X) = X \otimes X$,
and co-addition 
$\Delta_+ (X) = X \otimes 1 + 1 \otimes X$, 
is a co-ring object in the tensor category of (commutative) rings and 
$\bigotimes\limits_{\mathbb{Z}}$).

\noindent {\bf Example:}
\setcounter{equation}{15}
Taking for $W =[1]$,
the unit set,
a $[1]$-labeled-tree
is just a ''ladder''
$\{x_0, x_1, \ldots , x_n \}$,
$S(x_j) = x_{j-1}$,
and is determined by its length $n$.
Thus the element $F = ([F_1]; [\bar{F}_x] ; \sigma ) \in \Delta_X^{[1]}$,
is determined by the length of $n$ of $F_1$, 
the point 
$x \in X$ such that $\bar{F}_x \neq 0$, and the length $m$ of 
$\bar{F}_x$.
We write
$F = (x, z^n \cdot (z^t)^m)$,
and we have an isomorphism
\be
\Delta^{[1]} \xrightarrow{\sim}
{\mathbb F} [z^{\mathbb N} \cdot (z^t)^{\mathbb N} ]
\ee
with the generalized ring associated with the free 
monoid on one element,
$M = z^{\mathbb N} \cdot (z^t)^{\mathbb N} \cup \{ 0 \}$.

The self-adjoint quotient $\Delta_+^{[1]}$
of $\Delta^{[1]}$,
representing the functor
\be
\cG \cR_C^+ \rightarrow Mon^+ , \quad
A \mapsto \cG \cR_C^+ (\Delta_+^{[1]}, A) = A_{[1]}
\ee
is isomorphic to the generalized ring associated with 
the free self-adjoint monoid on one element
$M = z^{\mathbb N} \cup \{0 \}$,
\be
\label{eq2.6.18}
\Delta_+^{[1]} \xrightarrow{\sim}
{\mathbb F} [z^{\mathbb N}]
\ee

\section{Limits}
Given a partially ordered set $I$, a functor $A \in (\cG \cR)^I$ is given by objects 
$A^{(i)}\in \cG \cR $ for $i \in I$,
and homomorphism $\varphi_{i,j}: A^{(j)} \rightarrow A^{(i)}$ for $i \leq j , i,j \in I$,
such that $\varphi_{i,j} \circ \varphi_{j,k} = \varphi_{i,k}$ for $i \leq j \leq k $, and 
$\varphi_{ii} = id_{A^{(i)}}$.
The inverse limit of $A$ exists,
and can be computed in $Set_0$. We put
\begin{equation}
(\lim\limits_{\stackrel{\longleftarrow}{I}} A^{(i)})_X =
\{ a= (a^{(i)} ) \in \prod\limits_{i \in I } A_X^{(i)} , 
     \varphi_{i,j}(a_j ) = a_i \mbox{for all} \, i \leq j  \}
\end{equation}
With the operations of componentwise multiplication and contraction 
$\lim\limits_{\longleftarrow} A^{(i)}$ is a generalized-ring (sub-ring of 
$\prod\limits_{i \in I} A^{(i)}$).
We have the universal property
\begin{equation}
\cG \cR (B, \lim\limits_{\stackrel{\longleftarrow}{I}} A^{(i)} ) =
\lim\limits_{\stackrel{\longleftarrow}{I}} \cG \cR (B, A^{(i)}) =
\end{equation}
\[= 
\{ (\psi_i ) \in \prod\limits_{I} \cG \cR  (B, A^{(i)}) , 
\varphi_{i,j} \circ \psi_j = \psi_i \, \mbox{for} \, i \leq j \}
\]

\vspace{10pt}

If the set $I$ is \emph{directed}
(for $j_1 , j_2 \in I$ there is $i \in I$ with $i \leq j_1, j_2$) 
the direct limit of $A$ exists, and can be comoputed in $Set_0$.
We have
\begin{equation}
\label{eq2.7.3}
(\lim\limits_{\stackrel{\longrightarrow}{I}} A^{(i)})_X =
\lim\limits_{\stackrel{\longrightarrow}{I}} A^{(i)}_X =
(\coprod\limits_{i \in I} A_X^{(i)}) / \{a \sim \varphi_{i,j} (a)\}
\end{equation}
There are well defined operations of multiplication and contraction making 
$\lim\limits_{\stackrel{\longrightarrow}{I}} A^{(i)}$ 
into a generalized ring. We have the universal property
\begin{equation}
\cG \cR (\lim\limits_{\stackrel{\longrightarrow}{I}} A^{(i)} , B) =
\lim\limits_{\stackrel{\longleftarrow}{I}} \cG \cR (A^{(i)}, B) =
\end{equation}
\[=
\{ (\psi_i ) \in \prod\limits_{I} \cG \cR  (A^{(i)}, B) , 
\psi_i \circ \varphi_{i,j} = \psi_j  \, \mbox{for} \, i \leq j \}
\]

\chapter{Ideals}
\section{Equivalence ideals}
%\begin{definition}
\noindent (3.1.1) {\bf Definition:}
For  $A \in \cG \cR$ an \emph{equivalence ideal} $\ca$  is a collection of subsets $\ca_X \subseteq A_X \times A_X$,
such that $\ca_X$ is an equivalence relation on $A_X$,
(we write $a \sim a'$ for 
$(a,a') \in \ca_X$;
and for $a = (a_y)$, $a' = (a'_y) \in A_f$, we write $a \sim a'$ for
$(a_y , a'_y ) \in \ca_{f^{-1}(y)}$
 for all $y \in Y$),
and $\ca$  respects the operations: if $a \sim a'$ than for all $b \in A$, 
\setcounter{equation}{1}
\begin{equation}
\label{eq3.1.2}
\begin{array}{ll}
a \circ b \sim a' \circ b , & b \circ a \sim b \circ a'\\
(a,  b) \sim (a',  b) , & (b ,  a)  \sim ( b ,  a')
\end{array}
\end{equation}
(whenever these operations are defined.)

We let $eq(A)$ denote the set of equivalence ideals of $A$.

For $ \ca \in eq(A)$ we can form the quotient $A/\ca$, with 
\be
(A/\ca)_X = A_X/\ca_X = \ca_X\mbox{-equivalence classes in}
\, A_X
\ee
There is a natural surjection $\pi_X : A_X \rightarrow A_X/\ca_X$,
and there is a unique structure of generalized rings on 
$A/\ca$ such that $\pi \in \cG\cR (A,A/\ca)$.

\vspace{10pt}

\noindent (3.1.4) {\bf Definition} For $\varphi \in \cG\cR (A,A')$, we let
\[KER(\varphi)_X = \{ (a_1, a_2) \in A_X \times A_X , \varphi (a_1 ) = \varphi_X(a_2) \}
\]
$KER (\varphi )$ is an equivalence ideal.

\vspace{10pt}

\noindent (3.1.5) We have the universal property of the quotient
\[
\cG\cR (A/\ca, A' ) = \{ \varphi \in \cG\cR (A, A') ,  
KER ( \varphi)  \supseteq \ca \}
\]

\vspace{10pt}

\noindent (3.1.6) Every homomorphism 
$ \varphi \in \cG\cR (A, A')$ has a

 canonical factorization
\[
\varphi = j \circ \tilde{\varphi} \circ \pi
\]
\setcounter{equation}{6}

\begin{center}
\begin{diagram}
A & \rTo^{\varphi}& A'\\
\dOnto^{\pi}& &\uInto_{j}\\
A/KER(\varphi) & \rTo_{\sim}^{\widetilde{\varphi}}& \varphi(A) 
\end{diagram}
\end{center}

with $\pi$ surjection, $j$ injection, and $\tilde{\varphi}$ an isomorphism.

\section{Ideals}
\noindent (3.2.1) {\bf Definition :} For $A \in \cG\cR$, 
an \emph{ideal} $\sqa$ is a collection of subsets $\sqa_X \subseteq A_X$,
with $0_X \in \sqa_X$, and with $A \circ \sqa$, $\sqa \circ A$, $(A,\sqa)$, $(\sqa,A) \subseteq \sqa$.
We let $il(A)$ denote the set of ideals of $A$.

\vspace{10pt}

For $\ca \in eq(A)$, we have the associated ideal $Z(\ca)$:
\setcounter{equation}{1}
\be
Z(\ca)_X = \left\{a \in A_X , (a,0_X) \in \ca_X \right\}
\ee
For $\sqa \in il(A)$, we have the equivalence ideal $E(\sqa)$ generated by $\sqa$, 
it is the intersection of all equivalence ideals containing $(a,0)$ for all 
$a \in \sqa$.

These give a Galois correspondence:
\be
il(A)
\begin{array}{c}
  \stackrel{\scriptstyle E}{\longrightarrow} \\
  \stackrel{\textstyle \longleftarrow}{\scriptstyle Z}
\end{array}
eq(A)
%il(A) \stackrel{\textstyle \xrightarrow{E}}{\stackrel{\textstyle\leftarrow}{\scriptsize Z}} eq(A)
\ee
It is monotone,
\be
\begin{array}{ll}
\sqa_1 \subseteq \sqa_2 & \Rightarrow E(\sqa_1) \subseteq E(\sqa_2)\\
\ca_1 \subseteq \ca_2 & \Rightarrow Z(\ca_1) \subseteq Z(\ca_2)
\end{array}
\ee
and we have
\be
\sqa \subseteq ZE(\sqa) \quad , \quad EZ(\ca) \subseteq \ca 
\ee
It follows that we have
\be
ZEZ(\ca) = Z(\ca) \quad , \quad E(\sqa) = EZE(\sqa) ,
\ee
and the maps $E$, $Z$ induce inverse bijections
\be
E\mbox{-}il(A) \stackrel{\textstyle \sim}{\longleftrightarrow} 
Z\mbox{-}eq(A)
\ee
with 
\be
E\mbox{-}il(A) = \{ \sqa \in il(A) , \sqa = ZE(\sqa)\} =
          \{ Z(\ca) , \ca \in eq(A)\} ,
\ee
the \emph{stable} ideals, and
\be
Z\mbox{-}eq(A) = \{ \ca \in eq(A) , \ca = EZ(\ca)\} =
         \{E(\sqa) , \sqa \in il(A)\}
\ee

Let $\sqa \in il(A)$, and let $\ca= (\ca_X)$ be defined by
\be
\label{eq3.2.10}
\ca_X=
\ee
\[=
\left\{\begin{array}{l}
  (a, a') \in A_X \times A_X, \mbox{there exists a \emph{path}}\,  
a = a_1 , a_2, \ldots , a_n =a',\\
 \mbox{with} \,
\{a_j , a_{j+1} \} \,  \mbox{of the form 
   either} \, \{ (d_j \circ c_j , b_j ) , (d_j \circ \bar{c}_j , b_j ) \},\\
\mbox{or} \ \{ (d_j , b_j \circ c_j ) , (d_j , b_j \circ \bar{c}_j ) \}  
\ \mbox{with} \ 
    d_j \in A_{Y_j} , b_j \in A_{f_j} ,\\ f_j \in Set_{\bullet}(Z_j , X) ,
    c_j , \bar{c}_j \in A_{g_j} ,
 g_j \in Set_{\bullet}(Z_j , Y_j ) \ \mbox{or resp.} \ \\
   % d_j \in A_{Y_j} , b_j \in A_{f_j} ,\\ f_j \in Set_{\bullet}(Z_j , X) ,
    %c_j , \bar{c}_j \in A_{g_j} , 
g_j \in Set_{\bullet}(Y_j , Z_j )\ 
   \mbox{and with} \ c_j^{(z)} = \bar{c}_j^{(z)}, \mbox{ or }\
     c_j^{(z)},  \bar{c}_j^{(z)} \in \sqa , \\
\mbox{for all} \ z \in Y_j  
    (\mbox{resp.} \, z \in Z_j) 
\end{array}\right\}
\]
\noindent (3.2.11) {\bf Claim:} $E(\sqa) = \ca$
\setcounter{equation}{11}

{\bf Proof:}
It is clear  that $\ca_X$ is an equivalence relation on $A_X$,
and that $\ca_X \subseteq E(\sqa)_X$,
and we need to show that $\ca$ respects the operations (\ref{eq3.1.2}).
If $(a,a') \in \ca$, 
so there is a path
$a = a_1 , \ldots , a_n =a'$ 
as in (\ref{eq3.2.10}),
then $h \circ a_j$ (resp. $a_j \circ h$, $(h, a_j)$, $(a_j , h)$) 
is a path, showing $(h \circ a , h \circ a')$
 (resp. $(a \circ h , a' \circ h)$, 
$((h,a), (h,a'))$, $((a,h),(a',h)))$) is in $\ca$.
This follows from the identities
\be
\label{eq3.2.12}
\begin{array}{clcl}
  & h \circ (d \circ c, b) = \left( (h \circ d) \circ c , b \right) & , &
   h \circ (d, b \circ c) = (h \circ d , b \circ c)\\
\mbox{resp.}& (d \circ c ,b) \circ h = (d \circ (c \circ \tilde{h}), \tilde{b})
&, & (d, b \circ c)\circ h = (d \circ \tilde{h} , \tilde{b} \circ \tilde{c})\\
& (h,(d \circ c,b)) = (h \circ b, d \circ c) & ,& (h,(d,b \circ c)) =
  ((h \circ b)\circ c , d)\\
& ((d \circ c,b), h) = (d \circ c , h \circ b) &,&
((d, b \circ c), h) = (d, (h \circ b)\circ c)
\end{array}
\ee

\vspace{10pt}

It follows that $\sqa \in E\mbox{-}il(A)$, $\sqa = ZE(\sqa)$, 
if and only if for all $b,d,c,\bar{c}$ as in (\ref{eq3.2.10}),
\be
\begin{array}{cc}
 & (d \circ c, b) \in \sqa \Leftrightarrow (d \circ \bar{c} , b) \in \sqa\\
\mbox{and} & (d, b \circ c) \in \sqa \Leftrightarrow (d, b \circ \bar{c}) \in \sqa
\end{array}
\ee

\section{Operations on ideals}

The intersection of ideals is an ideal,
\be
\sqa_i \in il(A) \Rightarrow \bigcap\limits_i \sqa_i \ \in il(A)
\ee

Given a collection $B = \{ B_X \subseteq A_X \}$,
the ideal generated by $B$ will be denoted by $\{B\}_A$,
it is the intersection of all ideals containing $B$.
If $B$ satisfies $B \circ A \subseteq B$,
it can be described explicitly as
\be
\label{eq3.4.2}
\{B\}_{A,X} =
\ee
\[= 
\left\{ a \in A_X , a=(d \circ c, b)\ 
\mbox{or} \ a=(d, b \circ c)\
\mbox{for some}\  d \in A_Y,
              b \in A_f , 
\right.\]
\[f \in Set_{\bullet} (Z,X),
 c \in A_g , \ 
              g \in Set_{\bullet}(Z,Y),\]
\[ \ \mbox{or resp.} \
 g \in Set_{\bullet}(Y,Z)
\]
\[ \left. \mbox{and with} \ c= (c^{(z)}) , c^{(z)} \in B_{g^{-1}(z)} \,  
           \mbox{for all} \, z \in Y \,  (\mbox{resp.} \, z \in Z)
\right\}
\]

It is clear that the set described in (\ref{eq3.4.2}), contains $B_X$, 
and is contained in $\{B\}_{A,X}$,
and we only have to check that it is an ideal - this follows from the identities (\ref{eq3.2.12}).

In particular, for $\sqa_i \in il(A)$, 
we can take $B=\bigcup\limits_i \sqa_i$, and we obtain the smallest ideal
containing all the $\sqa_i$'s.
\be
\sum\limits_i \sqa_i = \left\{ \bigcup\limits_i \sqa_i\right\}_A
\ee
Thus $il(A)$ is a complete lattice, with minimal element the zero ideal 
$0 = \{0_X\}$,
and maximal element the unit ideal $\{1\}_A = \{A_X\}$.

Note that for arbitrary $B$ we can similarly describe the ideal it generates as

\be
\label{eq3.4.4}
\{B\}_{A,X} = 
\ee
\[=\left\{a \in A_X , a= (d\circ c \circ e , b) 
           \,\mbox{or} \, a = (d,b \circ c \circ e) , \mbox{with} \,
            c^{(z)}\in B \, \mbox{for all }\, z \right\}
\]

Note that for a subset $B \subseteq A_{[1]}$, we have
\be
\{B\}_{A,X} = \left\{ a \in A_X , a= (d, b \circ c) , d \in A_Y , 
              b \in A_f , 
      \right.\ee
\[\left.          
f \in Set_{\bullet}(Y,X),
 c = (c^{(y)}) \in A_{id_Y} =
               (A_{[1]})^Y, c^{(y)}\in B \cup B^t 
               \right\}
\]
We do not need the $e$'s in (\ref{eq3.4.4})
because we have commutativity $c \circ e = e \circ \tilde{c}$,
c.f. (\ref{eq1.10.3}),
and we do not need the two shapes of (\ref{eq3.4.2}) because the 
$(A_{[1]})^Y$-action is self-adjoint
$(d, b \circ c) = (d \circ c^t , b)$, 
c.f. (1.10.6).

\section{Homogeneous ideals} 
An ideal $\sqa \in il(A)$ is called \emph{homogeneous}
if it is generated by $\sqa_{[1]}$.
The subset $\sqa_{[1]} \subseteq A_{[1]}$ satisfies for all $X \in \Fbb$,
\be
\label{eq3.5.1}
\left( A_X , A_X \circ (\sqa_{[1]})^X \right) \subseteq \sqa_{[1]}
\ee
Conversly, if a subset $\sqa_{[1]} \subseteq A_{[1]}$ satisfies (\ref{eq3.5.1}),
then $\sqa_{[1]}^t = \sqa_{[1]}$ (because for $a \in \sqa_{[1]}$, 
$a^t = (1, 1 \circ a) \in \sqa_{[1]}$),
and $\sqa_{[1]} \circ A_{[1]} = \sqa_{[1]}$ (because for $a \in \sqa_{[1]}$, 
$b \in A_{[1]}$, $a \circ b = (b,1 \circ a^t ) \in \sqa_{[1]}$).
Moreover, the ideal $\sqb$ generated by
$\sqa_{[1]}$, has
\be
\sqb_X = \bigcup\limits_{f \in Set_{\bullet}(Y,X)}
      \left( A_Y , A_f \circ (\sqa_{[1]})^Y\right)
\ee
and in particular $\sqb_{[1]} = \sqa_{[1]}$.
Thus we identify the set of homogeneous ideals with the collection of 
subsets $\sqa_{[1]} \subseteq A_{[1]}$ satisfying (\ref{eq3.5.1}),
and we denote this set by $[1]\mbox{-}il(A)$,
\be
\label{eq3.5.3}
[1]\mbox{-}il(A) = \{ \sqa \subseteq A_{[1]} , (A_X , A_X \circ (\sqa)^X)
       \subseteq \sqa \}
\ee
The set $[1]\mbox{-}il(A)$ is a complete lattice, with minimal element $\{0\}$,
maximal element $\{1\}_A = A_{[1]}$.
For $\sqa_i \in [1]\mbox{-}il(A)$, we have 
$\bigcap\limits_i\sqa_i \in [1]\mbox{-}il(A)$, and 
$\sum\limits_i \sqa_i \in [1]\mbox{-}il(A)$,
where
\be
\label{eq3.5.4}
\sum\limits_i \sqa_i =
\ee
\[=
 \left\{ a \in A_{[1]}, a = (b, d \circ c) , 
              b,d \in A_X , c=(c^{(x)}) \in (\bigcup\limits_i \sqa_i)^X
           \subseteq A_{id_X} , X \in \Fbb \right\}
\]
Note that the homogeneous ideal generated by elements 
$a_i \in A_{[1]}$, $i \in I$,
can be described as 
\be
\label{eq3.5.5}
\{a_i\}_A =
\ee
\[=
\left\{ a \in A_{[1]} , a= (b,d \circ c) ,
b,d \in A_X , c=(c^{(x)}) \in A_{id_X},
c^{(x)}=a_{i(x)}\right. 
\]
\[\left.\mbox{or} \, c^{(x)}= a^t_{i(x)}
\mbox{for} \, x \in X  
\right\}
\]
We have also the operation of multiplication of homogeneous ideals. 
For $\sqa_1 , \sqa_2 \in [1]\mbox{-}il(A)$,
we let $\sqa_1 \cdot \sqa_2$ denote the homogeneous ideal 
generated by the product 
$\sqa_1 \circ \sqa_2 = \{ a_1 \circ a_2 , a_i \in \sqa_i \}$.

Thus
\be
\sqa_1 \cdot \sqa_2 =
\ee
\[= \left\{ a \in A_{[1]} , a=(b,d \circ c) , b,d \in A_X ,
c=(c^{(x)}) \in (\sqa_1 \circ \sqa_2)^X \subseteq A_{id_X} \right\}
\]
This operation is associative, we have for 
$\sqa_1, \sqa_2, \sqa_3 \in [1]\mbox{-}il(A)$,
\be
(\sqa_1 \cdot \sqa_2) \cdot  \sqa_3 =
 \sqa_1 \cdot \sqa_2 \cdot \sqa_3 = 
\sqa_1\cdot ( \sqa_2 \cdot \sqa_3) 
\ee
with $\sqa_1 \cdot \sqa_2 \cdot \sqa_3 = $\\
$=\left\{ a \in A_{[1]} , a = 
(b,d \circ c) , b,d \in A_X ,
c=(c^{(x)}) \in (\sqa_1 \circ \sqa_2 \circ \sqa_3)^X \subseteq A_{id_X} \right\}$
(use $(b, d \circ a_1 \circ a_2)\circ a_3 = 
(b, d \circ a_1 \circ a_2 \circ \tilde{a}_3^t)$ and
$a_1 \circ (b, d \circ a_2 \circ a_3) =
(b, d \circ a_2 \circ a_3 \circ \tilde{a}_1^t)$).

The multiplication of homogeneous ideals is clearly commutative, 
$\sqa_1 \circ \sqa_2 = \sqa_2 \circ \sqa_1$,
has unit element $\{1\}_A= A_{[1]}$,
$\sqa \cdot \{1\}= \sqa$,
and has zero element $\{0\}$, $\sqa \cdot \{0\} = \{0\}$.
Thus 
$[1]\mbox{-}il(A) \in Mon^+$.

For a homomorphism $\varphi \in \cG\cR (A,B)$,
and for $\sqb \in il(B)$, 
(resp. $\sqb \in [1]\mbox{-}il(B)$),
its inverse image $\varphi^{*} (\sqb)_X = \varphi_X^{-1}(\sqb_X)$,
(resp. 
$\varphi^*_{[1]}(\sqb) = \varphi^{-1}_{[1]}(\sqb) \subseteq A_{[1]}$) is clearly a (resp. homogeneous) 
ideal of $A$.
For $\sqa \in il(A)$,
(resp. $[1]\mbox{-}il(A)$)
let $\varphi_*(\sqa) \subseteq B$ denote the (homogeneous) 
ideal generated by the image $\varphi(\sqa)$ (resp. 
$\varphi_{[1]}(\sqa)$).
We have  Galois correspondences
\be
[1]\mbox{-}il(A) 
\begin{array}{c}
\stackrel{\textstyle \varphi^*_{[1]}}{\longleftarrow}\\
\stackrel{\textstyle\longrightarrow}{\varphi_{*[1]}}
\end{array}
[1]\mbox{-}il(B)
\quad
\mbox{and} \quad
il(A)
\begin{array}{c}
\stackrel{\textstyle \varphi^*}{\longleftarrow}\\
\stackrel{\textstyle \longrightarrow}{\varphi_*}
\end{array}
il(B)
\ee
The maps $\varphi^*$, $\varphi_*$ are monotone, and satisfy
\be
\label{eq3.5.9}
\sqa \subseteq \varphi^* \varphi_* (\sqa) ,
\quad
 \varphi_* \varphi^* (\sqb) \subseteq \sqb
\ee
It follows that we have
\be
\varphi_*(\sqa) = \varphi_*\varphi^*\varphi_*(\sqa) \quad , \quad
\varphi^*(\sqb) =  \varphi^*\varphi_* \varphi^* (\sqb)
\ee
and $\varphi^*$, $\varphi_*$ induce inverse bijections,
\[\hspace{-6.5cm}\left\{ \sqa \in [1]\mbox{-}il(A) ,
\sqa = \varphi^*\varphi_* (a) \right\} =\]
\[=
\left\{ \varphi^* (\sqb) , \sqb \in [1]\mbox{-}il(B) \right\}
\stackrel{\sim}{\longleftrightarrow}
\left\{ \varphi_* (\sqa) , \sqa \in[1]\mbox{-}il(A) \right\}=
\]
\[\hspace{5cm}
=
\left\{ \sqb \in [1]\mbox{-}il(B) , \sqb = \varphi_*\varphi^* (\sqb) 
 \right\}
\]
and similarly with $il(A)$ and $il(B)$.

\vspace{10pt}

\noindent (3.4.11) {\bf Definition:} 
\setcounter{equation}{11}
For an equivalence ideal $\varepsilon \in eq(A)$, 
and for an ideal or a homogeneous ideal $\sqa$,
we say $\sqa$ is $\varepsilon$-\emph{stable} if for all 
$(a,a') \in \varepsilon : a \in \sqa \Leftrightarrow a' \in \sqa$.
We denote by $il(A)^{\varepsilon}$ 
(resp. $[1]\mbox{-}il(A)^{\varepsilon}$) 
the set of $\varepsilon$-stable (homogeneous) ideals.

Letting $\pi_{\varepsilon} : A \twoheadrightarrow A/\varepsilon$ 
denote the canonical homomorphism, we have bijections
\be
\begin{array}{c}
il(A)^{\varepsilon} \stackrel{ \sim}{\longleftrightarrow}
il(A/\varepsilon) \ \mbox{and} \
{[1]}\mbox{-}il(A)^{\varepsilon} \stackrel{\sim}{\longleftrightarrow}
{[1]}\mbox{-}il(A/\varepsilon)\\
\sqa \mapsto \pi_{\varepsilon}(\sqa)\\
\pi^{-1}_{\varepsilon}(\sqb) \leftmapsto \sqb
\end{array}
\ee
\noindent (3.4.12) {\bf Definition:}
\setcounter{equation}{12}
For an (resp. homogeneous) ideal $\sqa$,
we say $\sqa$ is \emph{stable} if it is $E(\sqa)$-stable.
We denote by $E\mbox{-}il(A)$, (resp. $E[1]\mbox{-}il(A)$), the set of stable (homogeneous) ideals.
Note that by the explicit description of $E(\sqa)$ in  
(\ref{eq3.2.10}),
a subset $\sqa \subseteq A_{[1]}$ is a stable homogeneous ideal if and only if
\be
\mbox{for} \; X,Y \in \Fbb , b, d \in A_{X \oplus Y},
   c, \bar{c} \in (A_{[1]})^{X \oplus Y},
\ee
\[
\mbox{with} \;
c^{(x)}= \bar{c}^{(x)}
\mbox{for} \; x \in X
\mbox{and} \; c^{(y)} , \bar{c}^{(y)} \in \sqa
\mbox{for} \; y \in Y,
\]
\[
 \mbox{have}: ( b,d \circ c) \in \sqa \Leftrightarrow 
(b, d \circ \bar{c}) \in \sqa
\]

(taking $X = [0]$, $\bar{c}^{(y)}\equiv 0$, 
we see that this condition includes $\sqa$ being a homogeneous ideal).

\section{$h$-ideals}
\noindent(3.6.1) {\bf Definition:}
\setcounter{equation}{1}
For $A \in \cG\cR$, a subset $\sqa \subseteq A_{[1]}$ will be called \emph{$h$-ideal}
if for all $X \in \Fbb$,
$b,d \in A_X$, $c = (c^{(x)}) \in (\sqa)^X \subseteq A_{id_X}$, 
we have
\be
(b \circ c, d) \in \sqa
\ee
We denote by $h\mbox{-}il(A)$ the set of $h$-ideals of $A$.

Comparing with the definition of homogeneous ideals 
(\ref{eq3.5.3}),
we have
\be
[1]\mbox{-}il(A) = \{ \sqa \in h\mbox{-}il(A) , \sqa = \sqa^t \}
\ee
The set $h\mbox{-}il(A)$ is a complete lattice,
with minimal element $(0)$,
maximal element $(1)=A_{[1]}$.
For $\sqa_i \in h\mbox{-}il(A)$, we have 
$\bigcap\limits_i \sqa_i \in h\mbox{-}il(A)$,
and $\sum\limits_i \sqa_i \in h\mbox{-}il(A)$,
where (c.f. (\ref{eq3.5.4}))
\be
\sum\limits_i \sqa_i =
\left\{ a \in A_{[1]}, a = (b \circ c , d) , b,d \in A_X ,
c = (c^{(x)}) \in (\bigcup\limits_i a_i )^X \subseteq A_{id_X} \right\}
\ee
Note that the $h$-ideal generated by elements $a_i \in A_{[1]}$,
$i \in I$, can be described as (c.f. (\ref{eq3.5.5}))
\[(a_i)_A=
\left\{ a \in A_{[1]}, a = (b \circ c , d) , b,d \in A_X ,
c = (c^{(x)}) \in \{a_i\}^X \subseteq A_{id_X} \right\}
\]
In particular, for $a \in A_{[1]}$ the \emph{principal} 
$h$-ideal generated by $a$ is just
\be
(a)_A = a \circ A_{[1]}
\ee
Indeed, if $c^{(x)}= a \circ e_x$, $e_x \in A_{[1]}$ for 
$x \in X$,
than for $b,d \in A_X$,
\[(b \circ c , d ) = a \circ (b \circ e , d) \in a \circ A_{[1]}
\]
(while the homogeneous ideal generated by $a$ is the $h$-ideal 
generated by $a$ and $a^t$,
c.f. (\ref{eq3.5.5}) for $\{a\}_A$).

We have multiplication of $h$-ideals, for 
$\sqa_1 , \sqa_2 \in h\mbox{-}il(A)$,
\be
\sqa_1 \cdot \sqa_2 = 
\ee
\[=\left\{ a \in A_{[1]},
a = (b \circ c, d ), \, b,d \in A_X ,
c = (c^{(x)}) \in (\sqa_1 \circ \sqa_2)^X \subseteq A_{id_X}
\right\}
\]
It is associative,
\be
(\sqa_1 \cdot \sqa_2) \cdot \sqa_3 =
\sqa_1 \cdot \sqa_2 \cdot \sqa_3 =
\sqa_1 \cdot (\sqa_2 \cdot \sqa_3 )
\ee
\[ \mbox{with} \; 
\sqa_1 \cdot \sqa_2 \cdot \sqa_3 =
\left\{ a \in A_{[1]},
a = (b \circ c, d ),  b,d \in A_X ,\right.\]
\[\left.
c = (c^{(x)}) \in (\sqa_1 \circ \sqa_2 \circ \sqa_3)^X \subseteq A_{id_X}
\right\}\]

(Use now:
$(b \circ a_1 \circ a_2 , d ) \circ a_3 =
(b \circ a_1 \circ a_2 \circ \tilde{a}_3 , d)$,
and
$a_1 \circ (b \circ a_2 \circ a_3 , d ) =
 (b \circ \tilde{a}_1 \circ a_2 \circ a_3 , d )$).

It is clearly commutative, $\sqa_1 \cdot \sqa_2 = \sqa_2 \cdot \sqa_1$;
has unit $(1)$,
$\sqa \cdot (1) = \sqa$;
has zero $(0)$, $\sqa \cdot (0) = (0)$.

We can also divide $h$-ideals. For $\sqa_0 , \sqa_1 \in h\mbox{-}il(A)$,
we let
\be
(\sqa_0 : \sqa_1 ) = 
\{ c \in A_{[1]} , c \circ \sqa_1 \subseteq \sqa_0 \}
\ee
This is an $h$-ideal,
$(\sqa_0 : \sqa_1 ) \in h\mbox{-}il(A)$.
Indeed, for $c^{(x)} \in (\sqa_0 : \sqa_1 )$,
and for any $b,d \in A_X$,
and any $a_1 \in \sqa_1$,
we have
\[(b \circ (c^{(x)}), d ) \circ a_1 =
( b\circ (c^{(x)} \circ a_1 ) , d ) \in \sqa_0
\]
and so $(b \circ (c^{(x)}), d ) \in (\sqa_0 : \sqa_1 )$.

\vspace{10pt}

For elements $m_1 , m_2 \in A_X$, we have their \emph{annihilator}
\be
ann_A (m_1, m_2 ) = 
\{a \in A_{[1]} , a \circ m_1 = a \circ m_2 \}
\ee
This is an $h$-ideal, 
$ann_A(m_1, m_2 ) \in h\mbox{-}il(A)$.
Indeed for $c^{(y)} \in ann_A (m_1, m_2 )$,
and for any $b,d \in A_Y$ we have
\[( b \circ (c^{(y)}),d) \circ m_1 =
(b \circ (c^{(y)}\circ m_1 ) , \tilde{d}) =
(b \circ (c^{(y)}\circ m_2), \tilde{d}) =
(b \circ (c^{(y)}) , d ) \circ m_2
\]
and so 
$(b \circ (c^{(y)}), d) \in ann_A (m_1, m_2)$

\vspace{10pt}

Let $\varphi \in \cG\cR(A,B)$.
For an $h$-ideal $\sqb \in h\mbox{-}il(B)$,
its inverse image $\varphi^*(\sqb) = \varphi^{-1}_{[1]}(\sqb) \subseteq A_{[1]}$
is clearly an $h$-ideal of $A$.
For an $h$-ideal $\sqa \in h\mbox{-}il(A)$,
we let $\varphi_*(\sqa) \subseteq B_{[1]}$ denote the
$h$-ideal generated by the image $\varphi_{[1]}(\sqa)$,
$\varphi_*(\sqa) = 
\lim\limits_{\stackrel{\textstyle \rightarrow}{X}}
(B_X \circ ( \varphi_{[1]}(\sqa))^X , B_X )$.

We have a Galois correspondence
\be
h\mbox{-}il(A)
\begin{array}{c} 
\stackrel{\textstyle \varphi^*}{ \longleftarrow}\\
\stackrel{\textstyle \longrightarrow}{\varphi_*}
\end{array}
h\mbox{-}il(B)
\ee
The maps $\varphi^*$, $\varphi_*$ are monotone, and satisfy
\be
\sqa \subseteq \varphi^* \varphi_* (\sqa) 
\quad, \quad
\varphi_*\varphi^*(\sqb) \subseteq \sqb
\ee
It follows that we have
\be
\varphi_*(\sqa) = \varphi_*\varphi^*\varphi_*(\sqa) 
\quad ,\quad
\varphi^*(\sqb)=\varphi^*\varphi_*\varphi^*(\sqb)
\ee
and $\varphi^*$, $\varphi_*$ induce inverse bijections,
\be
\hspace{-4cm}
\{\sqa \in h\mbox{-}il(A) , \sqa=\varphi^*\varphi_* (\sqa) \}=
\ee
\[=
\{\varphi^*(\sqb) , \sqb \in h\mbox{-}il(B) \}
\stackrel{\sim}{\longleftrightarrow}
\{\varphi_*(\sqa) , \sqa \in h\mbox{-}il(A)\} =
\]\[\hspace{4cm}
=
\{\sqb\in h\mbox{-}il(B) , \sqb = \varphi_* \varphi^*\varphi_* (\sqb) \}
\]

\vspace{10pt}

In summary,
for a general $A \in \cG\cR$, we have defined the sets
\be
\label{eq3.6.14}
\begin{array}{ccccccc}
h\mbox{-}il(A)& \supseteq & [1]\mbox{-}il(A) & \subseteq &il(A) & & eq(A)\\
       &           & |\hspace{-2pt}\bigcup &                &|\hspace{-2pt}\bigcup& & |\hspace{-2pt}\bigcup\\
        &           &E\mbox{-}[1]\mbox{-}il(A)&\subseteq&E\mbox{-}il(A)&
        \stackrel{\sim}{\leftrightarrow}&Z\mbox{-}eq(A)
\end{array}
\ee

For a self adjoint $A \in \cG\cR^+$ we have equality 
$h\mbox{-}il(A) =[1]\mbox{-}il(A)$.
It is easy to check that for $A=\cG (B)$, $B$ a commutative ring,
all the inclusions in (\ref{eq3.6.14})
are equalities, and are identified with the set of ideals of $B$.

\chapter{Primes and Spectra}
\section{Maximal ideals and primes}

We say that an equivalence ideal $\varepsilon \in eq(A)$
is \emph{proper} if $(1,0) \not\in \varepsilon$,
or equivalently 
$\varepsilon_X  
\subsetneqq
 A_X \times A_X$ 
for some/all $X \in \Fbb$,
or equivalently $A/\varepsilon  \neq 0$.
We say that an ideal, or an $h$-ideal,
$\sqa$ is \emph{proper} if $1 \not\in \sqa$, 
or equivalently 
$\sqa_{[1]} \subsetneqq A_{[1]}$.
Since a union of a chain of proper $h$-ideals is again a proper $h$-ideal,
an application of Zorn's lemma gives

\noindent (4.1.1) {\bf Proposition:} For $A \in \cG\cR$, there exists maximal
proper $h$-ideal.

We let $Max(A) \subseteq h\mbox{-}il(A)$ denote the set of maximal ideals.

\vspace{10pt}

\noindent (4.1.2) {\bf Definition:} A (proper) $h$-ideal 
$\sqp \in h\mbox{-}il(A)$ is called \emph{prime} if 
$A_{[1]}\setminus \sqp$ is closed with respect to multiplication, 
i.e. if for all $a,b \in A_{[1]}$,
\[a \circ b \in \sqp \quad \mbox{implies} \quad  a \in \sqp 
\quad \mbox{or} \quad b \in \sqp
\]
We let $spec(A) \subseteq h\mbox{-}il(A) $ denote the set of primes of $A$.

\vspace{10pt}

\noindent (4.1.3) {\bf Proposition:} $Max(A) \subseteq spec(A)$.

\noindent {\bf Proof:} Let $\sqp \in Max(A)$,
and take any elements $a, a' \in A_{[1]} \setminus \sqp$.
Since $\sqp$ is maximal, the $h$-ideals $(\sqp, a)_A$
and $(\sqp, a')_A$ are the unit ideal.
We have therefore $1=(b \circ c, d)$, and $1=(b' \circ c' , d')$,
with $b,d \in A_X$, $b' ,d' \in A_{X'}$,
$c \in (\sqp \cup \{a\})^X$, 
$c' \in (\sqp \cup \{a'\})^{X'}$. 

Thus we have
\[1=1\circ 1= (b \circ c, d)\circ (b' \circ c', d')=
(b \circ c \circ \widetilde{b'} \circ \widetilde{c'}, d' \circ \widetilde{d})=
(b \circ \widetilde{b'} \circ \widetilde{c} \circ \widetilde{c'} , d' \circ \widetilde{d})
\]
But $\widetilde{c} \circ \widetilde{c'} \in ((\sqp \cup \{a\}) \circ (\sqp \cup \{a'\} ))^{X\prod X'}
\subseteq (\sqp \cup \{a\circ a' \})^{X \prod X'}$, and so
$1 \in (\sqp , a \circ a')_A$, 
hence $a \circ a' \not\in \sqp$.

\section{The Zariski topology}
\label{sec4.2}
The \emph{closed sets} in $spec(A)$ are defined to be 
the set of the form
\be
V(\sqa) = \{ \sqp \in spec(A) , \sqp \supseteq \sqa \},
\ee
with $\sqa \subseteq A_{[1]}$, 
which we may take to be an $h$-ideal $\sqa \in il(A)$.

We have

\noindent $
\begin{array}{ccl}
(4.2.2)& (i) & V(\sum\limits_i \sqa_i) = \bigcap\limits_i V(\sqa_i),\\
 &(ii) & V(\sqa \cdot \sqa') = V(\sqa) \cup V(\sqa'),\\
 &(iii) & V(0) = spec(A) \quad , \quad V(1) = \emptyset
\end{array}$

\setcounter{equation}{2}

This shows the sets $V(\sqa)$ define a topology on $spec(A)$,
the \emph{Zariski topology}.

\vspace{10pt}

For a subset $C \subseteq spec(A)$, we have the $h$-ideal,
\be
I(C) = \bigcap\limits_{\sqp \in I} \sqp
\ee
We have a Galois corrspondence,
\be
h\mbox{-}il(A) 
\begin{array}{c}
\stackrel{\textstyle V}{\longrightarrow}\\
\stackrel{\textstyle \longleftarrow}{I}
\end{array}
\{C \subseteq spec(A)\}
\ee
The maps $V$, $I$ are monotone
\be\begin{array}{l}
\sqa_1 \subseteq \sqa_2 \Rightarrow V(\sqa_1) \supseteq V(\sqa_2)\\
C_1 \subseteq C_2 \Rightarrow I(C_1) \supseteq I(C_2)
\end{array}
\ee
and we have
\be
\sqa \subseteq IV(\sqa) \quad , \quad 
C \subseteq VI(C)
\ee
It follows that we have
\be
V(\sqa) = VIV(\sqa) \quad \mbox{and} \quad
I(C) = IVI(C)
\ee
and the maps $V$, $I$ induce inverse bijections
\be
\label{eq4.2.8}
\{\sqa \in h\mbox{-}il(A) , \sqa = IV(\sqa)\}=
\hspace{2in}
\ee
\[=\{I(C) , C \subseteq spec(A)\} 
\stackrel{\textstyle \sim}{\longleftrightarrow}
\{ C \subseteq spec(A) , C=VI(C)\} =
\]
\[\hspace{2in}=\{ V(\sqa) , \sqa \in h\mbox{-}il(A) \}
\]

\noindent (4.2.9) {\bf Lemma:} For $\sqa \in h\mbox{-}il(A)$, we have
\[IV(\sqa) =
\{ a \in A_{[1]}, a^n \in \sqa \, \mbox{for some} \,
n > 0 \}
\stackrel{def}{=}
\sqrt{\sqa}
\]
\noindent {\bf Proof:} If $a \in \sqrt{\sqa}$, say $a^n \in \sqa$, 
then for all $\sqp \supseteq \sqa$,
$a \in \sqp $,
and so $\sqrt{\sqa} \subseteq \bigcap\limits_{\sqa \subseteq \sqp} 
\sqp = IV(\sqa)$.

Assume $a \not\in \sqrt{\sqa}$, so $a^n \not\in \sqa$ for all $n$.
An application of Zorn's lemma gives that there exists a 
maximal element $\sqp$ in the set

\setcounter{equation}{9}

\be
\label{eq4.2.10}
\{ \sqb \in h\mbox{-}il(A), \, \sqb \supseteq \sqa , \,
a^n \not\in \sqb \, \mbox{for all} \, n\}
\ee
We claim $\sqp$ is prime. If $x, x' \in A_{[1]}\setminus \sqp$, 
then the $h$-ideals $(\sqp,x)_A$,
$(\sqp, x')_A$ properly contain $\sqp$,
and by maximality of $\sqp$ in the set 
(\ref{eq4.2.10}),
we must have $a^n \in (\sqp, x)_A$,
$a^{n'} \in (\sqp, x')_A$,
for some $n$, $n'$.
We get 
\[
a^{n +n'}= a^n \circ a^{n'}=
(b \circ c , d) \circ (b' \circ c' , d' )=
(b \circ \widetilde{b'}\circ \widetilde{c} \circ \widetilde{c'} , 
 d' \circ \widetilde{d})
\]
with $b, d \in A_X$, $b', d' \in A_{X'}$, 
$c \in (\sqp \cup \{x \})^X$,
$c' \in (\sqp \cup \{x' \})^{X'}$.

But $\widetilde{c}\circ \widetilde{c'} \in \left(
(\sqp \cup \{x\})\circ(\sqp \cup \{x'\})
\right)^{X \prod X'} \subseteq
\left( \sqp \cup \{x \circ x' \}\right)^{X \prod X}$,
and since $a^{n+n'} \not\in \sqp$,
we must have $x \circ x' \not\in \sqp$;
and $\sqp$ is indeed prime.
Now $\sqp \supseteq \sqa$, and $a \not\in \sqp$, 
so $a \not\in \bigcap\limits_{\sqa \subseteq \sqp} \sqp =
IV(\sqa)$.

\vspace{10pt}

\noindent (4.2.11) {\bf Lemma:} For a subset $C \subseteq spec(A)$, 
$VI(C) = \bar{C}$ the closure of $C$.
\setcounter{equation}{11}

{\bf Proof:} We have $C \subseteq VI(C)$, 
and $VI(C)$ is closed. If $C \subseteq V(\sqa)$, 
where we may assume $\sqa = \sqrt{\sqa}$,
then $VI(C) \subseteq VIV(\sqa) = V(\sqa)$, and so 
\[
VI(C) = \bigcap\limits_{C \subseteq 
V(\sqa)} V(\sqa) = \bar{C}
\]

\vspace{10pt}

We can restate the bijection (\ref{eq4.2.8}),
as a bijection between the \emph{radical} $h$-ideals
and the closed subsets of $spec(A)$,
\be
\label{eq4.2.12}
\{\sqa \in h\mbox{-}il(A), \sqa = \sqrt{\sqa} \}
\stackrel{\textstyle \sim}{\longleftrightarrow}
\{C \subseteq spec(A) , C=\bar{C} \}
\ee

\section{Basic open sets}
A basis for the open sets of $spec(A)$ is given by the
\emph{basic open sets},
these are defined for $a \in A_{[1]}$ by
\be
D_a = spec(A) \setminus V(a) =
\{ \sqp \in spec(A) , a \not\in \sqp \}
\ee
We have,
\be
\begin{array}{l}
D_{a_1} \cap D_{a_2} = D_{a_1 \circ a_2}\\
D_1 = spec(A) \quad  , \quad D_0 = \emptyset
\end{array}
\ee

That every open set is the union of basic open sets, is shown by
\be
spec(A) \setminus V(\sqa) = \bigcup\limits_{a \in \sqa} D_a
\ee

\vspace{10pt}

Note that we have,
\be
D_a = spec(A) \Leftrightarrow a \circ A_{[1]} = \{a \}_A = (1)
\ee
\[
\Leftrightarrow  \, \mbox{there exists a (unique)} \,
a^{-1} \in A_{[1]} , a \circ a^{-1} =1
\]
We say that such $a$ is \emph{invertible}, and we let $A^*$ 
denote the set of invertible elements.
Note that $A^*$ is an abelian group (with involution), and
$A \mapsto A^*$ is a functor $\cG\cR \rightarrow Ab$ (=abelian groups).

\vspace{10pt}

Note that we have,
\be
D_a = \emptyset \Leftrightarrow 
a \in \bigcap\limits_{\sqp \in spec(A)} \sqp = \sqrt{0}
\ee
\[ 
\Leftrightarrow 
\, \mbox{there exists} \, n > 0 \, \mbox{with} \, a^n =0
\]
We say that such $a$ is \emph{nilpotent}.

\vspace{10pt}

\noindent (4.3.6) {\bf Lemma:} Let $\sqa = \sqrt{\sqa} 
\in h\mbox{-}il(A)$ 
be a radical $h$-ideal. Then
\[
V(\sqa) \, \mbox{is irreducible} \,
\Leftrightarrow 
\sqa \, \mbox{is prime}
\]
\setcounter{equation}{6}

{\bf Proof:} $(\Leftarrow)$: If $\sqa$ is prime, 
$V(\sqa)= VI\{\sqa\} = \overline{\{\sqa\}}$ 
is the closure of a point,
hence irreducible.

$(\Rightarrow)$: For any $a \in A_{[1]}$, we have
\be
\label{eq4.3.7}
V(\sqa) \cap D_a \neq \emptyset
\Leftrightarrow \exists \sqp \in spec(A) ,
\sqp \supseteq \sqa , \sqp \not\ni a  
\ee
\[
\Leftrightarrow 
a \not\in \bigcap\limits_{\sqa \subseteq \sqp} \sqp =
\sqrt{\sqa} = \sqa
\]
Hence for any basic open sets $D_a$, $D_b$, $a,b \in A_{[1]}$, 
we have
\[
V(\sqa) \cap D_a \neq \emptyset 
\, \mbox{and}\,
V(\sqa) \cap D_b \neq \emptyset
\Leftrightarrow 
a \not\in \sqa \, \mbox{and} \, b \not\in\sqa
\]
If $V(\sqa)$ is irreducible this implies
\[
\emptyset \neq V(\sqa) \cap D_a \cap D_b =
V(\sqa) \cap D_{a \circ b}
\Leftrightarrow 
a \circ b \not\in \sqa
\]

\vspace{10pt}

Thus the bijection (\ref{eq4.2.12}) induces a bijection
\be
\begin{array}{rl}
spec(A) \stackrel{\textstyle \sim}{\longleftrightarrow}&
\{ C \subseteq spec(A) , C= \bar{C} \, 
\mbox{closed and irreducible} \}\\
\sqp \mapsto & V(\sqp) = \widebar{\{ \sqp \}}
\end{array}
\ee
and the space $spec(A)$ is a Sober space (or a Zariski space):
every closed irreducible subset has a unique generic point.

\vspace{10pt}

\noindent (4.3.9) {\bf Proposition:} For $a \in A_{[1]}$, 
the basic open set $D_a$ is compact.
\setcounter{equation}{9}

In particular, $D_1 = spec(A)$ is compact.

{\bf Proof:} We have to show that in every covering of $D_a$ 
by basic open sets $D_{g_i}$,
there is always a finite subcovering. We have
\be
\label{eq4.3.10}
\begin{array}{ccl}
D_a \subseteq \bigcup\limits_{i}D_{g_i} &
\Leftrightarrow & 
V(a) \supseteq \bigcap\limits_i V(g_i) = V(\sum\limits_i
     g_i \circ A_{[1]})\\
&\Leftrightarrow & 
\sqrt{\{a\}_A} = IV(a) \subseteq IV\left(
\sum\limits_i g_i \circ A_{[1]}\right) = 
\sqrt{\sum\limits_i g_i \circ A_{[1]}}\\
&\Leftrightarrow &
\mbox{for some} \, n , a^n \in \sum\limits_i g_i \circ A_{[1]}\\
&\Leftrightarrow &
\mbox{for some} \, n , X \in \Fbb, b, d, \in A_X ,
a^n = (b \circ c, d) ,\\
&& \hspace{2cm} 
\mbox{with} \, c= (c^{(x)}) \in \left( 
\{g_i \} \right)^X
\end{array}
\ee
Thus $c^{(x)} = g_{i(x)}$, 
and going backwards in the above equivalences we get 
$D_a \subseteq \bigcup\limits_{x \in X} D_{g_{i(x)}}$,
 a finite subcovering.

\section{Functoriality}

For a homomorphism of general rings $\varphi \in \cG\cR (A,B)$,
the pull-back of a prime is a prime, and we have a map
\be
\begin{array}{cl}
\varphi^* = spec (\varphi): &spec(B) \rightarrow spec(A)\\
& \sqq \mapsto \varphi^* (\sqq) = \varphi^{-1}_{[1]} (\sqq)
\end{array}
\ee
The inverse image under $\varphi^*$ of a closed set is closed, we have
\be
\varphi^{*-1} (V_A(\sqa)) =
\{\sqq \in spec(B) , \varphi^{-1}_{[1]}(\sqq) \supseteq \sqa \}=
\ee
\[=
\{ \sqq \in spec(B) , \sqq \supseteq \varphi_{[1]}(\sqa) \}=
V_B(\varphi_{[1]}(\sqa))
\]
Also the inverse image under $\varphi^*$ of a basic open set is a basic open set, we have
\be
\varphi^{*-1}(D_a) =
\{\sqq \in spec(B) ,  \varphi^{-1}_{[1]}(\sqq) \not\ni a\}=
\ee
\[=
\{\sqq \in spec(B) , \varphi_{[1]}(a) \not\in \sqq\}=
D_{\varphi_{[1]}(a)}
\]
Thus the map $\varphi^*= spec(\varphi)$ is continuous, 
and we see that $spec$ is a contravariant
functor from $\cG\cR$ to the caegory $Top$,
whose objects are (compact, sober)
topological spaces, and continuous maps,
\be
spec: (\cG\cR)^{op} \rightarrow Top
\ee

\vspace{10pt}

\noindent (4.4.5) {\bf Lemma:} For $\varphi \in \cG\cR(A,B)$,
and for $\sqb \in h\mbox{-}il(B)$, we have
\[
V_A(\varphi_{[1]}^{-1}(\sqb)) =
\overline{\varphi^*(V_B(\sqb))}
\]
\setcounter{equation}{5}

{\bf Proof:} We may assume without loss of generality that 
$\sqb = \sqrt{\sqb}$ is radical 
(noting that 
$\sqrt{\varphi^{-1}_{[1]}(\sqb)} =\varphi^{-1}_{[1]}(\sqrt{\sqb})$).
Put $\sqa = I\varphi^*(V(\sqb))$, 
so that $V(\sqa) = \overline{\varphi^*(V(\sqb))}$ by (4.2.11).

We have for any $a \in A_{[1]}$,
\[\begin{array}{cl}
a \in \sqa&\Leftrightarrow 
a \in \sqp , \, \mbox{for every prime} \, \sqp \in \varphi^*(V(\sqb))\\
&\Leftrightarrow a \in \varphi^* (\sqq) = \varphi^{-1}_{[1]}(\sqq) ,
  \, \mbox{for every prime} \, \sqq \in V(\sqb)\\
&\Leftrightarrow
\varphi_{[1]}(a) \in \bigcap\limits_{\sqb \subseteq \sqq} \sqq =
       \sqrt{\sqb} = \sqb\\
&\Leftrightarrow a \in \varphi_{[1]}^{-1}(\sqb)
\end{array}
\]
Thus $\sqa = \varphi_{[1]}^{-1}(\sqb)$, and the lemma is proved.

\section{The stable spectrum}
In this section  assume $A \in \cG\cR^+$ is self-adjoint.

An application of Zorn's lemma, noting that $0 \in E\mbox{-}[1]\mbox{-}il(A)$,
and that for a chain $\sqa_i \in E\mbox{-}[1]\mbox{-}il(A)$ 
also the union is stable,
$\bigcup\sqa_i \in E\mbox{-}[1]\mbox{-}il(A)$, gives

\noindent (4.5.1) {\bf Proposition:} There exists maximal proper \emph{stable} $h$-ideals.
\setcounter{equation}{1}

\vspace{10pt}

More generally, for a stable $h$-ideal $\sqa \in E\mbox{-}[1]\mbox{-}il(A)$,
and for $f \in A_{[1]}$, such that $f^n \not\in \sqa$
for all $n$,
an application of Zorn's lemma gives that there exists a 
maximal element $\sqm$ in the set
\be
\{\sqb \in E\mbox{-}[1]\mbox{-}il(A) , \sqb \supseteq \sqa , f^n \not\in \sqb 
\, \mbox{for all} \, n \}
\ee
\noindent (4.5.3) {\bf Claim:} Such a maximal $\sqm$ is prime, 
$\sqm \in Spec(A)$.
\setcounter{equation}{3}

{\bf Proof:} For $x \in A_{[1]} \setminus \sqm$, 
the stable $h$-ideal generated by $\sqm$ and $x$,
$ZE\{m,x\}_A$,
properly contains $m$,
hence contains some power $f^n$.
Therefore, there exists a path, 
$f^n= a_1, a_2, \ldots , a_l = 0$, with 
$\{a_j, a_{j+1} \}$ of the form 
$\{(b \circ c , d ) , ( b \circ \bar{c} , d) \}$
with
\[ 
b,d \in A_{X_j \oplus Y_j \oplus Z_j} , \quad
c, \bar{c} \in \left(A_{[1]} \right)^{X_j \oplus Y_j \oplus Z_j} , \quad
c^{(x)} = \bar{c}^{(x)} \,\mbox{for} \,
x \in X_j , \]
\[
c^{(y)}, \bar{c}^{(y)} \in \sqm 
\, \mbox{for} \,
y \in Y_j , \quad
c^{(z)} = x , \bar{c}^{(z)}=0 \, \mbox{for} \,
z \in Z_j
\]
For $x' \in A_{[1]} \setminus \sqm$, 
there exists a similiar path 
$f^{n'}= a'_1, a'_2, \ldots , a'_l = 0$,
with $\{a'_j, a'_{j+1} \}$ of the form as above with 
$x$ replaced by $x'$.
Assume $\sqm$ is not prime and $x \cdot x' \in \sqm$.
The path 
$\{x' \cdot f^n , x' \cdot a_2 , \ldots , x' \cdot a_l =0 \}$ 
has the form 
$\{( b \circ \widetilde{x'}\circ c , d ) , ( b \circ \widetilde{x'} 
   \circ \bar{c} , d ) \}$ $=$ 
$\{ ( b \circ e , d ) , ( b \circ \bar{e} , d ) \}$
with $e^{(x)} = \bar{e}^{(x)}$ for 
$x \in X_j$,
and $e^{(y)}, \bar{e}^{(y)} \in \sqm$ for
$y \in Y_j \coprod Z_j$
(because $x' \cdot x \in \sqm$),
hence $x'\cdot f^n$
is in $ZE(\sqm)=\sqm$
(because $\sqm$ is stable).

The path 
$\{f^{n +n'}, f^n \cdot a'_n , \ldots , f^n\cdot a'_{l'}=0\}$
has the form\\
$\{(b' \circ \widetilde{f^n} \circ c' , d' ) ,
( b' \circ \widetilde{f^n} \circ \bar{c'} , d') \}$,
with
$(\widetilde{f^n} \circ c' )^{(x)} =
(\widetilde{f^n} \cdot \bar{c'})^{(x)}$
for $x \in X'_j$,
and $(\widetilde{f^n} \circ c')^{(y)},
(\widetilde{f^n} \circ \bar{c'})^{(y)}\in \sqm$
for $y \in Y'_j \coprod Z'_j$
(because $f^n \cdot x' \in \sqm$),
hence $f^{n +n'}$
is in $ZE(\sqm) =\sqm$,
contradiction.

\vspace{10pt}

\noindent (4.5.4) {\bf Corollary:} For a \emph{stable} $h$-ideal 
$\sqa \in E\mbox{-}[1]\mbox{-}il(A)$, we have
\[
\sqrt{\sqa} = \bigcap\limits_{\sqa \subseteq \sqp \in ESpec(A)} \sqp
\]
the intersection taken over all the \emph{stable} primes 
containing $\sqa$.
\setcounter{equation}{4}

\vspace{10pt}

We have $ESpec(A) \subseteq spec(A)$,
the subset of stable primes,
and we give it the induced Zariski topology.
The closed sets are
\be
\widetilde{V}(\sqa) = V(\sqa) \cap ESpec(A), \quad \sqa \subseteq A_{[1]}
\ee
Note that this set depends only on $ZE(\sqa)_{[1]}$, 
and we may take $\sqa \in E\mbox{-}[1]\mbox{-}il(A)$.
For a subset $C \subseteq ESpec(A)$, 
$I(C) = \bigcap\limits_{\sqp \in C} \sqp$ is stable,
$I(C) \in E\mbox{-}[1]\mbox{-}il(A)$.
The formulas, and lemmas, of (\ref{sec4.2})
carry over to the stable setting.
We have a Galois correspondence
\be
 E\mbox{-}[1]\mbox{-}il(A) 
\begin{array}{c}
\stackrel{\textstyle \widetilde{V}}{\longrightarrow}\\
\stackrel{\textstyle \longleftarrow}{I}
\end{array}
\{C \subseteq ESpec(A) \}
\ee
inducing a bijection.
\be
\label{eq4.5.7}
\{ \sqa \in E\mbox{-}[1]\mbox{-}il(A) , \, \sqa = \sqrt{\sqa} \}
\stackrel{\textstyle \sim}{\longleftrightarrow}
\{C \subseteq ESpec(A) , \, C =\bar{C} \, \mbox{closed} \}
\ee
A basis for the topology of $ESpec(A)$ is given by the 
\emph{basic stable open sets},
\be
\widetilde{D}_a = ESpec(A) \setminus \widetilde{V}(a) =
\{ \sqp \in ESpec(A) , a \not\in \sqp \}
=D_a \bigcap ESpec(A)
\ee
Note that we have
\be
\widetilde{D}_a = ESpec(A) \Leftrightarrow ZE(a \circ A)_{[1]} = A_{[1]}
\Leftrightarrow
\ee
there exists a path $1=a_1, a_2, \ldots , a_l=0$, 
with $\{a_j , a_{j+1} \}$ of the form 
$\{(b \circ c,d), (b \circ \bar{c}, d) \}$, 
with $b,d \in A_{X_j \oplus Y_j}$,
$c,\bar{c} \in \left( A_{[1]}\right)^{X_j \oplus Y_j}$,
and $c^{(x)} = \bar{c}^{(x)}$ for $x \in X_j$,
$c^{(y)}=a$, $\bar{c}^{(y)}=0$ for $y \in Y_j$.

We say that $a$ is  a \emph{unit}, and we let
$A^{(1)}$ denote
the set of units;
$A^* \subseteq A^{(1)} \subseteq A_{[1]}$.

Note that on the other hand,
\be
\widetilde{D}_a = \emptyset \Leftrightarrow
a \in \bigcap\limits_{\sqp \in ESpec(A)} \sqp = \sqrt{0}
\Leftrightarrow
a \, \mbox{is nilpotent} \, 
\Leftrightarrow
D_a = \emptyset
\ee
Thus we have

\noindent (4.5.11) {\bf Proposition:}
\setcounter{equation}{11}
$ESpec(A)$ is dense in $spec(A)$.

\vspace{10pt}

For a radical stable ideal,
$\sqa = \sqrt{\sqa} \in E\mbox{-}[1]\mbox{-}il(A)$, 
and for $a \in A_{[1]}$, we have cf. (\ref{eq4.3.7}),
\be
\widetilde{V}(\sqa) \bigcap \widetilde{D}_a \neq \emptyset
\Leftrightarrow
\exists \sqp \in ESpec(A) , \sqp \supseteq \sqa, \sqp \not\ni a 
\ee
\[
\Leftrightarrow
a \not\in \bigcap\limits_{\sqa \subseteq \sqp} \sqp = 
\sqrt{\sqa}=\sqa
\]
Hence (4.3.6) goes through:
$\widetilde{V}(\sqa)$ is irreducible if and only if $\sqa$ is prime.
Thus the bijection (\ref{eq4.5.7}) induces a bijection
\be
ESpec(A) \stackrel{\textstyle \sim}{\longleftrightarrow}
\{ C \subseteq ESpec(A) , C \, \mbox{closed and irreducible} \}
\ee
i.e. the space $ESpec(A)$ is sober too.

\vspace{10pt}

\noindent (4.5.14) {\bf Proposition:}
\setcounter{equation}{14}
The basic stable open sets $\widetilde{D}_a = D_a \bigcap ESpec(A)$
are compact.

In particular, $\widetilde{D}_1 = ESpec(A)$ is compact.

{\bf Proof:} We have, cf. (\ref{eq4.3.10}),
\be
\begin{array}{cl}
\widetilde{D}_a \subseteq \bigcup\limits_i \widetilde{D}_{g_i}&
\Leftrightarrow I\widetilde{V} (a) \subseteq I\widetilde{V}(\{g_i\})\\
&\Leftrightarrow a \in \sqrt{(ZE\{g_i\})_{[1]}}\\
&\Leftrightarrow \mbox{for some} \, n, a^n \in (ZE\{g_i\})_{[1]} 
\end{array}
\ee
$\Leftrightarrow$ there exists a path 
$a^n = a_1, a_2, \ldots , a_l =0$,
with $\{a_j , a_{j+1}\}$ of the form
$\{(b \circ c, d ), ( b \circ \bar{c}, d) \}$
with $b,d \in A_{X_j \oplus Y_j}$,
$c, \bar{c} \in \left( A_{[1]}\right)^{X_j \oplus Y_j}$,
and $c^{(x)} = \bar{c}^{(x)}$ for $x \in X_j$,
$c^{(y)}= g_{i(y)}$, 
$\bar{c}^{(y)}=0$ for $y \in Y_j$.

Going backwards in the above equivalences we have a finite subcovering,
\[\widetilde{D}_a \subseteq \bigcup\limits_{0 < j<l}
     \bigcup\limits_{y \in Y_j} \widetilde{D}_{g_{i(y)}}
\]

\vspace{10pt}

For a homomorphism of self-adjoint generalized rings 
$\varphi \in \cG \cR^+ (A, A')$,
the pull back of a stable $[1]$-ideal is stable,
\be
\varphi^* : E\mbox{-}[1]\mbox{-}il(A') \longrightarrow E\mbox{-}[1]\mbox{-}il(A),
\ee
and in particular, we get a continuous map
\be
ESpec(\varphi) = \varphi^*: ESpec(A') \longrightarrow ESpec(A).
\ee
The analogs of (4.4.2) , (4.4.3), (4.4.4) , (4.4.5)
hold:
\be
\varphi^{*-1}(\widetilde{V}_A(\sqa)) =
\widetilde{V}_{A'}(\varphi_{[1]}(\sqa)),
\ee
\be
\varphi^{*-1}(\widetilde{D}_a)=
\widetilde{D}_{\varphi_{[1]}(a)},
\ee
\be
\widetilde{V}_A(\varphi^{-1}_{[1]}(\sqb) ) =
\overline{\varphi^*(\widetilde{V}_{A'}(\sqb)},
\, \mbox{for} \, \sqb \in E\mbox{-}[1]\mbox{-}il(A'),
\ee
and $ESpec$ is a functor from $\cG\cR^+$ to 
(compact, sober) topological 
spaces and continuous maps
\be
ESpec: (\cG\cR^+)^{op} \rightarrow Top
\ee

\vspace{10pt}

For $A=\cG(B)$, $B$ a commutative ring, we have identifications,
\be
ESpec(A) = spec(A) =spec(B)
\ee

\chapter{Localization and sheaves}
\section{Localization}
For $A \in \cG\cR$, a subset 
$S \subseteq A_{[1]}$ is called \emph{multiplicative}
if
\be
1 \in S , \ \mbox{and} \
S \circ S \subseteq S , \ \mbox{and} \
S^t = S 
\ee
For such $S \subseteq A_{[1]}$, and for $X \in \Fbb$, 
we let $(S^{-1}A)_X = (A_X \times S) /\hspace{-5pt}\approx$
denote the equivalence classes of $A_X \times S$ with respect to the equivalence relation defined by
\be
(a_1 , s_1 ) \approx (a_2, s_2) \ \mbox{if and only if} \
s \circ s_2 \circ a_1 = s \circ s_1 \circ a_2
\ \mbox{for some} \ s \in S
\ee
We write $a/s$ for the equivalence class $(a,s) /\hspace{-5pt}\approx$.
Note that by taking "common denominator"
we can write any element 
$a = (a^{(y)}/s_y) \in (S^{-1}A)_f = \prod\limits_{y \in Y}
(S^{-1}A)_{f^{-1}(y)}$,
in the form
\[a = (\bar{a}^{(y)}/s), \left( \mbox{take} \ s = \prod s_y ,\
\bar{a}^{(y)}= (\prod\limits_{y' \neq y}s_{y'}) \cdot a^{(y)})\right)
\]
For $f \in Set_{\bullet}(X,Y)$, 
$g \in Set_{\bullet}(Y,Z)$,
we have well-defined operations of multiplication and contraction,
independent of the choice of representatives,
\be
\circ: (S^{-1}A)_g \times (S^{-1}A)_f \longrightarrow 
       (S^{-1}A)_{g \circ f} , \quad
 a/s_1 \circ b/s_2 = a \circ b / s_1 \circ s_2
\ee
\be
(\, ,\, ):(S^{-1}A)_{g \circ f } \times (S^{-1}A)_f \longrightarrow 
       (S^{-1}A)_{g} , \quad 
(a/s_1 , b/s_2) = (a , b) / s_1 \circ (s_2)^t
\ee

It is straightforward to check that these operations 
satisfy the axioms of a generalized ring
(commutative if $A$ is commutative).
The cannonical homomorphism
\[
\phi:A \rightarrow S^{-1}A , \quad \phi(a) = a/1
\]
satisfies the universal property
\be\cG\cR (S^{-1}A,B) =
\{\varphi \in \cG\cR (A,B) , \varphi(S) \subseteq B^*\}
\ee
\[\begin{array}{r}\widetilde{\varphi} \mapsto \widetilde{\varphi} \circ \phi\\
\varphi(s)^{-1} \circ \varphi (a) = \widetilde{\varphi} (a/s)
\leftmapsto \varphi
\end{array}
\]

\noindent (5.1.7) {\bf Example:} For $s \in A_{[1]}$,
take $S = \{ s^n \cdot (s^m)^t , \ m,n \geq 0 \}$.
We write $A_{s}$ for $S^{-1}A$,
and $\phi_{s} \in \cG\cR (A,A_s)$ satisfy 
\[\cG\cR (A_{s},B) = \{ \varphi \in \cG\cR(A,B) , \varphi (s) \in B^* \}
\]

\vspace{10pt}

\noindent (5.1.8) {\bf Example:}
\setcounter{equation}{8}
For $\sqp = \sqp^t \in spec(A)$,
take $S_{\sqp} = A_{[1]}\setminus \sqp$.
We write $A_{\sqp}$ for $S^{-1}_{\sqp} A$,
and $\phi_{\sqp} \in \cG\cR (A,A_{\sqp})$ satisfy 
\[\cG\cR (A_{\sqp},B) = \{ \varphi \in \cG\cR(A,B) , 
    \varphi (A_{[1]}\setminus \sqp) \subseteq B^* \}
\]

\section{Localization and $h$-ideals}
For an $h$-ideal $\sqa \in h\mbox{-}il(A)$, we let
\be
S^{-1}\sqa = \{ a/s \in (S^{-1}A)_{[1]} , s \in S , a \in \sqa \}
\ee
By using common denominator,
we see that $S^{-1}\sqa$ is an $h$-ideal of $S^{-1}A$,
$S^{-1}\sqa \in h\mbox{-}il(S^{-1}A)$:

For $b/s_1, d /s_2  \in (S^{-1}A)_X$, and for 
$a_x/s_x \in S^{-1}\sqa$, $x \in X$, we have
\be
(b/s_1 \circ (a_x/s_x), d/s_2)=
(b \circ (s'_x \circ a_x), d)/s_1 \circ (s_2)^t \circ
    \prod\limits_{x \in X}s_x \in S^{-1}\sqa ,
\ee
\[
 \mbox{with} \
    s'_x =\prod\limits_{x' \neq x} s_{x'}
\]

We have, therefore, the Galois correspondence
\be
h\mbox{-}il(A)
\begin{array}{c}
\stackrel{\textstyle S^{-1}}{\longrightarrow}\\
\stackrel{\textstyle \longleftarrow}{\phi^*}
\end{array}
h\mbox{-}il(S^{-1}A)
\ee
For $\sqb \in h\mbox{-}il(S^{-1}A)$, we have
\be
S^{-1}\phi^* \sqb = \sqb
\ee
Indeed, for an element $a/s \in \sqb$, 
we have $a/1 \in \sqb$, or $a \in \phi^*\sqb$,
and $a/s \in S^{-1}\phi^*\sqb$;
hence $\sqb \subseteq S^{-1} \phi^* \sqb$,
and the reverse inclusion is clear.
We have immediately from the definitions, 
for $\sqa \in h\mbox{-}il(A)$,
\be
\phi^* S^{-1}\sqa =
\{ a \in A_{[1]}, \ \mbox{there exists} \ s \in S
\ \mbox{with} \ s \circ a \in \sqa \} = 
\bigcup\limits_{s \in S} (\sqa:s)
\ee
In particular,
\be
S^{-1}\sqa = (1) 
\Leftrightarrow
\sqa \cap S \neq \emptyset
\ee
We say that $\sqa \in h\mbox{-}il(A)$ is $S$-\emph{saturated}
if $\phi^* S^{-1}\sqa = \sqa$, that is if 
\be
\mbox{for all}\ s \in S, 
a \in A_{[1]}: s \circ a \in \sqa 
\Rightarrow a \in \sqa
\ee
We get that $S^{-1}$ and $\phi^*$ induce inverse bijections,
\be
\{ \sqa \in h\mbox{-}il(A) , \quad \sqa \ \mbox{is} \ 
S\mbox{-saturated} \}
\stackrel{\textstyle \sim}{\longleftrightarrow}
h\mbox{-}il(S^{-1}A)
\ee
For an $S$-saturated $h$-ideal 
$\sqa \in h\mbox{-}il(A)$,
let $\pi_{\sqa}: A \twoheadrightarrow A/\sqa = A/E(\sqa)$
be the canonical homomorphism, and let
$\bar{S} = \pi_{\sqa}(S) \subseteq (A/\sqa)_{[1]}$, 
then we have canonical isomorphism
\be
\label{eq5.2.10}
\bar{S}^{-1}(A/\sqa) \cong S^{-1}A/S^{-1}\sqa
\ee
Note that for a prime $\sqp \in spec(A)$,
$\sqp$ is $S$-saturated if and only if $\sqp \cap S = \emptyset$, 
and
in this case $S^{-1}\sqp$ is a prime,
$S^{-1}\sqp \in spec(S^{-1}A)$.
Note that for a prime $\sqq \in spec(S^{-1}A)$,
$\phi^*(\sqq)$ is always an $S$-saturated prime.
We get the bijection,
\be
\label{eq5.2.11}
\{\sqp \in spec(A) , \sqp \cap S = \emptyset \}
\stackrel{\textstyle \sim}{\longleftrightarrow}
spec(S^{-1}A)
\ee
This is a homeomorphism for the Zariski topology.

\vspace{10pt}

From now on let $A \in \cG\cR^+$ be self-adjoint.

For $s \in A_{[1]}$, the homeomorphism (\ref{eq5.2.11})
gives for example (5.1.7),
\be
\phi_{s}^*: spec(A_{s})
\stackrel{\textstyle \sim}{\longrightarrow}
D_{s} \subseteq spec(A)
\ee

For a prime $\sqp \in spec(A)$, the homeomorphism (\ref{eq5.2.11})
with $S = S_{\sqp}$,
c.f. (5.1.8), reads
\be
\phi_{\sqp}^*: spec(A_{\sqp})
\stackrel{\textstyle \sim}{\longrightarrow}
\{ \sqq \in spec(A) , \sqq \subseteq \sqp \}
\ee
The generalized ring $A_{\sqp}$
is a \emph{local-generalized-ring}
in the sense that it has a unique maximal $h$-ideal 
$m_{\sqp} = S_{\sqp}^{-1} \sqp = (A_{\sqp})_{[1]} \setminus A_{\sqp}^{*}$.
The \emph{residue field at $\sqp$} is defined by
\be
\mathbb{F}_{\sqp} = A_{\sqp}/m_{\sqp} = A_{\sqp}/E(m_{\sqp}).
\ee
We have the canonical homomorphism 
$\pi_{\sqp}: A \rightarrow A/\sqp$,
and putting $\bar{S}_{\sqp} = \pi_{\sqp}(S_{\sqp})$,
we have (\ref{eq5.2.10})
\be
\mathbb{F}_{\sqp} = \bar{S}_{\sqp}^{-1}(A/\sqp)
\ee
The square diagram
\be
\label{eq5.2.17}
%\mbox{PICTURE}
\ee
\begin{diagram}
A & \rTo^{\phi_{\sqp}}& A_{\sqp}\\
\dOnto^{\pi_{\sqp}}& & \dOnto_{\pi_{m_{\sqp}}}\\
A/\sqp& \rTo^{\phi_0}& {\mathbb F}_{\sqp}
\end{diagram}

is cartesian,
\be
\mathbb{F}_{\sqp} = (A/\sqp) \bigotimes\limits_{A} A_{\sqp}
\ee
\[\cG \cR^+ (\mathbb{F}_{\sqp}, B) =
\left\{\varphi \in \cG\cR^+(A,B) , \
       \varphi(\sqp) \equiv 0 , \
       \varphi(A_{[1]}\setminus \sqp) \subseteq B^* 
\right\}
\]

\vspace{10pt}

\noindent (5.2.17) {\bf Proposition:}
\setcounter{equation}{17}
We have for $\sqp \in spec(A)$,
\[\mathbb{F}_{\sqp} \neq 0 
\Leftrightarrow
\sqp \ \mbox{is stable:} \
\sqp \in E[1]\mbox{-}il(A)
\]

{\bf Proof:} We have
\[\mathbb{F}_{\sqp}=0 
\Leftrightarrow
1 =0 \ \mbox{in} \
\mathbb{F}_{\sqp,[1]}
\]
$\Leftrightarrow$
there exists a path in 
$(S_{\sqp}^{-1}A)_{[1]}$,
$1=a_1/s_1 , a_2/s_2 , \ldots , a_l/s_l = 0$,
with
$\{a_j/s_j , a_{j+1}/s_{j+1} \}$
of the form
$\{ (b \circ c, d) , (b \circ \bar{c},d) \}$
with
$b,d \in (S^{-1}A)_{X_j\oplus Y_j}$,
$c, \bar{c} \in (S^{-1}A_{[1]})^{X_j\oplus Y_j}$,
and
$c^{(x)}=\bar{c}^{(x)}$
for
$x \in X_j$,
$c^{(y)}, \bar{c}^{(y)} \in m_{\sqp}$
for
$y \in Y_j$

$\Leftrightarrow$
there exists a path in
$A_{[1]}$,
$s = a_1,a_2, \ldots , a_l =0$,
with
$s \in S_{\sqp}$,
and
$\{a_j, a_{j+1}\}$
of the form
$\{ (b \circ c, d) , (b \circ \bar{c},d) \}$
with
$b,d \in A_{X_j\oplus Y_j}$,
$c, \bar{c} \in (A_{[1]})^{X_j\oplus Y_j}$,
and
$c^{(x)}=\bar{c}^{(x)}$
for
$x \in X_j$,
$c^{(y)}, \bar{c}^{(y)} \in \sqp$
for
$y \in Y_j$

$\Leftrightarrow$
$\sqp \subsetneqq \left(ZE\sqp \right)_{[1]}$

\vspace{10pt}

For a homomorphism of self-adjoint generalized rings
$\varphi \in \cG\cR^+ (A,B)$,
and for $\sqq \in spec(B)$
with
$\sqp = \varphi^*(\sqq) \in spec(A)$,
the square diagram (\ref{eq5.2.17})
is functorial,
and we have a commutative cube diagram

\noindent (5.2.18) 
\setcounter{equation}{18}

\begin{diagram}
 & & A &\rTo^{\hspace{0.5cm}\phi_{\sqp}} & & & A_{\sqp}\\
 &\ldTo^{\varphi} &\dTo& &   & \ldTo^{\varphi_{\sqp}}&\dTo\\
B&\rTo^{\hspace{2cm}\phi_{\sqq}}& & & B_{\sqq}& & \\
\dTo &              & & &\dTo     & & \\
     &              &A/\sqp& \rTo& & &{\mathbb F}_{\sqp}\\
  &\ldTo &                 &   &  & \ldTo& \\
B/\sqq& \rTo &             &  &{\mathbb F}_{\sqq}   & &
\end{diagram}

Note that the homomorphism $\varphi_{\sqp} \in \cG\cR^+(A_{\sqp},
 B_{\sqq} )$
is a \emph{local-homomorphism}
in the sense that
\be
m_{\sqp}= \varphi^{-1}_{\sqp}(m_{\sqq}) , 
\ \mbox{or equivalently} \
\varphi_{\sqp}(m_{\sqp}) \subseteq m_{\sqq}
\ee

\noindent (5.2.20) {\bf Definition:} We let $\cL\cG\cR^+$
denote the subcategory of $\cG\cR^+$
with objects the local generalized rings,
and with maps
\[
\cL\cG\cR^+(A,B) = \{ \varphi \in \cG\cR^+(A,B), \varphi^*(m_B) = m_A\}
\]

\section{The structure sheaf $\cO_A$}

\noindent(5.3.1) {\bf Definiton:}
\setcounter{equation}{1}
For a self-adjoint $A \in \cG\cR^+$,
$U \subseteq spec(A)$ open,
$X \in \Fbb$,
we denote by $\cO_A(U)_X$
the set of \emph{sections}
\[f:U \rightarrow \coprod\limits_{\sqp \in U} (A_{\sqp})_X , \quad
f(\sqp) \in (A_{\sqp})_X
\]
such that $f$ is \emph{locally a fraction}:

\noindent (5.3.2) for all  $\sqp \in U$,
there exists open $U_{\sqp} \subseteq U$, $\sqp \in U_{\sqp}$, 
and there exist $a \in A_X$, $s \in A_{[1]}\setminus 
     \bigcup\limits_{\sqq \in U_{\sqp}}\sqq$,
such that for all $\sqq \in U_{\sqp}$, 
$f(\sqq) \equiv a/s \in (A_{\sqq})_X$.
\setcounter{equation}{2}

\vspace{10pt}

Note that $\cO_A(U)$ is a (self-adjoint) generalized ring, 
and for $U' \subseteq U$ restriction gives a homomorphism of generalized rings
\be
\cO_A(U) \rightarrow \cO_A(U') , \quad
f \mapsto f|_{U'}
\ee
Thus $\cO_A$ is a pre-sheaf of generalized rings over $spec(A)$,
and by the local nature of the condition (5.3.2)
it is clear that it is a \emph{sheaf of generalized rings},
i.e. for $X \in \Fbb$,
$U \mapsto \cO_A(U)_X$ is a sheaf.
It is also clear that the \emph{stalks} are given by
\be
\cO_{A,\sqp}= \lim\limits_{\stackrel{\longrightarrow}{\sqp \in U}}
\cO_A(U) 
\stackrel{\textstyle \sim}{\longrightarrow} A_{\sqp}
\ee
\[\left(f/\approx\right) \mapsto f(\sqp)
\]

\noindent (5.3.5) {\bf Theorem:} For $s \in A_{[1]}$, 
we have a canonical isomorphism
\setcounter{equation}{5}
\[\Psi : A_{s} \stackrel{\textstyle \sim}{\longrightarrow}
\cO_A(D_{s}) , \quad
\Psi(a/s^n) =
\{f(\sqp) \equiv a/s^n \}
\]
In particular for $s=1$,
\[
A \stackrel{\textstyle \sim}{\rightarrow} \cO_A(spec(A))
\]

{\bf Proof:} The  map $\Psi$ which takes 
$a/s^n \in A_{s}$ to the constant section $f$
with $f(\sqp) \equiv a/s^n$ for all $\sqp \in D_{s}$,
is clearly well-defined, and is a homomorphism of generalized rings.

$\Psi$ \underline{is injective}: Assume $\Psi(a_1/s^{n_1})= \Psi(a_2 / s^{n_2})$,
and let $\sqa = ann_A(s^{n_2}\circ a_1 ,  s^{n_1}\circ a_2) \in 
h\mbox{-}il(A)$,
cf. (\ref{eq3.5.9}).
We have,

\be
\begin{array}{l}
a_1/s^{n_1} = a_2/s^{n_2} \ \mbox{in} \ A_{\sqp}
    \ \mbox{for all} \ \sqp \in D_{s} \\
\Rightarrow s_{\sqp} \circ s^{n_2} \circ a_1 =
            s_{\sqp} \circ s^{n_1} \circ a_2 
             \ \mbox{with} \ s_{\sqp} \in A_{[1]}\setminus \sqp 
             \ \mbox{for} \ \sqp \in D_{s}\\
\Rightarrow \sqa \not\subseteq \sqp \ \mbox{for} \
            \sqp \in D_{s} \\
\Rightarrow V(\sqa) \cap D_{s} = \emptyset \\
\Rightarrow V(\sqa) \subseteq V(s) \\
\Rightarrow s \in IV(\sqa) = \sqrt{\sqa} \\
\Rightarrow s^n \in \sqa \ \mbox{for some} \ n\\
\Rightarrow s^{n + n_2} \circ a_1 = s^{n+n_1} \circ a_2 \\
\Rightarrow a_1/s^{n_1} = a_2/s^{n_2} \ \mbox{in} \ A_{s}
\end{array}
\ee
$\Psi$ \underline{is surjective:} Fix 
  $f \in \cO_A(D_{s})_X$.
Since $D_{s}$ is compact (4.3.9),
we can cover $D_{s}$ by a finite collection of basic open sets, 
$D_{s} = D_{g_1}\cup \ldots \cup D_{g_N}$,
such that on $D_{g_i}$ the section $f$ is constant,
\[
f(\sqp) = a_i/s_i \ \mbox{for} \
          \sqp \in D_{g_i}, \ \mbox{with} \
          D_{s_i} \supseteq D_{g_i}
\]
We have $V(s_i) \subseteq V(g_i)$, hence 
$g_i \in IV(s_i) = \sqrt{s_i}$, hence for some $n_i$,
and some $c_i \in A_{[1]}$,
$g_i^{n_i}= c_i \circ s_i$.
Thus our section $f$ is given on $D_{g_i}$ by 
$a_i/s_i = c_i \circ a_i/g_i^{n_i}$.
Noting that $D_{g_i}= D_{g_i^{n_i}}$,
we may replace $g_i^{n_i}$ by $g_i$,
and replace $c_i \circ a_i$ by $a_i$,
and we have
\[f(\sqp) = a_i/g_i \ \mbox{for} \
          \sqp \in D_{g_i},
\]
On the set $D_{g_ig_i} = D_{g_i}\cap D_{g_j}$, $i \neq j$,
our section $f$ is given by both $a_i/g_i$ and $a_j/g_j$.
By the injectivity of $\Psi$, we have
\[
a_i/g_i = a_j/g_j \ \mbox{in} \ A_{g_ig_2}
\]
Thus for some $n$ we have
\[
(g_ig_j)^n \circ g_j \circ a_i = (g_ig_j)^n\circ g_i \circ a_j
\]
By finiteness we may assume one $n$ works for all $i,j \leq N$.
Replacing $g_i^n\circ a_i$ by $a_i$,
and replacing $g_i^{n+1}$ by $g_i$,
we may assume $f \equiv a_i/g_i$ on $D_{g_i}$,
and
\be
g_j \circ a_i = g_i \circ a_j \ \mbox{for all} \ i,j
\ee
We have $D_{s} \subseteq \bigcup D_{g_i}$, 
hence by (\ref{eq4.3.10})
we have 
\[s^M = (b\circ c , d)\]
with $b,d \in A_Y$,
$c =(c^{(y)}) \in \left(A_{[1]}\right)^Y$ with 
$c^{(y)} = g_{i(y)}$, $i(y):Y \rightarrow \{1, \ldots, N\}$.
Define $a \in A_X$ by 
\be
a = (b \circ e , \tilde{d}), \ \mbox{with} \
    \tilde{d} \in A_{\pi_X} = (A_Y)^X ,
    \tilde{d}^{(x)} \equiv d ,
\ee\[
    e \in A_{\pi_Y} = (A_X)^Y , e^{(y)} = a_{i(y)}
\]
Here the cartesian diagram is

\begin{diagram}
 & &X \otimes Y& & \\
 &\ldTo^{\pi_X}& &\rdTo^{\pi_Y}& \\
X & & & &Y\\
 & \rdTo & & \ldTo & \\
 & & [1] & & 
\end{diagram}

We have for $j = 1, \ldots, N$
\be
\begin{array}{ll}
g_j \circ a &= g_j \circ (b \circ e , \tilde{d}) = 
               (b \circ (g_j e) , \tilde{d})= 
               (b \circ (g_j \circ a_{i(y)} ), \tilde{d} )\\
  & = (b \circ (g_{i(y)} \circ a_j ), \tilde{d} ) =
   (b \circ c \circ \tilde{a}_j , \tilde{d} ) =
(b \circ c , d) \circ a_j = s^M \circ a_j
\end{array}
\ee
Thus we have in $A_{s}$,
$a_j/g_j = a/s^M$ for all $j$,
and our section $f$ is constant 
$f=\Psi(a/s^M)$, and
$\Psi$ is surjective.

\chapter{Schemes}

\section{Locally generalized ringed spaces}
All generalized rings are assumed to be self-adjoint
(and can be assumed to be commutative).

\noindent (6.1.1) {\bf Definition:}
For a topological space $\cX$, we let $\cG\cR^+/\cX$ 
denote the category
of sheaves of generalized rings over $\cX$.
Its objects are pre-sheaves $\cO$ of (self-adjoint) 
generalized rings,
i.e. functors $U \mapsto \cO(U): \cC_{\cX} \rightarrow \cG\cR$,
(with $\cC_{\cX}$ the category of open sets of $\cX$, 
with $\cC_{\cX}(U,U')=\{j^U_{U'}\}$
for $U' \subseteq U$,
otherwise
$\cC_{\cX}(U,U') = \emptyset$),
such that for all $X \in \Fbb$,
$U \mapsto \cO(U)_X$
is a sheaf.
The maps $\cG\cR/\cX(\cO,\cO')$
are natural transformations of functors $\varphi = \{\varphi(U)\}$,
$\varphi(U) \in \cG\cR(\cO(U), \cO'(U))$.

\vspace{10pt}

\noindent (6.1.2) {\bf Definition:} 
We denote by $\cG\cR\cS$
the category of generalized ringed spaces.
Its objects are pairs $(\cX, \cO_{\cX})$,
with $\cX \in Top$,
and $\cO_{\cX} \in \cG\cR^+/\cX$.
The maps $f \in \cG\cR\cS (\cX,\cY)$ 
are pairs of a continuous function $f \in Top(\cX, \cY)$, 
and a map of sheaves of generalized rings over $\cY$,
$f^{\sharp} \in \cG\cR/\cY (\cO_{\cY}, f_*\cO_{\cX})$;
explicitly, for all open subsets $U \subseteq \cY$,
we have a homomorphism of generalized rings 
\[f^{\sharp}_U = \{f^{\sharp}_{U,X}\} \in \cG\cR (\cO_{\cY}(U), 
\cO_{\cX}(f^{-1}U))\] 
and these homomorphisms are compatible with restictions:
for $U' \subseteq U \subseteq \cY$ open, and for 
$a \in \cO_{\cY}(U)_X$,
we have 
$f^{\sharp}_{U,X}(a) |_{f^{-1}(U')} = f^{\sharp}_{U',X} (a|_{U'})$
in $\cO_{\cX} (f^{-1} U')_X$.

\vspace{10pt}

\noindent (6.1.3) {\bf Remark:} For a continuous map 
$f \in Top(\cX, \cY)$,
we have a pair of adjoint functors
\[
\cG\cR^+/\cX
\begin{array}{c}
\stackrel{\textstyle f^*}{\curvearrowleft}\\
\stackrel{\textstyle \curvearrowbotright}{f_*}
\end{array}
\cG\cR^+/\cY
\]
For sheaves of genralized rings $\cO_{\cX} \in \cG\cR^+/\cX$,
$\cO_{\cY} \in \cG\cR^+/\cY$,
we have
\[
f_*\cO_{\cX}(U)= \cO_{\cX}(f^{-1}U) \ , \ U \subseteq \cY \ 
\mbox{open;}
\]
\[
f^*\cO_{\cY} (U)_X = \ 
\mbox{sheaf associated to the pre-sheaf} \]
\[
U \mapsto 
\lim\limits_{\begin{array}{c}
\stackrel{\textstyle \longrightarrow}{V \subseteq \cY \ \mbox{open}}\\
f(U) \subseteq V
\end{array}}
\cO_{\cY}(V)_X;
\]
and we have adjunction,
\[
\cG\cR^+/\cY (\cO_{\cY} , f_*\cO_{\cX}) = 
\cG\cR^+/\cX(f^*\cO_{\cY} ,\cO_{\cX})
\]

\noindent (6.1.4) {\bf Remark:}
For a map of generalized ringed spaces 
$f \in {\cal GRS(X,Y)}$,
and for a point $x \in \cX$,
we get the induced homomorphism on \emph{stalks}
$f^{\sharp}_x \in \cG\cR(\cO_{\cY,f(x)}, \cO_{\cX,x})$,
via
\[
f^{\sharp}_x: \cO_{\cY,f(x)}=
\lim\limits_{\begin{array}{c}
\stackrel{\textstyle \longrightarrow}{V \subseteq \cY \ \mbox{open}}\\
f(x) \in V
\end{array}}
\hspace{-1cm}
\cO_{\cY}(V)
\xrightarrow{\lim\limits_{\rightarrow} f_V^{\sharp}}
\lim\limits_{\begin{array}{c}
\stackrel{\textstyle \longrightarrow}{V \subseteq \cY \ \mbox{open}}\\
x \in f^{-1}V
\end{array}}
\hspace{-1cm}
\cO_{\cX}(f^{-1}V)
\rightarrow\]\[
\hspace{2in}
\rightarrow
\lim\limits_{\begin{array}{c}
\stackrel{\textstyle \longrightarrow}{U \subseteq \cX \ \mbox{open}}\\
x \in U
\end{array}}
\hspace{-1cm}
\cO_{\cX}(U)=\cO_{\cX,x}
\]

\noindent (6.1.5) {\bf Definition:}
We let $\cal LGRS \subseteq GRS$
denote the subcategory of $\cal GRS$
of \emph{locally generalized ringed spaces}.
Its objects are the objects $(\cX, \cO_{\cX})\in \cG\cR\cS$
such that for all points $x \in \cX$ the stalk
$\cO_{\cX,x} \in \cL\cG\cR$
is a \emph{local} self-adjoint generalized ring.
The maps 
$f \in \cL\cG\cR\cS (\cX, \cY)$
are the maps
$(f, f^{\sharp})\in \cG\cR\cS(\cX, \cY)$,
such that for all points $x \in \cX$,
the induced homomorphism on stalks is a \emph{local}
homomorphism,
(5.2.20),
$f_x^{\sharp} \in \cL\cG\cR (\cO_{\cY,f(x)},
\cO_{\cX,x})$.

\vspace{10pt}

\noindent (6.1.6) {\bf Theorem:}
The functor of global sections
\[
\Gamma: \cL\cG\cR\cS \rightarrow (\cG\cR^+)^{op} , \quad
\Gamma(\cX, \cO_{\cX})=\cO_{\cX}(\cX) , \]\[
\Gamma(f,f^{\sharp}) = f_{\cY}^{\sharp} \ \mbox{for} \
f \in \cL\cG\cR\cS (\cX, \cY)
\]
and the spectra functor
\[
spec:(\cG\cR^+)^{op} \rightarrow {\cal LGRS} ,\quad
spec(A) = (spec(A), \cO_A) ,\]\[
spec(\varphi) = \varphi^* \ \mbox{for} \
\varphi \in \cG\cR^+(A,B)
\]
are an adjoint pair:
\[
{\cal LGRS}(\cX , spec(A)) =
\cG\cR^+ (A,\cO_{\cX}(\cX) ) \]\[ \mbox{functorially in} \
\cX \in \cL\cG\cR\cS  , \quad
A \in \cG\cR^+
\]

{\bf Proof:}
For a point $x\in \cX$ we have the canonical 
homomorphism of taking the stalk at $x$ of a global section,
$\phi_x \in \cG\cR(\cO_{\cX}(\cX) , \cO_{\cX,x})$.
Since $\cO_{\cX,x}$ is local with a unique maximal $h$-ideal 
$m_{\cX,x}$,
we get by pullback a prime 
$\sqp_x = \phi_x^*(m_{\cX,x}) \in spec(\cO_{\cX}(\cX))$.
Thus we have a canonincal map
\[
\sqp:\cX \rightarrow  spec(\cO_{\cX}(\cX)) ,\quad
x \mapsto \sqp_x
\]

The map $\sqp$ is continuous:
For a global section
$g \in \cO_{\cX}(\cX)_{[1]}$,
we have the basic open set 
$D_g \subseteq spec(\cO_{\cX}(\cX))$,
and
\setcounter{equation}{6}
\be
\sqp^{-1}(D_g) = 
\{x \in \cX , \ \sqp_x \in D_g \}=
\{x \in \cX , \phi_x(g) \not\in m_{\cX,x}\}
\ee

This set is open in $\cX$, because if $\phi_x(g) \not\in m_{\cX,x}$
we have in $\cO_{\cX,x}$ some $v_x$
with $v_x \circ \phi_x(g) =1$,
hence there is an open set $U \subseteq \cX$,
with $x \in U$,
and an element $v \in\cO_{\cX}(U)_{[1]}$ with $v \circ g|_{U}=1$,
and for all $x'\in U$,
$v_{x'}\circ \phi_{x'}(g) =1$,
and
$\phi_{x'}(g) \not\in m_{\cX,x'}$.
This shows $\sqp$ is continuous.
The uniqueness of the inverse $v_x= \phi_x(g)^{-1}$ for 
$x \in \sqp^{-1}(D_g)$ 
shows we have a well defined inverse 
$v = (g|_{\sqp^{-1}(D_g)})^{-1} \in \cO_{\cX}(\sqp^{-1}(D_g))_{[1]}$.
Thus we have a homomorphism of generalized rings
\be
\label{eq6.1.8}
\sqp^{\sharp}_{D_g}: \cO_{\cX}(\cX)_g 
\rightarrow
\cO_{\cX}(\sqp^{-1}(D_g)) , \quad
a/g^n \mapsto v^n \circ \left(a|_{\sqp^{-1}(D_g)}\right)
\ee
The collection of homomorphisms
$\{ \sqp^{\sharp}_{D_g} , g \in \cO_{\cX}(\cX)_{[1]} \}$,
are compatible with restrictions, 
and the sheaf property gives homomorphisms
$\sqp_U^{\sharp} \in \cG\cR \left(\cO_{spec \cO_{\cX}(\cX)}(U) ,
\cO_{\cX}(\sqp^{-1}(U))\right)$.
Thus we have a map of generalized ringed spaces 
$\sqp=  (\sqp, \sqp^{\sharp}) \in \cG\cR\cS (\cX , 
spec(\cO_{\cX}(\cX)))$.
For a point $x \in \cX$, we can take the direct limit of 
$\sqp^{\sharp}_{D_g}$,
(\ref{eq6.1.8}),
over all global sections
$g \in \cO_{\cX}(\cX)_{[1]}$
with
$\phi_x(g) \not\in m_{\cX,x}$,
and we get a local homomorphism
$\sqp_x^{\sharp} \in \cL\cG\cR 
(\cO_{\cX}(\cX)_{\sqp_x}, \cO_{\cX,x})$.
This shows $\sqp$ is a map of locally-ringed spaces,
$\sqp \in {\cal LGRS}(\cX, spec(\cO_{\cX}(\cX)))$.

Given a homomorphism of generalized rings 
$\varphi \in \cG\cR^+(A, \cO_{\cX}(\cX))$,
we get the map in $\cal LGRS$
\be
(spec \varphi)\circ \sqp: \cX \rightarrow spec (\cO_{\cX}(\cX)) 
\rightarrow spec(A)
\ee
Given a map of locally ringed spaces 
$f = (f, f^{\sharp}) \in {\cal LGRS}(\cX, spec(A))$,
we get a homomorphism in $\cG\cR^+$,
\be
\Gamma (f) = f^{\sharp}_{spec(A)}:A = \cO_A
(spec(A)) \rightarrow \cO_{\cX}(\cX)
\ee
These correspondences give the functorial bijection of
(6.1.6),
we need only show they are inverses of each other.
First for $\varphi \in \cG\cR^+(A, \cO_{\cX}(\cX))$,
we have
\be
\Gamma(spec(\varphi)\circ \sqp ) =
\Gamma(\sqp) \circ \Gamma (spec(\varphi)) =
id_{\cO_{\cX}(\cX)} \circ \varphi = \varphi
\ee
Fix a map $f=(f, f^{\sharp}) \in {\cal LGRS}(\cX, spec(A))$.
For a point $x \in \cX$,
we have a commutative square in $\cG\cR^+$
\be
\begin{array}{ccc}
A=\cO_A(spec(A))& \stackrel{\textstyle \Gamma (f)}{\longrightarrow}&
\cO_{\cX}(\cX)\\
\phi_{f(x)} \downarrow & &\downarrow \phi_x\\
A_{f(x)}= \cO_{A,f(x)}&
 \stackrel{\textstyle f^{\sharp}_x}{\longrightarrow}&
\cO_{\cX,x} 
\end{array}
\ee
Since the homomorphism $f_x^{\sharp}$ is assumed to be local, we get
\be
\Gamma (f)^{-1}(\sqp_x) = \Gamma (f)^{-1}(\phi_{x}^{-1}(
  m_{\cX,x})) = \phi_{f(x)}^{-1}
(f_x^{\sharp -1} (m_{\cX,x}))=
\ee\[=
\phi_{f(x)}^{-1}(m_{A_{f(x)}})=f(x)
\]
This shows that 
$(spec\ \Gamma(f)) \circ \sqp = f$
as continuous maps.

For an element $s \in A_{[1]}$,
we have the commutative square in $\cG\cR^+$,
\be
\begin{array}{ccl}
A=\cO_A(spec(A))&
\stackrel{\textstyle \Gamma(f)}{\longrightarrow}&
\cO_X(X)\\
\downarrow& &\downarrow\\
A_{s}=\cO_A(D_{s})&
\stackrel{\textstyle f_{D_{s}}^{\sharp}}{\longrightarrow}&
\cO_{\cX}(f^{-1}(D_{s}))=\cO_{\cX}(D_{f^{\sharp}(s)})
\end{array}
\ee
Thus for $a/s^n \in A_{s}$,
we must have
\be
f_{D_{s}}^{\sharp}(a/s^n)=
\left(\Gamma(f)(s^n)|_{f^{-1}(D_{s})}\right)^{-1}
\circ
\left(\Gamma(f)(a)\right)|_{f^{-1}(D_{s})}
\ee
\[
=\sqp^{\sharp}_{D_{f^{\sharp}(s)}}\circ \Gamma(f) (a/s^n)
\]
This shows that 
$f = (spec \ \Gamma(f) ) \circ \sqp$
also as maps of generalized-ringed spaces.

\section{Schemes}
\noindent (6.2.1) {\bf Definition:}
An object $\cX =(\cX, \cO_{\cX}) \in \cL\cG\cR\cS$
will be called a 
\emph{Grothendieck-generalized-scheme}
if it is locally isomorphic to 
$spec(A)$'s:
there exists a covering of $\cX$ by open sets $U_i$,
$X = \bigcup\limits_{i}U_i$,  
such that the canonical maps are isomorphisms
\[
\sqp:(U_i, \cO_{\cX}|_{U_i})
\stackrel{\sim}{\longrightarrow}
spec (\cO_{\cX}(U_i))
\]
We let $\cG\cG\cS$ denote the full sub-category of 
$\cL\cG\cR\cS$,
with objects the Grothendieck-genralized-schemes.

\vspace{10pt}

\noindent (6.2.2) {\bf Open subschemes:}
Note that for $\cX \in \cG\cG\cS$,
and for an open set $U \subseteq \cX$,
we have the \emph{open subscheme} of $\cX$
given by $(U, \cO_{\cX}|_{U})$.
That this is again a scheme,
$(U, \cO_{\cX}|_{U})\in \cG\cG\cS$,
follows from the existence of affine basis for the 
Zariski topology on $spec(A)$, 
$A \in \cG\cR^+$,
namely $(D_{s}, \cO_A|_{D_{s}}) \cong spec(A_{s})$
for $s \in A_{[1]}$.

\vspace{10pt}

\noindent (6.2.3) {\bf Gluing shemes:}
The local nature of the definition of Grothendieck-generalized-scheme implies that $\cG\cG\cS$
admits gluing:

Given $\cX_i \in \cG\cG\cS$,
and open subsets $U_{ij} \subseteq \cX_i$,
and maps $\varphi_{ij} \in \cG\cG\cS (U_{ij}, U_{ij})$,
satisfying the consistency conditions
\begin{itemize}
\item[$(i)$] $U_{ii}=\cX_i$, and $\varphi_{ii}=id_{\cX_i}$,
\item[$(ii)$] $\varphi_{ij}(U_{ij}\cap U_{ik})= U_{ji}\cap U_{jk}$,
and $\varphi_{jk} \circ \varphi_{ij}= \varphi_{ik}$
on $U_{ij}\cap U_{ik}$,
\end{itemize}
 there exists $\cX \in\cG\cG\cS $, and maps 
$\varphi_i \in \cG\cG\cS (\cX_i , \cX)$ such that
\begin{itemize}
\item[$(i)$] $\varphi_i$ is an isomorphism of $\cX_i$ onto an 
    open subset $\varphi_i(\cX_i) \subseteq \cX$
\item[$(ii)$] $\cX = \bigcup\limits_i \varphi_i(\cX_i)$
\item[$(iii)$]$\varphi_i(\cX_i) \cap \varphi_j(\cX_j) =
               \varphi_i(U_{ij})$,
and $\varphi_j\circ\varphi_{ij}=\varphi_i$ on $U_{ij}$.
\end{itemize}

\vspace{10pt}

\noindent (6.2.4) {\bf Ordinary Schemes:}
For an ordinary scheme 
$(\cX, \cO_{\cX})$,
with $\cO_{\cX}$
 a sheaf of commutative rings, there is a covering by open sets
$\cX= \bigcup\limits_i U_i$,
with $(U_i,\cO_{\cX}|_{U_i}) \cong spec(A_i)$,
the ordinary spectrum of the commutative ring
$A_i = \cO_{\cX}(U_i)$.
We then have Grothendieck-generalized schemes
$\cX_i = spec(\cG(A_i)) = (U_i, \cG(\cO_{\cX})|_{U_i})$.
These can be glued along $U_{ij}= U_i \cap U_j$,
to a Grothendieck-generalized scheme denoted by
$\cG(\cX) = (\cX, \cO_{\cG(\cX)} = \cG(\cO_{\cX}))$.
It is just the underlying topological space $\cX$
with the sheaf of generalized rings $\cG(\cO_{\cX})$
associated to the sheaf of commutative rings $\cO_{\cX}$ 
via the functor $\cG : Ring \rightarrow \cG\cR^+$.
(\ref{eq2.2.3}).
Denoting by $\cR\cS$ the category of 
(ordinary, commutative) ringed spaces,
the functor 
$\cG$ applied to a sheaf of commutative rings $\cO$,
gives a sheaf of self-adjoint generalized rings $\cG(\cO)$,
and we have a functor $\cG:\cR\cS \rightarrow \cG\cR\cS$.
Denoting by $\cL\cR\cS$ (resp. by $\cS$) the
category of
locally-(commutative)-ringed
spaces (resp. the full subcategory of ordinary schemes),
the fact that $\cG$ is fully-faithful implies that we have
full-embeddings of categories.

\begin{diagram}
\cL \cR \cS & \rInto^{\cG} & {\cal LGRS}\\
\uInto & & \uInto\\
\cS & \rInto^{\cG}& \cG\cG \cS\\
\uInto^{spec}& & \uInto_{spec}\\
Ring^{op}& \rInto^{\cG}& (\cG\cR^+)^{op}
\end{diagram}

\vspace{10pt}

\noindent (6.2.5) The affine line over $\mathbb F$ is given by,
cf. (\ref{eq2.6.18}),
\[
{\mathbb A}^1 = spec \, \Delta_+^{[1]} = spec \, {\mathbb F} 
[z^{\mathbb N}]
\]
We have $F[z^{\mathbb N}]_{[1]} = z^{\mathbb N} \cup \{ 0 \}$;
$(0)$ is a prime, the generic point of ${\mathbb A}^1$;
and $(z)$ is a prime, the closed point of 
${\mathbb A}^ 1 = \{ (0), (z) \}$.

\vspace{10pt}

\noindent (6.2.6) The multiplicative group over $\mathbb F$ 
is given by
\[
{\mathbb G}_m = spec \, {\mathbb F}[z^{\mathbb Z}]=
\{(0)\} \subseteq {\mathbb A}^1
\]
For a (self-adjoint) generalized ring $A$,
\[
\cG\cR ({\mathbb F} [z^{\mathbb Z}], A ) =
A^* =
\{a \in A_{[1]}, \ \mbox{there is} \ a^{-1} \in A_{[1]} , 
a \circ a^{-1}=1 \}
\]

\noindent (6.2.7) The projective line over $\mathbb F$ 
is obtained by 
gluing two affine lines along
${\mathbb G}_m$
\[
{\mathbb P}^1 = spec  \, {\mathbb F}[z^{\mathbb N}]
\prod\limits_{spec  \, {\mathbb F}[z^{\mathbb Z}]}
spec \, {\mathbb F}[(z^{-1})^{\mathbb N}] =
\{ m_1 , m_0, m_{\infty}\}
\]
It has a generic point $m_1 = (0)$,
and two closed points $m_0 = (z)$, $m_{\infty}= (z^{-1})$.

Interchanging $z$ and $z^{-1}$ we get an involutive automorphism

\vspace{10pt}

\noindent (6.2.8) $\displaystyle{I: {\mathbb P}^1
\xrightarrow{\sim} {\mathbb P}^1}$, \quad
$I \circ I = id_{\mathbb P'}$

interchanging $m_0$ and $m_{\infty}$.

Every rational number $f \in {\mathbb Q}^*$,
defines a geometric map 
$f_{\mathbb Z} \in \cG\cG\cS (spec \, \cG({\mathbb Z}), {\mathbb P'})$.

If $f = \pm 1$ this is given by the constant map
\[
\mathbb F [z^{\mathbb Z}] \twoheadrightarrow {\mathbb F}[\pm 1]
\subseteq \cG (\mathbb Z) , 
z \mapsto f = \pm 1
\]
If $f \neq \pm 1$,
let $N_0 = \prod\limits_{\nu_p (f) >0} p$,
$N_{\infty} = \prod\limits_{\nu_p (f) <0} p$,
then
\setcounter{equation}{8}
\be
spec \, \cG (\mathbb Z) = spec \,  \cG (\mathbb Z [\frac{1}{N_0}])
\coprod\limits_{spec \, \cG(\mathbb Z [\frac{1}{N_0N_{\infty}}])}
spec \, \cG (\mathbb Z [\frac{1}{N_{\infty}}])
\ee
and the geometric map $f_{\mathbb Z}$ is given by the spec-map associated
to the homomorphisms:
\be
\begin{array}{ccc}
\mathbb F [z^{\mathbb N}] & \longrightarrow & 
       \cG(\mathbb Z [\frac{1}{N_{\infty}}])\\
| \bigcap & & | \bigcap\\
\mathbb F [z^{\mathbb N}] & \longrightarrow & 
       \cG(\mathbb Z [\frac{1}{N_0 \cdot N_{\infty}}])\\
| \bigcup & & | \bigcup\\
\mathbb F [(z^{-1})^{\mathbb N}] & \longrightarrow & 
       \cG(\mathbb Z [\frac{1}{N_{0}}])
\end{array}
\ee
\[
\begin{array}{c}
z \mapsto f\\
z^{-1} \mapsto f^{-1}
\end{array}
\]

\section{Projective limits}

\noindent (6.3.1) The category of locally generalized ringed spaces
$\cal LGRS$
admits directed inverse limits. For a partially ordered set 
$J$,
which is directed
(for $j_1, j_2 \in J$, have $j \in J$ with $j \geq j_1$, 
$j \geq j_2$)
and for a functor $\cX:J\rightarrow \cL\cG\cR\cS$,
$J \ni j \mapsto \cX_j$,
$j_1 \geq j_2 \mapsto \pi_{j_2}^{j_1}\in \cL\cG\cR\cS (\cX_{j_1} ,
 \cX_{j_2})$,
we have the inverse limit
$\lim\limits_{\stackrel{\longleftarrow}{J}}\cX \in {\cal LGRS}$.
The underlying topological space of 
$\lim\limits_{\stackrel{\longleftarrow}{J}}\cX$ 
is the inverse limit of the sets $\cX_j$,
with basis for the topology given by the sets
$\pi_j^{-1}(U_j)$,
with $U_j \subseteq \cX_j$
open, and where 
$\pi_j : \lim\limits_{\stackrel{\longleftarrow}{j \in J}}
\cX_j \rightarrow \cX_j$
denote the projection.
The sheaf of generalized rings 
$\cO_{\lim\limits_{\leftarrow}\cX}$ over 
$\lim\limits_{\stackrel{\longleftarrow}{J}}\cX_j$,
is the sheaf associated to the pre-sheaf
$U \mapsto \lim\limits_{\stackrel{\longrightarrow}{J}}
\pi_j^* \cO_{\cX_j}(U)$.
For a point 
$x =(x_j) \in \lim\limits_{\longleftarrow}\cX_j$,
the stalk $\cO_{\lim\limits_{\longleftarrow}\cX,x}$
is the direct limit of the local-generalized-rings 
$\cO_{\cX_j, x_j}$,
and hence is local, and
$(\lim\limits_{\longleftarrow}\cX_j,
\cO_{\lim\limits_{\longleftarrow}\cX}) 
\in {\cal LGRS}$.
An alternative explicit description of the sections 
$s \in \cO_{\lim\limits_{\longleftarrow}\cX}(U)$,
for $U \subseteq \lim\limits_{\longleftarrow}\cX_j$ 
open,
are as maps
\setcounter{equation}{1}
\be
s:U \rightarrow \coprod\limits_{x \in U}
    \cO_{\lim\limits_{\longleftarrow}\cX,x},
\ \mbox{with} \
s(x) \in \cO_{\lim\limits_{\longleftarrow}\cX,x}
\ee
such that for all $x \in U$,
there exists $j \in J$, and open subset $U_j \subseteq \cX_j$,
with $x \in \pi_j^{-1}(U_j) \subseteq U$ 
and there exists a section $s_j \in \cO_{\cX_j}(U_j)$,
such that for all $y \in \pi_j^{-1}(U_j)$,
we have 
$s(y) = \pi_j^{\sharp}(s_j) |_{y}$.

We have the universal property
\be
\label{eq6.3.3}
{\cal LGRS}(Z, \lim\limits_{\stackrel{\longleftarrow}{J}}\cX)=
\lim\limits_{\stackrel{\longleftarrow}{j \in J}} {\cal LGRS} (Z, \cX_j)
\ee
Note that if 
$\cX_j = spec(A_j)$ are affine generalized schemes, then the inverse limit
\be
\lim\limits_{\stackrel{\longleftarrow}{J}}(spec(A_j)) =
spec(\lim\limits_{\stackrel{\longrightarrow}{J}}A_j)
\ee
is the affine generalized scheme associated to 
$\lim\limits_{\stackrel{\longrightarrow}{J}}A_j$
the direct limit of the $A_j$'s computed in $\cal GR^+$.
(Hence in $Set_0$, cf. (\ref{eq2.7.3})).

Note on the other hand that $\cal GGR$ is not closed 
under directed inverse limits
(just as in the "classical" counterparts, the category 
$\cal LRS$
(resp. $Ring$) is closed under directed inverse (resp. direct) limits , 
while the category $\cS$ of schemes is not closed under directed inverse limits).

\vspace{10pt}

\noindent (6.3.5) {\bf Defintion:}
\setcounter{equation}{5}
The category of generalized schemes $\cal GS$
is the category of pro-objects of the category of 
Grothendieck-generalized schemes,
$\cG\cS = pro\mbox{-}\cG\cG\cS$.

Thus the objects of $\cal GS$
are inverse systems 
$\cX=(\{\cX_j\}_{j \in J}, \{\pi_{j_2}^{j_1}\}_{j_1 \geq j_2})$,
where $J$ is a directed partially ordered set,
$\cX_j \in {\cal GGS}$ for $j \in J$,
and $\pi_{j_2}^{j_1}\in {\cal GGS}(\cX_{j_1}, \cX_{j_2})$ 
for $j_1 \geq j_2$, $j_1, j_2 \in J$, 
with $\pi_{j}^{j} = id_{\cX_j}$,
and $\pi_{j_3}^{j_2} \circ \pi_{j_2}^{j_1}
=\pi_{j_3}^{j_1}$ for $j_1 \geq j_2 \geq j_3$.
The maps from such an object to another object
$\cY = (\{\cY_i\}_{i \in I}, \{\pi_{i_2}^{i_1}\}_{i_1 \geq i_2})$ 
are given by
\be
\cG\cS (\cX, \cY )=
\lim\limits_{\stackrel{\longleftarrow}{I}}
\lim\limits_{\stackrel{\longrightarrow}{J}}
\cG\cG\cS (\cX_j , \cY_i )
\ee
i.e. the maps $\varphi \in \cG\cS(\cX, \cY)$ 
are a collection of maps
$\varphi_i^j \in \cG\cG\cS (\cX_j, \cY_i)$ defined for all 
$i \in I$,
and for $j \geq \tau(i)$ sufficiently large 
(depending on $i$),
these maps satisfy:
\begin{itemize}
\item[(6.3.6.1)]
for all $i \in I$, and for $j_1 \geq j_2$ sufficiently large in $J$:
\[\varphi_i^{j_1} =\varphi_i^{j_2} \circ \pi_{j_2}^{j_1}\]
\item[(6.3.6.2)] for all $i_1 \geq i_2$ in $I$,
and for $j \in J$ sufficiently large:
\[\pi_{i_2}^{i_1} \circ \varphi_{i_1}^{j} = \varphi_{i_2}^{j}\]
\end{itemize}
 The maps $\varphi = \{\varphi_i^j \}_{j \geq \tau(i)}$,
and $\tilde{\varphi} = \{\tilde{\varphi}_i^j\}_{j \geq 
       \tilde{\tau}(i)}$,
are considered equivalent if
\begin{itemize}
\item[(6.3.6.3)] for all $i \in I$, and for $j \in J$
sufficiently large:
\[\varphi_i^j = \tilde{\varphi}_i^j
\]
\end{itemize}
The composition of $\varphi = \{\varphi_i^j\}_{j \geq \tau(i)} \in
\cG\cS(\cX, \cY )$,
with $\psi = \{\psi^i_k \}_{i \geq \sigma(k)} \in
\cG\cS(\cY, \cZ)$,
is given by
$\psi \circ \varphi = \{ \psi^i_k \circ \varphi_i^j \}_{j \geq \tau(\sigma(k))} \in \cG\cS(\cX, \cZ)$.

There is a canonical map (which in general is not injective or surjective, but see \cite{G}).
\be
\label{eq6.3.7}
\lim\limits_{\stackrel{\longrightarrow}{J}} {\cal LGRS} 
(\cX_j , \cY_i) \longrightarrow {\cal LGRS}
(\lim\limits_{\stackrel{\longleftarrow}{J}} \cX_j, \cY_i)
\ee
By the universal property (\ref{eq6.3.3}) we have bijection
\be
\label{eq6.3.8}
\lim\limits_{\stackrel{\longleftarrow}{I}} {\cal LGRS}
(\lim\limits_{\stackrel{\longleftarrow}{J}} \cX_j, \cY_i)=
{\cal LGRS} 
(\lim\limits_{\stackrel{\longleftarrow}{J}} \cX_j ,
\lim\limits_{\stackrel{\longleftarrow}{I}} \cY_i) 
\ee
Composing (\ref{eq6.3.7}) and (\ref{eq6.3.8})
we obtain a map
\be
\cL: \lim\limits_{\stackrel{\longleftarrow}{I}}
     \lim\limits_{\stackrel{\longrightarrow}{J}}
{\cal LGRS} (\cX_j , \cY_i) \longrightarrow {\cal LGRS}
(\lim\limits_{\stackrel{\longleftarrow}{J}} \cX_j ,
\lim\limits_{\stackrel{\longleftarrow}{I}} \cY_i)
\ee
Thus we have a functor
\be
\cL: {\cal GS} \longrightarrow {\cal LGRS} , \quad
\cL(\{\cX_j\}_{j \in J}) = 
\lim\limits_{\stackrel{\longleftarrow}{J}} \cX_j 
\ee
We view the category $\cG\cG\cS$ as a full subcategory of 
$\cG\cS$
(consisting of the objects
$\cX = \{ \cX_j\}_{j \in J}$,
with indexing set $J$ reduced to a singleton).

\section{The compactified $\overline{spec \mathbb{Z}}$}
\label{sec6.4}
We denote by $\eta$ the \emph{real prime} of $\mathbb{Q}$,
so $|\ |: \mathbb{Q} \rightarrow [0, \infty )$
is the usual (nonarchimedian) absolute value, 
and we let $\cO_{\eta}$ denote the associated generalized ring
(\ref{sec2.3}),
$\cO_{\eta} \subseteq \cG(\mathbb{Q})$.
For a square-free integer $N \geq 2$, 
we have the sub-generalized-ring
\be
A_N = \cG(\mathbb{Z}[\frac{1}{N}]) \cap \cO_{\eta}
\subseteq \cG(\mathbb{Q})
\ee
The localization of $A_N$ with respect to 
$\frac{1}{N} \in A_{N,[1]}$
gives 
$(A_N)_{\frac{1}{N}}= \cG(\mathbb{Z}[\frac{1}{N}])$,
so the inclusion $j_N: A_N \hookrightarrow \cG(\mathbb{Z}[\frac{1}{N}])$
gives the basic open set
\be
j_N^*:spec(\mathbb{Z}[\frac{1}{N}])=
      spec \, \cG(\mathbb{Z}[\frac{1}{N}])
\stackrel{\textstyle \sim}{\longrightarrow}
D_{\frac{1}{N}} \subseteq spec(A_N)
\ee
The inclusion 
$i_N:A_N \hookrightarrow \cO_{\eta}$,
gives the real prime $\eta_N \in spec \,(A_N)$,
\be
\eta_N = i^*_N (\sqm_{\eta}) , \quad
(\eta_N)_X = \{ a = (a_x) \in (\mathbb{Z}[\frac{1}{N}])^X ,
||a||^2 = \sum\limits_{x \in X}|a_x|^2 < 1 \}
\ee
Note that $\eta_N$ is the unique maximal $h$-ideal of $A_N$,
and $A_N$ is a local generalized ring.
Let $\cX_N$ denote the Grothendieck generalized scheme obtained by gluing
$spec \, (A_N)$ with $spec \,  \cG(\mathbb{Z})$ along the common (basic) open set 
$spec(\cG(\mathbb{Z}[\frac{1}{N}]))$, cf. (6.2.3).
The open sets of $\cX_N$ are the open sets 
$U_{M} = spec(\mathbb{Z}[\frac{1}{M}]) \subseteq spec(\mathbb{Z})$,
(and $\cO_{\cX_N}(U_M)= \cG(\mathbb{Z}[\frac{1}{M}])$),
as well as the sets $\{\eta_N\} \cup U_M $,
with $M$ dividing $N$
(and $\cO_{\cX_N}(\{\eta_N\} \cup U_M)=A_M$, $M|N$).
For $ N_2$ dividing $N_1$,
we have a map $\pi_{N_2}^{N_1}\in \cG\cG\cS (\cX_{N_1}, \cX_{N_2})$
induced by the inclusions
$A_{N_2} \hookrightarrow A_{N_1}$,
and $\cG(\mathbb{Z}[\frac{1}{N_2}]) \hookrightarrow 
     \cG(\mathbb{Z}[\frac{1}{N_1}])$.
Note that $\pi_{N_2}^{N_1}$ is a bijection on points,
and that moreover,\\ 
$(\pi_{N_2}^{N_1})_* \cO_{\cX_{N_1}}= \cO_{\cX_{N_2}}$ and 
$(\pi_{N_2}^{N_1})^{\sharp}$ is the identity map of 
$\cO_{\cX_{N_2}}$.
But there are more open sets in $\cX_{N_1}$ 
then there are in $\cX_{N_2}$.
The compactified 
$\overline{spec \,  \mathbb{Z}}$
is the object of $\cG\cS = pro \cG\cG\cS$ given by 
$(\{\cX_N\}, \{\pi_{N_2}^{N_1}\}_{N_2|N_1})$,
\be
\overline{spec \, \mathbb{Z}} =
\left\{
\cX_N = spec(A_N) \coprod\limits_{spec \, \cG(\mathbb{Z}[\frac{1}{N}])}
        spec \, \cG(\mathbb{Z})
\right\}_{N \geq 2 \ \mbox{square free}}
\ee
Note that the associated locally-generalized-ring space
\be
\cX= \cL(\overline{spec \mathbb{Z}})=
\lim\limits_{\stackrel{\longleftarrow}{N}} \cX_N \in {\cal LGRS}
\ee
has underlying topological space
$\cX = \{\eta\} \coprod spec(\mathbb{Z})$, with open sets 
$U_M = spec(\mathbb{Z}[\frac{1}{M}])$
(and $\cO_{\cX}(U_M) = \cG(\mathbb{Z}[\frac{1}{M}])$), 
as well as the sets $\{\eta\} \coprod U_M$,
with no restrictions on $M$
and $\cO_{\cX}(\{\eta\} \coprod U_M) = A_M$
for $M \geq 2$,
while the global sections are 
$\cO_{\cX}(\cX)= \Fbb[\pm 1]$.

The stalks of $\cO_{\cX}$ are given by
\be
\begin{array}{lll}
\cO_{\cX, p} = \cG(\mathbb{Z}_{(p)}) , & p \in spec(\mathbb{Z}) , &
\mathbb{Z}_{(p)}= \{ \frac{m}{n} \in \mathbb{Q} , p \not\mid n\},\\
\cO_{\cX, \eta} = \cO_{\eta} & &
\end{array}
\ee
Similarly for a number field $K$,
with ring of integers $\cO_K$, and with real primes $\eta_i$,
$i = 1, \ldots , \gamma= \gamma_{\mathbb{R}} + \gamma_{\mathbb{C}}$,
we have the sub-generalized-ring of $\cG(K)$ given by
\be
A_{N,i}= \cG(\cO_K[\frac{1}{N}]) \cap \cO_{K,\eta_i}
\subseteq \cG(K)
\ee
Let $\cX_N$ be the Grothendieck generalized scheme obtained by gluing
$\{spec (A_{N,i})\}_{i \leq \gamma }$ and
$\{spec (\cG(\cO_K))\}$ along the common (basic) open set 
$spec \left(\cG(\cO_K[\frac{1}{N}])\right)$.
For $N_2|N_1$,
we have $\pi_{N_2}^{N_1} \in \cG\cG\cR(\cX_{N_1}, \cX_{N_2})$
induced by the inclusions $A_{N_2,i} \hookrightarrow A_{N_1,i}$.
We get the compactified $\overline{spec(\cO_K)}$,
it is the object of $\cG\cS$ given by the $\cX_N$'s and 
$\pi_{N_2}^{N_1}$'s.
The space
\be
\cX_K = \cL(\overline{spec(\cO_K)})= 
\lim\limits_{\stackrel{\longleftarrow}{N}}\cX_N \in {\cal LGRS}
\ee
has for points the sets 
$spec(\cO_K) \coprod \{\eta_i\}_{i \leq \gamma}$,
and for open subsets the sets $U \coprod \{\eta_i\}_{i \in I}$,
$U \subseteq spec(\cO_K)$ open,
$I \subseteq \{1, \ldots ,\gamma\}$, where
\be
\cO_{\cX_K}(U \coprod \{\eta_i\}_{i \in I}) =
\bigcap\limits_{p \in U}\cG(\cO_{K,p})\cap
\bigcap\limits_{i \in I} \cO_{K,\eta_i}
\ee
In particular, the global sections are
\be
\cO_{\cX_K}(\cX_K) = \bigcap\limits_{p \in spec\cO_K}
\cG(\cO_{K,p}) \cap \bigcap\limits_{i \leq \gamma} \cO_{K,\eta_i}=
\mathbb{F}[\mu_K]
\ee
with $\mu_K \subseteq \cO^*_K$ 
the group of roots of unity in $\cO^*_K$.

Every rational number $f \in \mathbb Q^*$,
defines a geometric map 
$\underline{f} \in \cG\cS (\cX , \mathbb P')$,
i.e. a collection of maps
$\underline{f}_N \in \cG\cG\cS (\cX_N , \mathbb P')$,
for $N$ divisible by $N_0 \cdot N_{\infty}$,
$N_0 = \prod\limits_{\nu_p(f) >0} p$,
$N_{\infty} = \prod\limits_{\nu_p(f) <0} p$,
with $\underline{f} \circ \pi_N^M = \underline{f}_M$.
For $f = \pm 1$ it is the constant map given by
\be
\mathbb F[z^{\mathbb Z}] \twoheadrightarrow \mathbb F [\pm 1] =
\cG(\mathbb Z) \cap A_N , \ \mbox{for any} \ N
\ee
\[z \mapsto f = \pm 1
\]
For $f \neq \pm 1$,
 we may assume $|f|_{\eta} < 1$, by the commutativity of
\be
 I \ \mbox{the inversion (6.2.8)}  
\ee
\begin{diagram}
 & & \mathbb P'\\
& \ruTo^{\underline{f}}& \\
\cX & &\uTo^{\wr}_{I}\\
  &\rdTo_{\underline{f^{-1}}} &\dTo\\
& & \mathbb P'
\end{diagram}

Thus for $N$ divisble by $N_0 \cdot N_{\infty}$
we have $f \in A_N$,
and the map $\underline{f}_N$ is given by
\be
\ee
\begin{diagram}
\cX_N = &spec \, \cG (\mathbb Z ) 
\coprod\limits_{spec \, \cG (\mathbb Z [\frac{1}{N}])} A_N& & \cG (\mathbb Z [\frac{1}{N}])& \supseteq & A_N &f \\
\dTo_{\underline{f}_N}= &f_{\mathbb Z} \coprod f_{\eta , N}& \mbox{, with} & \uTo & & \uTo^{f^{\sharp}_{\eta, N}}& \uMapsto\\
\mathbb P' = & spec \, \mathbb F [(z^{-1})^{\mathbb N}]
\coprod\limits_{spec \, \mathbb F [z^{\mathbb Z}]}
spec \, \mathbb F [z^{\mathbb N}]& &\mathbb F [z^{\mathbb Z}]& \supseteq& \mathbb F [z^{\mathbb N}] & z
\end{diagram}

and $f_{\mathbb Z}$ is as in (6.2.10).

Similarly for a number field $K$,
every element $f \in K^*$ defines a geometric map
\[
\underline{f} \in \cG\cS (\overline{spec \, \cO_K}, \mathbb P^1)
\]

\chapter{Products}
\section{Tensor product}
The category $\cG\cR_C$ of commutative (or the self-adjoint part 
$\cG\cR_C^+$)
generalized rings has tensor-products, i.e. fibred sums:
Given homomorphisms $\varphi^j \in \cG\cR_C (A,B^j)$
$j=0,1$, there exists $B^0 \bigotimes\limits_{A}B^1 \in \cG\cR_C$,
and homomorphisms 
$\psi^j \in \cG\cR_C (B^j,B^0 \bigotimes\limits_{A}B^1)$,
such that $\psi^0 \circ \varphi^0 = \psi^1 \circ \varphi^1$,
and for any $C \in \cG\cR_C$,
\be
\cG\cR (B^0 \bigotimes\limits_{A}B^1, C) = 
\cG\cR(B^0,C) \prod\limits_{\cG\cR(A,C)} \cG\cR (B^1,C)
\ee
So given homomorphisms $f^j \in \cG\cR (B^j,C)$ 
with $f^0 \circ \varphi^0 = f^1 \circ \varphi^1$,
there exists a unique homomorphism
$f^0 \otimes f^1 \in \cG\cR (B^0 \bigotimes\limits_{A}B^1 ,C)$,
such that $(f^0 \otimes f^1) \circ \psi^j = f^j$.
The construction of $B^0 \bigotimes\limits_{A}B^1$ goes as follows.
First for a finite set
$\{b_1^0, \ldots , b_n^0,b_1^1, \ldots , b_m^1 \}$,
where $b_i^j \in B^j_{X_i^j}$, 
we have the free commutative generalized ring on the sets
$\{X_1^0, \ldots , X_n^0,X_1^1, \ldots , X_m^1 \}$
(\ref{eq2.6.7}),
and we write $\underline{b}_1^0, \ldots , \underline{b}_n^0, 
\underline{b}_1^1, \ldots ,\underline{b}_m^1 $
for its canonical generators.
Taking the direct limit over such finite subsets 
(\ref{eq2.7.3}),
we have the free commutative generalized ring $\Delta$
with generators $\underline{b}$,
with $b \in B^0$ or $b \in B^1$.
We divide $\Delta$ by the eqivalence ideal 
$\varepsilon_A$ generated by
\be
\begin{array}{ccl}
\underline{b} \circ \underline{b'} \sim \underline{b \circ b'}&,&
b,b' \in B^j \ , \ j=0,1;\\
(\underline{b} , \underline{b'}) \sim \underline{(b ,b')}&,&
b,b' \in B^j \ , \ j=0,1;\\
\underline{1^j} \sim 1 &,&
\mbox{where} \ 1^j \in B^j_{[1]} \ \mbox{is the unit};\\
\underline{\varphi^0(a) } \sim \underline{\varphi^1(a)}&,&
\mbox{for} \ a \in A
\end{array}
\ee
The quotient generalized ring 
$\Delta/\varepsilon_A$ is the tensor product 
$B^0 \bigotimes\limits_{A}B^1$,
the homomorphism $\psi^j$ is given by 
$\psi^j(b) = \underline{b} \  mod \ \varepsilon_A$,
$b \in B^j$.

Note that every element of $(B^0 \bigotimes\limits_{A}B^1)_X$ 
can be expressed (non-uniquely) as
\be
(\underline{a},\underline{b}) = (\underline{a}_1 \circ 
\underline{a}_2 \circ \cdots \circ \underline{a}_n ,
\underline{b}_1 \circ \cdots \circ \underline{b}_m)\ mod \ \varepsilon_A
\ee
with $a_i \in B_{f_i}^{i(mod 2)}$,
$b_j \in B_{g_j}^{j(mod 2)}$,
and $f_1 \circ \cdots \circ f_n = c_X \circ g_1 \circ \cdots \circ g_m$
(where $c_X \in Set_{\bullet}(X,[1])$ is the canonical map,
$c_X(x) =1$ for all $x \in X$).
These elements  are multiplied and contracted by the formulas
(\ref{eq1.12.11}) and (\ref{eq1.12.12}).

\vspace{10pt}

\noindent (7.1.4) {\bf Example:}
For self-adjoint monoids $M_0$, $M_1$, $N$,
and homomorphisms $\psi^i \in Mon^+ (N, M_i)$,
$i=0,1$, we have (by adjunction (\ref{eq2.4.8})),
\[
\mathbb{F}[M_0] \bigotimes\limits_{\mathbb{F}[N]} 
\mathbb{F}[M_1] = \mathbb{F}[M_0 \bigotimes\limits_{N} M_1]
\]
where $M_0 \bigotimes\limits_{N} M_1$
is the fibered sum in the category $Mon^+$.
The monoid $M_0 \bigotimes\limits_{N} M_1$ is given by elements 
$m_0 \otimes m_1$, $m_i \in M_i$,
with relations
\[
m_0 \otimes 0 = 0 \otimes 0 = 0 \otimes m_1 \ , \ m_i \in M_i
\]
and
\[
m_0\cdot \psi^0(n) \otimes m_1 = 
m_0  \otimes \psi^1(n) \cdot m_1 \ , \ n \in N
\]

\vspace{10pt}

\noindent (7.1.5) {\bf Example:}
For a commutative (semi-) ring $B$, 
let $B^!$ denote the underlying multiplicative monoid of $B$ 
(i.e. forget addition),
and let $\mathbb{F}[B^!]$ denote the associated 
(commutative and self-adjoint) generalized ring, cf. (\ref{sec2.4}).
From the identity map $B^! 
\xrightarrow{=}\cG(B)_{[1]}$,
we obtain by adjunction (\ref{eq2.4.8}) 
the canonical injective homomorphism
$J_B \in \cG\cR(\mathbb{F}[B^!] , \cG(B))$.
The unique homomorphism of (semi) rings
$\mathbb{N} \rightarrow B$,
gives the unique homomorphism of generalized rings 
$I_B \in \cG\cR (\cG(\mathbb{N}), \cG(B) )$.
We get a canonical homomorphism of generalized rings,

\noindent (7.1.5.1) \quad
$\displaystyle{\Psi_B = I_B \otimes J_B \in \cG\cR\left(\cG(\mathbb{N}) 
\bigotimes\limits_{\mathbb{F}} \mathbb{F} [B^!] , \cG(B)
\right)}$

The homomorphism $\Psi_B$ is always surjective 
(as follows from (\ref{eq2.2.8})).
For any (self-adjoint) monoid $B$, the generalized ring
$n^B = \cG(\mathbb{N}) \bigotimes\limits_{\mathbb{F}} \mathbb{F}[B]$,
can be described as

\noindent (7.1.5.2) \quad $\displaystyle{n_X^B = \left\{
(\pi: \widetilde{X} \rightarrow X, \mu: \widetilde{X} \rightarrow B) \right\}/\approx
}$

The elements of $n_X^B$ are (isomorphism classes of) sets over 
$X \prod B$,
where the equivalence relation $\approx$
is generated by isomorphisms, i.e.

$(\pi:\widetilde{X} \rightarrow X, \mu: \widetilde{X} \rightarrow B)
\approx 
(\pi':\widetilde{X'} \rightarrow X, 
\mu': \widetilde{X'} \rightarrow B)$ if there is a bijection 
$\sigma:\widetilde{X} \rightarrow \widetilde{X}'$, 
$\pi=\pi' \circ \sigma$, 
$\mu= \mu' \circ \sigma$,

and by zero, i.e.

\noindent (7.1.5.3) \quad $\displaystyle{(\widetilde{X} , \pi , \mu)
\approx (\widetilde{X} \setminus \{x\} , 
\pi|_{\widetilde{X} \setminus \{x\}},
\mu|_{\widetilde{X} \setminus \{x\}})}$
if $\mu(x)=0$.

For $f \in Set_{\bullet}(X,Y)$, and for 
$(\widetilde{X} , \mu ) \in n_X^B$,
$(Z, \lambda) \in n_f^B$,
we have the contraction, cf. (\ref{eq2.2.12}),

\noindent (7.1.5.4) \quad 
$\displaystyle{\left( (\widetilde{X}, \mu),(Z, \lambda) \right) =
\left( \widetilde{X} \prod\limits_{X} Z,
(\mu, \lambda)\right)}$
\[
(\mu, \lambda)(x,z) = \mu(x)\cdot \lambda(z)
\]

For $(\widetilde{Y}, \mu ) \in n_Y^B$
we have the multiplication, cf. (\ref{eq2.2.14}),

\noindent (7.1.5.5) \quad $\displaystyle{(\widetilde{Y}, \mu) \circ
(Z, \lambda) =
(\widetilde{Y} \prod\limits_Y Z, \mu \circ \lambda)}$
\[\mu \circ \lambda (y,z) = \mu(y) \cdot \lambda(z)
\]
For a commutative (semi)ring $B$, the canonical homomorphism
(7.1.5.1)
$\Psi_B \in \cG\cR (n^{B!}, \cG(B))$ is given in this description as

\noindent (7.1.5.6) \quad $\displaystyle{\left(\Psi_B (\widetilde{X}, \mu )\right)_x = \sum\limits_{\stackrel{\textstyle{\tilde{x} \in \widetilde{X} }}{\pi(\tilde{x})=x}}\mu(\tilde{x})}$

To get such a surjective homomorphism we can use any multiplicative submonoid 
$B_0 \subseteq B^!$
such that $\mathbb{N}[B_0]=B$.
For example, for $B=\mathbb{Z}$ the integers,
we can take $B_0 = \{0, \pm 1 \}$,
and we get a surjective homomorphism

\noindent (7.1.5.7) \quad $\displaystyle{\Psi \in \cG\cR(\cG(\mathbb{N}) \bigotimes\limits_{\mathbb{F}} \mathbb{F} [\pm 1] , \cG(\mathbb{Z}))}$

\section{The arithmetical plane
$T^+ = \cG(\mathbb{N}) \bigotimes\limits_{\mathbb{F}}\cG(\mathbb{N})
$}

An \emph{oriented-tree} is a (rooted) tree $F$ together with a map
\[\varepsilon_F: F \setminus \partial F \rightarrow \{0,1 \}
\]
It is \emph{1-reduced} if $\nu(a) \neq 1$
for all $a \in F$.
If for some $a \in F$,
$S_F^{-1}(a) = \{a'\}$, we obtain by 1-reduction the tree
\be
\label{eq7.2.1}
 1_a(F) = F \setminus \{ a \}
\ee
with $S_{F'}(a') = S_F(a)$.

For every oriented tree $F$ there is a unique 1-reduced tree
$F_{1-\mbox{red}}$;
it is obtained from $F$ by a finite sequence of 1-reductions. 

The oriented tree $F$ is \emph{O-reduced} if for all 

$a \in F \setminus (\partial F \coprod \{0_F\})$,
$\varepsilon(a) \neq \varepsilon(S(a))$.
If for some $a  \in F \setminus (\partial F \coprod \{0_F\})$, 
$\varepsilon(a) = \varepsilon(S(a))$,
we obtain by \emph{O-reduction} the tree
\be
\label{eq7.2.2}
O_a(F) = F \setminus \{a\}
\ee
with $S_{O_a(F)}(a') = S_F(a)$
if $S_F(a')=a$.

For every oriented tree $F$ there is a unique O-reduced tree 
$F_{O-\mbox{red}}$;
it is obtained from $F$ by a finite sequence of O-reductions.
For a O-reduced oriented tree $F$, the orientation 
$\varepsilon_F$ is
completely determined by its value at the root 
$\varepsilon_F(0_F)$,
since $\varepsilon_F(x) \equiv \varepsilon_F(0_F)+ ht(x) (mod 2 )$.
Thus we view O-reduced oriented trees $F$ as ordinary trees together with an orientation of the root 
$\varepsilon_F = \varepsilon_F(0_F) \in \{0, 1 \}$.

The commutativity relation on oriented trees can 
be described as for ordinary trees, 
cf. (\ref{eq2.5.11}) - (\ref{eq2.5.15});
it is generated by $F \approx C^{\sigma}_{G,H}F$,
but now the isomorphisms 
$\sigma_a: H_a \xrightarrow{\sim} H$ has to preserve orientation.

It is also generated by the transpositions 
$F \approx C_{b}^{\sigma}F$, cf. (\ref{eq2.5.25}),
which are a special case of commutativity.
But note that even if $F$ is O -reduced, the tree 
$C_{G,H}^{\sigma}F$
need not be O-reduced;
in fact, already for transposition $C_b^{\sigma}F$
is (almost) always not O-reduced, and the associated O-reduced tree
$\bar{C}_bF = (C_b^{\sigma})_{0-\mbox{red}}$ 
can be described as follows:
for $b \in F \setminus \{ 0_F \}$, such that for all 
$a \in S_{F}^{-1}(b)$,
$\nu_F(a) =n$,
\be
\bar{C}_bF = F \setminus (\{ b \} \coprod S_F^{-1}(b) )
\ \mbox{with} \ 
S_{\bar{C}_bF}(x) = \left\{\begin{array}{ll} 
S_F(x) & \mbox{if} \ S_F^2(x) \neq b\\
S_F(b) & \mbox{if} \ S_F^2 (x) = b
\end{array}\right.
\ee

For O-reduced $F$ such that for all 
$a \in S_F^{-1}(0_F)$,
$\nu_F(a) =n$,
(the case of $b=0_F$ above), we have

\noindent $(7.2.3)_0$ \quad $\displaystyle{O_{0_F}F =
(F \setminus S_F^{-1}(0) ) \coprod [n]}$, 
as in (2.5.25)
but with the new orientations.

Note that the operations of 
$\varepsilon_{O_{0_F}}=1-\varepsilon_F$
 of 1-reduction, 
O-reduction, and transportation, do not alter the boundary of a tree.

We let $\approx$ denote the equivalence relation on oriented trees generated by
$1$-reductions, O-reductions, and commutativity relations.
We let $[F]$ denote the equivalence class
of the oriented tree $F$.
Thus $[F]=[F']$ if and only if there exist
$F=F_0, F_1, \ldots , F_l = F'$,
such that for $j=1, \ldots , l$,
the pair $\{F_j , F_{j-1} \}$
is related by $1$-reduction,
or O-reduction, or transposition;
it follows that there is a canonical identification of the 
boundaries:
$\partial F = \partial F'$.

\vspace{10pt}

For a finite set $X \in \Fbb$,
let $\Upsilon_X$ denote the collection of isomorphism classes of data
\be
\Upsilon_X = \{F=([F_1]; [\bar{F}_x], x \in X ; \sigma_F) \}/\cong
\ee
where $F_1$, $\bar{F}_x$
are oriented trees taken modulo $\approx$-equivalence,
and $\sigma_F$ is a bijection 
$\sigma_F: \partial F_1 \xrightarrow{\sim} \coprod\limits_{x \in X}
\partial \bar{F}_x$.
Thus explicitly, the data $F$ is equivalent to the data $F'$,
if and only if there exists
$F=F^0, F^1, \ldots , F^l=F'$
such that for $j=1, \ldots ,l$
the pair $\{F^j, F^{j-1} \} = \{ G,G' \}$
is related by either:

\vspace{10pt}

\noindent (7.2.5) {\bf Isomorphism:}
have isomorphism
$\tau_1: G_1 \xrightarrow{\sim} G'_1, 
\tau_x: \bar{G}_x \xrightarrow{\sim} \bar{G}'_x, x \in X$
such that 
$\sigma_{G'} \circ \tau_1(b) = \tau_x \circ \sigma_G(b)$
for
$b \in \partial G_1$,
$\sigma_G(b) \in \partial \bar{G}_x$.

\vspace{10pt}

\noindent (7.2.6) {\bf $1$-reduction:}
have $G' = 1_aG$,
for some
$a \in G_1 \coprod {\displaystyle \coprod\limits_{x \in X} } \bar{G}_x$
with $\nu(a) =1$,
cf. (\ref{eq7.2.1}).

\vspace{10pt}

\noindent (7.2.7) {\bf O-reduction:}
have $G' =O_a G$,
for some \\$a \in \left(G_1 \setminus (\partial G_1 \coprod \{ 0\} )
\right)\coprod \displaystyle{\coprod\limits_{x \in X}}
\bar{G}_x \setminus (\partial \bar{G}_x \coprod \{0 \})$
with $\varepsilon (a) = \varepsilon (S(a))$, cf. (\ref{eq7.2.2}).

\vspace{10pt}

\noindent (7.2.8) {\bf Transposition:}
\setcounter{equation}{8}
have $G' = C_b^{\tau}G$,
for some $b \in G_1 \coprod \displaystyle{\coprod\limits_{x \in X}} 
\bar{G}_x$,
such that for all $a \in S^{-1}(b)$,
$\varepsilon(a) \equiv \varepsilon$, and
$\tau_a : S^{-1}(a) \xrightarrow{\sim} [n]$
bijection, cf. (\ref{eq2.5.25}).

\vspace{10pt}

The operations of multiplication (\ref{eq2.5.7}),
and of contraction (\ref{eq2.5.9}),
induce well defined operations on equivalent classes of data,
and make $\Upsilon$ into a commutative (but non-self-adjoint) 
generalized ring.
It is straightforward to check that
\be
\begin{array}{cl}
F \circ G & \cong F \circ (1_a G) \cong F \circ (O_a G) 
            \cong F \circ (C_b^{\tau} G)\\
          &\cong (1_a F) \circ G \cong (O_a F) \circ G 
            \cong (C_b^{\tau} F) \circ G\\
(F,G) & \cong (F , 1_a G) \cong (F , O_a G) 
            \cong (F,  (C_b^{\tau} G)\\
      &\cong (1_a F,  G) \cong (O_a F, G) 
            \cong (C_b^{\tau} F,  G)
\end{array}
\ee

whenever the operations
$1_a$, $O_a$, $C_b^{\tau}$ are relevant,
and that $\Upsilon$ satifies the axioms of a commutative generalized ring.
Note that for $\varepsilon =0, 1$,
we have the elements
\be
\delta_X^{\varepsilon} = \left( [X \coprod \{ 0 \}];
      [0_x] , x \in X ; \sigma \right) \in \Upsilon_X
\ee
where
$X \coprod \{ 0 \}$ is the oriented tree with
$\varepsilon (0) = \varepsilon$,
$S(x) =0$ for $x \in X$,
and $\sigma :X \xrightarrow{\sim} 
    \coprod\limits_{x \in X} \{ 0_x \}$ is the natural bijection 
$\sigma (x) = 0_x$.

For $f \in Set_{\bullet}(X,Y)$,
and $(\delta_{f}^{\varepsilon})^{(y)} =
      \delta_{f^{-1}(y)}^{\varepsilon} $,
$y \in Y$,
we have via O-reduction
\be
\label{eq7.2.11}
\delta_Y^{\varepsilon} \circ \delta_f^{\varepsilon}
\cong
\delta_{D(f)}^{\varepsilon}
\ee
we also have by $1$-reduction
\be
\label{eq7.2.12}
\delta_{[1]}^{\varepsilon} = 
([\{0\} \coprod \{1\}]; [0_1];\sigma)\cong
([0_1];[0_1];id)=
1 \in \Upsilon_{[1]}
\ee
Thus we get homomorphisms, cf. (2.2.16),
$\Psi^{\varepsilon} \in \cG\cR(\cG(\mathbb{N}), \Upsilon)$
with 
$\Psi^{\varepsilon}(\mathbbm{1}_X) = \delta_X^{\varepsilon}$.
It is clear that $\Upsilon$ is generated by the 
$\delta_X^{\varepsilon}$, and the only relations they satisfy
are (\ref{eq7.2.11}), (\ref{eq7.2.12}), and commutativity.
It follows that $\Upsilon$ is the sum of $\cG(\mathbb{N})$ 
with itself in the category of commutative generalized rings:
for any commutative $A \in \cG\cR_C$,
\be
\begin{array}{ccl}
\cG\cR(\cG(\mathbbm{N}), A) \times \cG\cR(\cG(\mathbb{N}), A)&
\stackrel{\sim}{\leftrightarrow}&
\cG\cR(\Upsilon, A)\\
(\varphi \circ \psi^0 , \varphi \circ \psi^1)&
\leftmapsto&
\varphi\\
\varphi^0, \varphi^1 & \mapsto &
\varphi^0 \otimes \varphi^1 (\delta_X^{\varepsilon}) :=
\varphi^{\epsilon}(\mathbbm{1}_X)
\end{array}
\ee

The tensor product 
$\cG(\mathbb{N}) \bigotimes\limits_{\mathbb{F}} \cG ( \mathbb{N})$
in the category of commutative 
\underline{and self-adjoint}
generalized rings
$\cG\cR^+_C$ is therefore the self-adjoint quotient of $\Upsilon$,
$\Upsilon^+ = \Upsilon / \varepsilon^+$.
Here $\varepsilon^+$ is the equivalence ideal of $\Upsilon$ 
generated by
\[
F = ([F_1] ; [\bar{F}_1]; \sigma ) \sim
F^t = ([\bar{F}_1]; [F_1] ; \sigma^{-1})
\]
Passing to this quotient is equivalent to adding the 
following identifications on elements
$F = ([F_1]; \{[\bar{F}_x]\}_{x \in X}; \sigma) \in \Upsilon_X$,

(for a tree $F$, and $b \in F$, we write 
$F/b= \coprod\limits_{n \geq 0} S_{F}^{-n}(b)$
for the full subtree of $F$ with root $b$):

if for some $b_1 \in F_1$,
and some $b_{x_0} \in \bar{F}_{x_0}$,
$x_0 \in X$, 
the bijection $\sigma$ satisfy 
$\sigma(\partial \left(F_1/b_1\right)) 
=
\partial \left(\bar{F}_{x_0}/ b_{x_0}\right)$,
then $F$ is identified with 
\be
\label{eq7.2.14}
(b_1,b_{x_0})F :=
\ee
\[
([(F_1\setminus \left(F_1/b_1 \right) \coprod
\left(\bar{F}_{x_0}/ b_{x_0}\right)];
\{\bar{F}_{x}\}_{x \neq x_0},
\{(\bar{F}_{x_0}\setminus \left(\bar{F}_{x_0}/ b_{x_0}\right))
\coprod
\left(F_1/b_1 \right) \}; 
\widetilde{\sigma})
\]
with 
$\widetilde{\sigma} = \sigma|_{\partial F_1 \setminus \partial 
\left(F_1/b_1 \right)}
\coprod
\sigma^{-1}|_{\partial \left(\bar{F}_{x_0}/ b_{x_0}\right)}$

Unfortunately, the self-adjunction relation
(\ref{eq7.2.14})
is not compatible with the transposition relation (7.2.8).
We do not have canonical representatives for the elements of 
$\Upsilon_X^+$.

The diagaonal homomorphism
\be
\nabla \in \cG\cR(\Upsilon^+, \cG(\mathbb{N}))
\ee
is determined by 
$\nabla_X (\delta_X^{\varepsilon}) = \mathbbm{1}_X$,
and is given explicitly by
\be
\label{eq7.2.16}
\nabla_X ([F_1]; [\bar{F}_x]_{x \in X}; \sigma )/\hspace{-4pt}\cong \  =
(\underline{\sigma}: \partial F_1 \rightarrow X)/ \cong
\ee
(realizing $\cG(\mathbb{N})_X$
as isomorphism classes of sets over $X$).
The homomorphism $\nabla$ is surjective,
but it is not injective.
For calculations of $\nabla^{-1}_{[1]}(N)$ for $N=1,2,3$
see the appendix.

For a (self-adjoint) monoid $B$,
the tensor product 
$\Upsilon^+ \bigotimes\limits_{\mathbb{F}} \mathbb{F}[B]$
can be described as isomorphism classes of data
\be
(\Upsilon^+ \bigotimes\limits_{\mathbb{F}}\mathbb{F}[B])_X
:=
\{F=([F_1]; \{[\bar{F}_x]\}_{x \in X}; \sigma_F; \mu_F)\}/\cong
\ee
Here the data $([F_1]; \{[\bar{F}_x]\}_{x \in X}; \sigma_F)$
is the data for $\Upsilon^+_X$, 
and $\mu_F$
is a map $\mu_F: \partial F_1 \rightarrow B$,
and isomorphisms are required to preserve the $B$-valued maps,
and the zero law
cf. (7.1.5.3) holds in the form:
\[
\mu_F(b) = 0 , \ \underline{\sigma_F}(b) = x_0 \Rightarrow
F \approx ([F_1 \setminus \{b\}]; \{[\bar{F}_x ]\}_{x \neq x_0}
\cup \{[\bar{F}_{x_0}] \setminus \{\sigma _F(b)\})
\]
The operations of multiplication and contraction are the 
given ones on the $\Upsilon^+$-part of the data 
(i.e. given by (\ref{eq2.5.7}) and (\ref{eq2.5.9}),
and are given on the $B$-valued maps by 
(using the notations of (2.5.6) and (2.5.8)):
\be
\begin{array}{ll}
\mu_{G \circ F}(b,a) = \mu_G(b) \cdot 
\mu_{F_{\underline{\tau}(b)}}(a)&
, b \in \partial G_1 , a \in \partial F_{\underline{\tau}(b)}\\
\mu_{(G, F)}(b,a) = \mu_G(b) \cdot 
\mu_{F_{f \circ \underline{\tau}(b)}}(\sigma^{-1}a)&
, b \in \partial G_1 , a \in \partial \widebar{F}_{\underline{\tau}(b)}
\end{array}
\ee

For commutative rings $B_0$, $B_1$, 
taking $B = B^!_0 \otimes B_1^!$
(the sum in $Mon^+$, cf. (7.1.4),
we get the generalized ring
\[
\Upsilon^+ \bigotimes\limits_{\mathbb{F}} \mathbb{F}[B_0^! \otimes
B_1^!]= \cG(\mathbb{N}) \bigotimes\limits_{\mathbb{F}} 
\mathbb{F}[B_0^!] \bigotimes\limits_{\mathbb{F}} 
\cG(\mathbb{N}) \bigotimes\limits_{\mathbb{F}} 
\mathbb{F}[B_1^!]
\]
which maps surjectively onto
$\cG(B_0)\bigotimes\limits_{\mathbb{F}} \cG(B_1)$.

For the integers $\mathbb{Z}$, taking $B=\{0, \pm 1\}$,
we get the generalized ring
\be
\Upsilon^+ \bigotimes\limits_{\mathbb{F}} \mathbb{F}[\pm 1] =
(\cG(\mathbb{N}) \bigotimes\limits_{\mathbb{F}} \mathbb{F}[\pm 1])
\bigotimes\limits_{\mathbb{F}[\pm 1]}
(\cG(\mathbb{N}) \otimes \mathbb{F}[\pm 1])
\ee
which maps surjectively onto 
$\cG(\mathbb{Z}) \bigotimes\limits_{\mathbb{F}[\pm 1]} \cG
(\mathbb{Z}) = \mathbb{T}$

\section{Products of Grothendieck-Generalized-schemes}
\label{sec7.3}

The categor $\cal{GGS}$
has fibred products:

Given maps $f^j \in \cG\cG\cS(X^j, Y)$,
there exists $X^0 \Pi_YX^1 \in \cG\cG\cS$,
and maps $\pi_j \in \cG\cG\cS(X^0 \Pi_YX^1, X^j)$,
with $f^0 \circ \pi_0 = f^1\circ\pi_1$,
and for any $g^j \in \cG\cG\cS(Z,X^j)$,
with $f^0 \circ g^0 = f^1\circ g^1$,
there exists a unique map
$g^0 \pi g^1 \in \cG\cG\cS(Z,X^0\Pi_YX^1)$,
such that $\pi_j \circ (g^0 \pi g^1) = g^j$, $j=0,1$.

Writing $Y=\bigcup\limits_i spec (A_i)$,
$(f^j)^{-1}(spec(A_i))= \bigcup\limits_k spec (B^j_{i,k})$,
the fibred product $X^0 \Pi_Y X^1$
is obtained by gluing 
$spec (B^0_{i,k_0} \bigotimes\limits_{A_i} B^1_{i,k_1})$.
See the construction of fibred product of ordinary schemes,
e.g. \cite[Theorem 3.3, p. 87]{Hart}.

\section{Products of Generalized-schemes}
\label{sec7.4}

The category $\cal{GS}$
has fibred products.
This is an immediate corollary of 
(\ref{sec7.3}). Given maps
\[
\varphi = \{\varphi_i^j\}_{j \geq \sigma(i)} \in
\cG\cS ( \{X_j\}_{j \in J}, \{Y_i\}_{i \in I}),
\]
and
\[
\varphi' = \{\varphi_i^{j'}\}_{j' \geq \sigma'(i)} \in
\cG\cS ( \{X'_{j'}\}_{j' \in J'}, \{Y_i\}_{i \in I})
\]
the fibred product of $\varphi$
and $\varphi'$ in $\cG\cS$ is given by the inverse system
$\{X_j \Pi_{Y_i} X'_{j'}\}$,
the indexing set is
\[
\{(j,j',i) \in J \times J' \times I\  | \  j \geq \sigma(i) , \ j' \geq \sigma'(i) \}
\]

\section{The Arithmetical plane \\
$ \mathbb{X} = \overline{spec \mathbb{Z}} 
\prod\limits_{\mathbb{F}[\pm 1]}
\overline{spec \mathbb{Z}} $}

This is a special case of (\ref{sec7.4}):
The (compactified) arithmetical plane $\mathbb{X}$
is given by the inverse system
$\{ \cX_N \prod\limits_{spec \mathbb{F}[\pm 1]} \cX_M \}$,
with indexing set
$\{ (N,M) \in \mathbb{N} \times \mathbb{M} |
N, M \ \mbox{square-free} \}$
and with
\[
\cX_N = spec \, \cG(\mathbb{Z}) 
\coprod\limits_{spec \cG( \mathbb{Z} [\frac{1}{N}])} 
spec (\cG( \mathbb{Z} [\frac{1}{N}]) \cap \cO_{\eta})
\]
as in (\ref{sec6.4}).
This generalized scheme $\mathbb{X}$ contains the open dense subset,
\[
spec \cG(\mathbb{Z}) \prod\limits_{\mathbb{F}[\pm 1]}
spec \cG(\mathbb{Z}) =
spec (\cG(\mathbb{Z}) \bigotimes\limits_{\mathbb{F} [\pm 1]} 
\cG(\mathbb{Z}))
\]
A basis for neighborhoods of $(p, \eta)$
is given by 
\[\cG(\mathbb{Z}[\frac{1}{N}]) 
\bigotimes\limits_{\mathbb{F} [\pm 1]} 
(\cG(\mathbb{Z}[\frac{1}{M}])\cap \cO_{\eta})\]
where $p$ does not divide $N$, and $M$ is arbitrary.

Similarly, for any number field $K$ we 
have the compactified surface
\[
\overline{spec \cO_K} \prod\limits_{spec \mathbb{F}[\mu_K]}
\overline{spec \cO_K}
\]
It contains the open dense subset
$spec (\cG(\cO_K) \bigotimes\limits_{\mathbb{F} [\mu_K]} 
\cG(\cO_K))
$.

\section*{Appendix}
We enumerate the isomorphism classes of oriented trees $F$, 
which are 1-reduced and o-reduced,
 for $N = \sharp \partial F = 0,1,2, \ldots, 6$,
and we show the commutativity classes.
Note that every such tree gives a partition:
$N -1 = \sum\limits_{a \in F \setminus \partial F} \nu(a) -1$.
We mark with a circle the points of 
$F \setminus (\partial F \coprod S(\partial F))$.
The orientation of the tree is denoted by
$\varepsilon = \varepsilon_F(0) \in \{ 0,1 \}$. 

\newpage

\noindent {\bf \underline{$N=1$}:}
$\bullet_{\varepsilon} \approx 
\bullet_{1-\varepsilon} \hspace{-.8cm} -$
\quad , \quad
$\varepsilon = 0,1$

\vspace{10pt}

\noindent {\bf \underline{$N=2$}:} $(1):$
\includegraphics[scale=0.5]{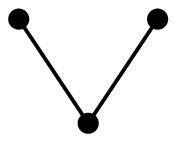}\hspace{-.5cm} $\varepsilon$

\vspace{10pt}

\noindent {\bf \underline{$N=3$}:} $(2):$
\includegraphics[scale=0.5]{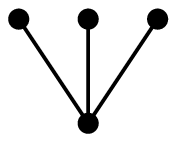}\hspace{-.5cm} 
   $\varepsilon$
\quad
$(1,1):$ 
\includegraphics[scale=0.5]{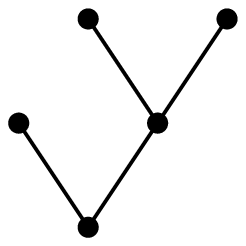}\hspace{-.9cm} 
   $\varepsilon$

\noindent {\bf \underline{$N=4$}:} $(3):$
\includegraphics[scale=0.5]{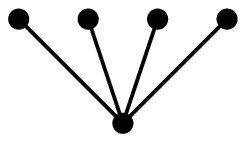}\hspace{-.7cm} $\varepsilon$

$(2,1):
$ \includegraphics[scale=0.5]{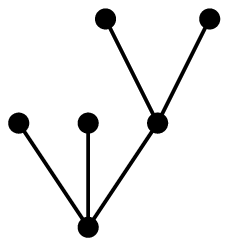}\hspace{-.7cm} $\varepsilon$
 \quad , \quad
\includegraphics[scale=0.5]{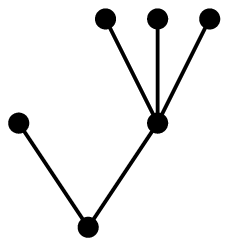}\hspace{-.8cm} $\varepsilon$

$(1,1,1):$ \framebox{\includegraphics[scale=0.5]{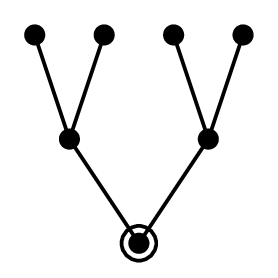}\hspace{-.7cm} $\varepsilon$ \hspace{.5 cm}
$\stackrel{\approx}{\leftrightarrow}$
\includegraphics[scale=0.5]{fig7_N4_4.eps}\hspace{-.8cm} 
$1-\varepsilon$}
\quad , \quad
\includegraphics[scale=0.5]{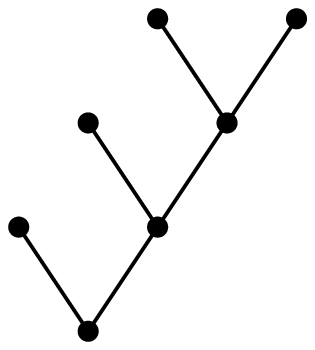}\hspace{-1.3cm} $\varepsilon$

\vspace{10pt}

\noindent {\bf \underline{$N=5$}:} $(4):$
\includegraphics[scale=0.5]{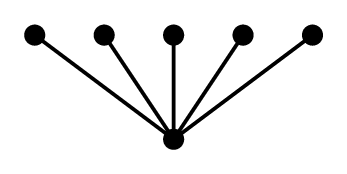}\hspace{-.9cm} $\varepsilon$
\quad
$(2,2):$
\includegraphics[scale=0.5]{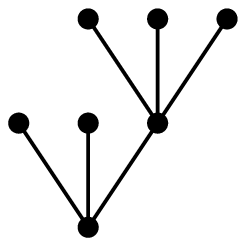}\hspace{-.9cm} $\varepsilon$

$(3,1):$
\includegraphics[scale=0.5]{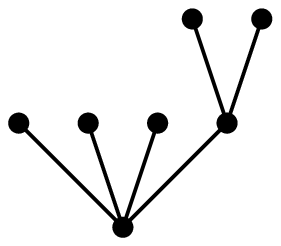}\hspace{-.9cm} $\varepsilon$
\quad , \quad
\includegraphics[scale=0.5]{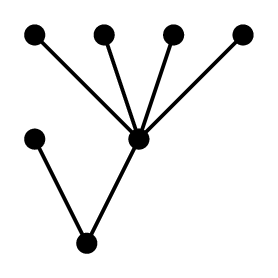}\hspace{-.9cm} $\varepsilon$

\begin{tabular}{c|c|ccc|c|}\cline{2-2} \cline{6-6}
$(2,1,1):$&\includegraphics[scale=0.5]{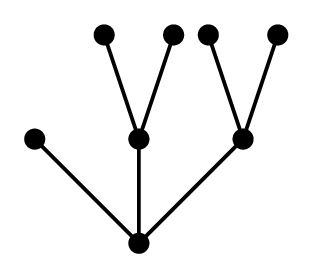}\hspace{-.9cm} $\varepsilon \quad $&,\includegraphics[scale=0.5]{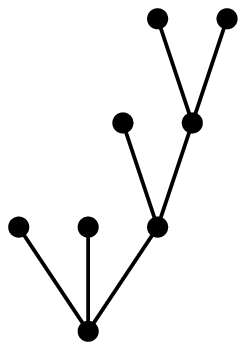}\hspace{-.9cm} $\varepsilon$,&\includegraphics[scale=0.5]{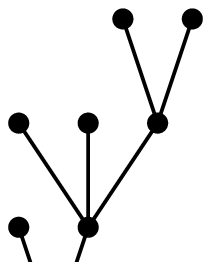}\hspace{-.9cm} $\varepsilon$, &\includegraphics[scale=0.5]{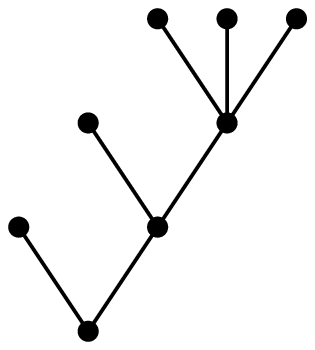}\hspace{-1.2cm} $\varepsilon \quad \quad$,&\includegraphics[scale=0.5]{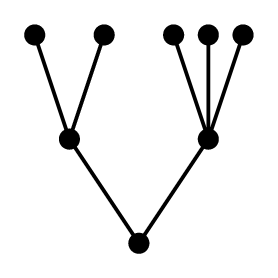}\hspace{-.8cm} $\varepsilon$\\
          &$\uparrow \wr \wr$&  & & & $\uparrow \wr \wr$\\
$(1,1,1,1):$&\includegraphics[scale=0.5]{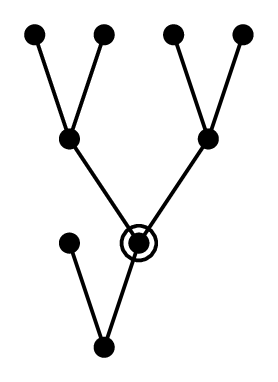}\hspace{-.9cm} $\varepsilon \quad $&  &\includegraphics[scale=0.5]{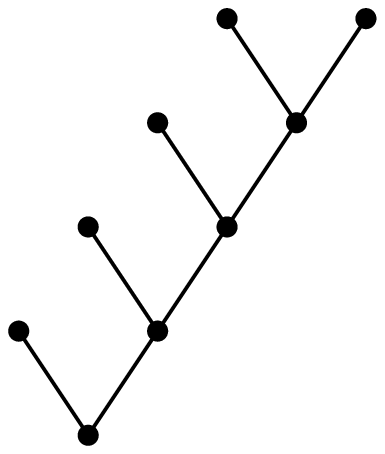}\hspace{-1.5cm} $\varepsilon$& &\includegraphics[scale=0.5]{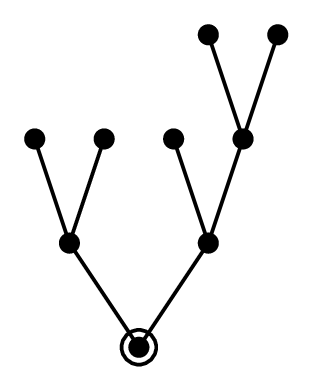}\hspace{-.9cm} $\varepsilon \quad$\\
\cline{2-2} \cline{6-6}
\end{tabular}

\newpage

\noindent {\bf \underline{$N=6$}:} $(5):$
\includegraphics[scale=0.5]{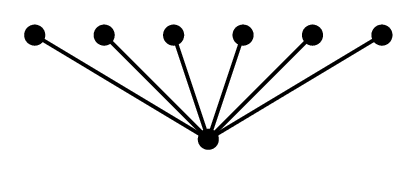}\hspace{-1cm} $\varepsilon$

$(4,1):$
\includegraphics[scale=0.5]{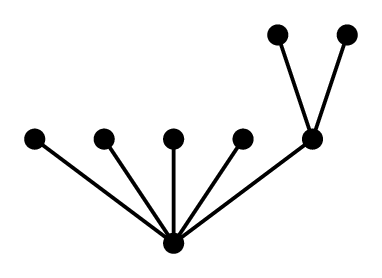}\hspace{-1cm} $\varepsilon$
\quad
\includegraphics[scale=0.5]{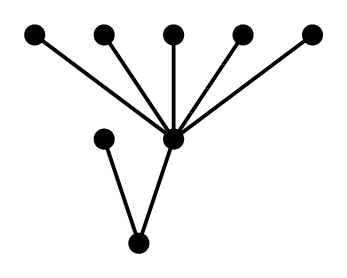}\hspace{-1cm} $\varepsilon$

$(3,2):$
\includegraphics[scale=0.5]{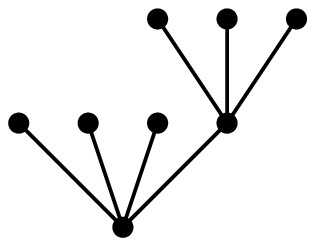}\hspace{-1cm} $\varepsilon$
\quad
\includegraphics[scale=0.5]{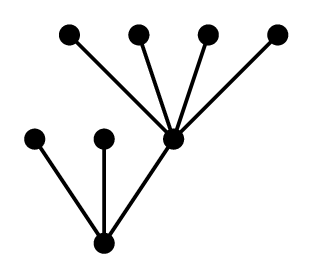}\hspace{-1cm} $\varepsilon$

$(3,1,1):$

\begin{tabular}{c|c|c|c|c|c|}\cline{2-2}\cline{4-4}\cline{6-6}
$(3,1,1):$&
\includegraphics[scale=0.5]{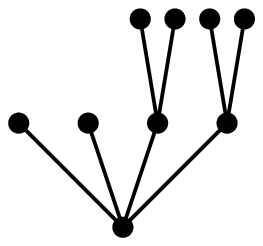}\hspace{-.9cm} $\varepsilon$\quad \quad&
\includegraphics[scale=0.5]{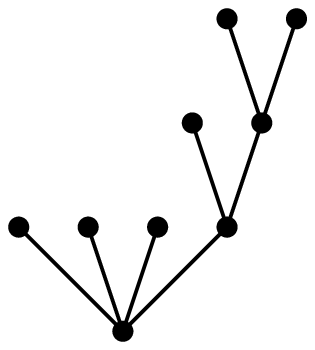}\hspace{-1.1cm} $\varepsilon$
\quad
\includegraphics[scale=0.5]{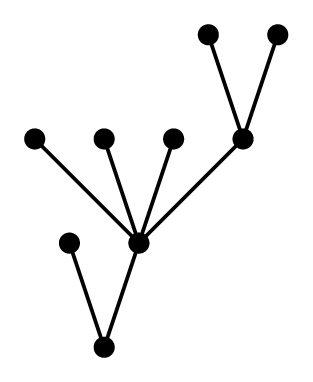}\hspace{-1.1cm} $\varepsilon$
\quad \quad
\includegraphics[scale=0.5]{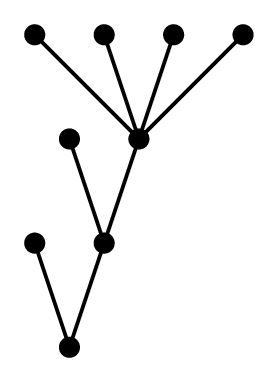}\hspace{-1.1cm} $\varepsilon$&
\includegraphics[scale=0.5]{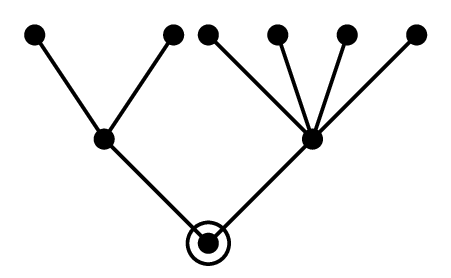}\hspace{-1.1cm} $\varepsilon$&\quad& \\ $(2,2,1):$&
$\uparrow$& 
\includegraphics[scale=0.5]{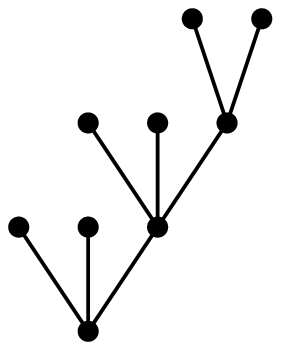}\hspace{-1.1cm} $\varepsilon$
\quad
\includegraphics[scale=0.5]{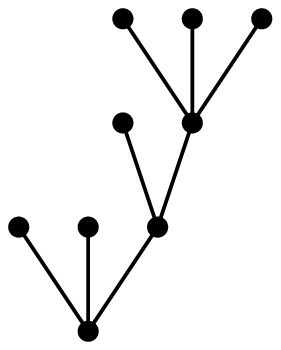}\hspace{-1.1cm} $\varepsilon$
\quad \quad
\includegraphics[scale=0.5]{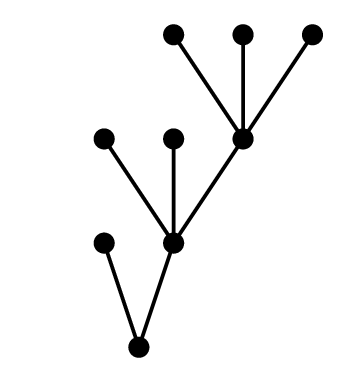}\hspace{-1.1cm} $\varepsilon$&
\includegraphics[scale=0.5]{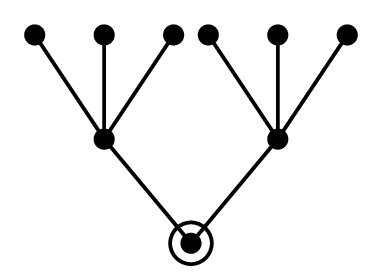}\hspace{-.9cm} $\varepsilon$&&
\includegraphics[scale=0.5]{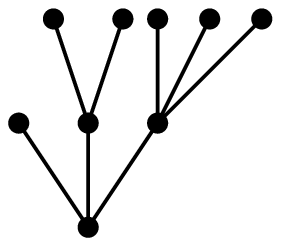}\hspace{-1.1cm} $\varepsilon$ \quad \quad \\
$(2,1,1,1):$&
\includegraphics[scale=0.5]{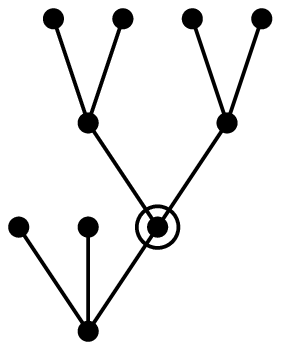}\hspace{-1.1cm} $\varepsilon$\quad \quad &
\includegraphics[scale=0.5]{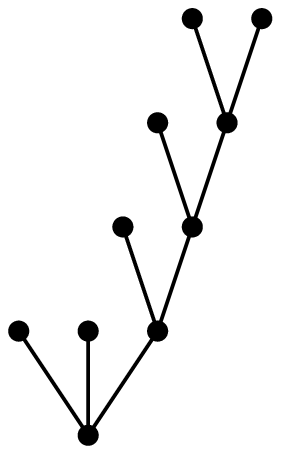}\hspace{-1.1cm} $\varepsilon$
\quad
\includegraphics[scale=0.5]{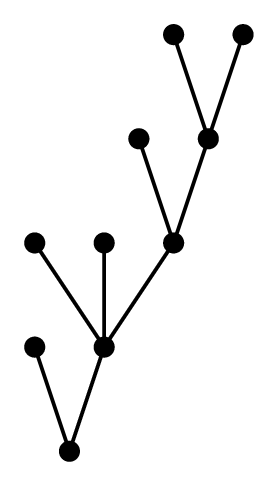}\hspace{-1.1cm} $\varepsilon$
\quad
\includegraphics[scale=0.5]{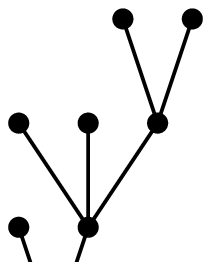}\hspace{-1.1cm} $\varepsilon$
\quad
\includegraphics[scale=0.5]{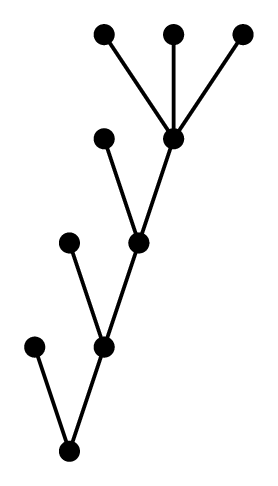}\hspace{-1.1cm} $\varepsilon$\quad \quad \quad&
\includegraphics[scale=0.5]{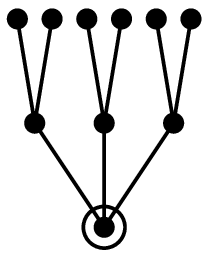}\hspace{-1.1cm} $\varepsilon$
\quad \quad
\includegraphics[scale=0.5]{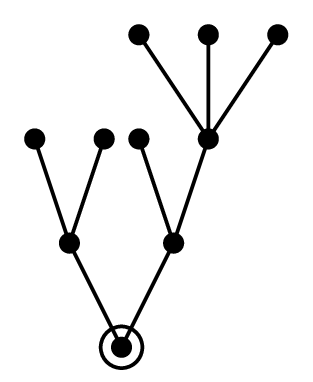}\hspace{-.9cm} $\varepsilon$\quad \quad \quad& & $\uparrow$\\ \cline{2-2}
\end{tabular}
\newpage

We enumerate the elements of
$\nabla_{[1]}^{-1}(N) \subseteq \left(\cG(\mathbb N) \bigotimes\limits_{\mathbb F} \cG(\mathbb N)\right)_{[1]}$ for $N=1,2,3$.
The bijection $\sigma$ is described by the numbering at the boundaries of the trees.
We take only the 1-reduced and o-reduced trees - the orientation is given by the number next to the root.

\noindent \underline{$\bf N=1:$} $1 = \left(\bullet , \bullet \right)$

\vspace{10pt}

\noindent \underline{$\bf N=2$:} \quad {\Huge (}\includegraphics[scale=0.5]{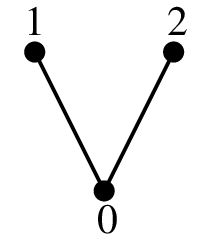},
\includegraphics[scale=0.5]{fig_appen120.eps}{\huge )},
{\Huge (}\includegraphics[scale=0.5]{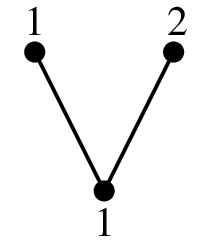},
\includegraphics[scale=0.5]{fig_appen121.eps} {\Huge )},

\vspace{10pt}

\hspace{1cm} \framebox{{\Huge (} \includegraphics[scale=0.5]{fig_appen120.eps},
\includegraphics[scale=0.5]{fig_appen121.eps} {\Huge )} $\approx$
{\Huge (} \includegraphics[scale=0.5]{fig_appen121.eps},
\includegraphics[scale=0.5]{fig_appen120.eps} {\Huge )}}

\vspace{10pt}

\noindent \underline{$\bf N=3$:} \quad 
{\Huge (}\includegraphics[scale=0.5]{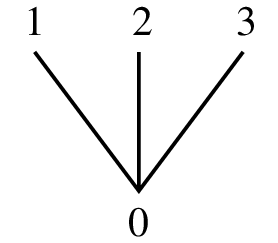},
\includegraphics[scale=0.5]{fig_appen1230.eps} {\Huge )}, \quad
{\Huge (}\includegraphics[scale=0.5]{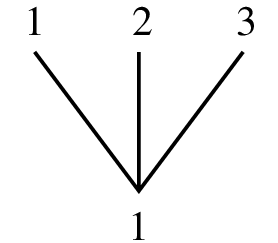},
\includegraphics[scale=0.5]{fig_appen1231.eps} {\Huge )},

\vspace{10pt}

\begin{tabular}{|ccc|} \hline
{\Huge (}\includegraphics[scale=0.5]{fig_appen1230.eps},
\includegraphics[scale=0.5]{fig_appen1231.eps} {\Huge )}&
$\approx$ &
{\Huge (}\includegraphics[scale=0.5]{fig_appen1231.eps},
\includegraphics[scale=0.5]{fig_appen1230.eps} {\Huge )}\\ 
$\wr \wr $ & & $\wr \wr $\\ 
{\Huge (}\includegraphics[scale=0.5]{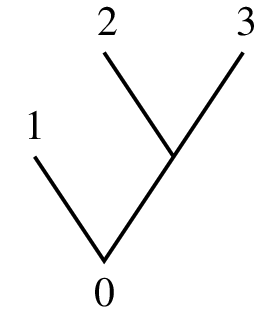},
\includegraphics[scale=0.5]{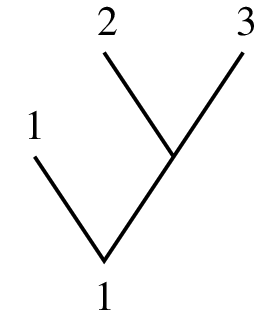} {\Huge )}&
$\approx$ &
{\Huge (}\includegraphics[scale=0.5]{fig_tree1231.eps},
\includegraphics[scale=0.5]{fig_tree1230.eps} {\Huge )}\\ \hline
\end{tabular}

\vspace{10pt}

{\Huge (}\includegraphics[scale=0.5]{fig_tree1230.eps},
\includegraphics[scale=0.5]{fig_tree1230.eps} {\Huge )},
{\Huge (}\includegraphics[scale=0.5]{fig_tree1231.eps},
\includegraphics[scale=0.5]{fig_tree1231.eps} {\Huge )}

\vspace{10pt}

{\Huge (}\includegraphics[scale=0.5]{fig_tree1230.eps},
\includegraphics[scale=0.5]{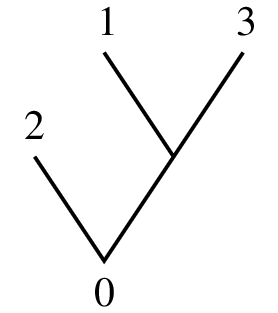} {\Huge )},
{\Huge (}\includegraphics[scale=0.5]{fig_tree1231.eps},
\includegraphics[scale=0.5]{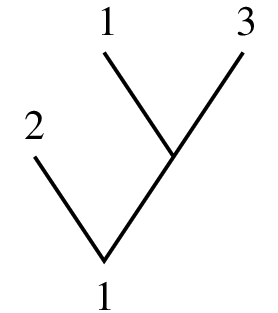} {\Huge )}

\vspace{10pt}

\framebox{{\Huge (}\includegraphics[scale=0.5]{fig_tree1230.eps},
\includegraphics[scale=0.5]{fig_tree2131.eps} {\Huge )} $\approx$
{\Huge (}\includegraphics[scale=0.5]{fig_tree2131.eps},
\includegraphics[scale=0.5]{fig_tree1230.eps} {\Huge )}
}

\vspace{10pt}

\begin{tabular}{|c|c|c|}\cline{1-1} \cline{3-3}
{\Huge (}\includegraphics[scale=0.5]{fig_appen1230.eps},
\includegraphics[scale=0.5]{fig_tree1230.eps} {\Huge )}& ,&
{\Huge (}\includegraphics[scale=0.5]{fig_appen1231.eps},
\includegraphics[scale=0.5]{fig_tree1231.eps} {\Huge )}
\\
$\wr \wr$ & & $\wr \wr$\\
{\Huge (}\includegraphics[scale=0.5]{fig_tree1230.eps},
\includegraphics[scale=0.5]{fig_appen1230.eps} {\Huge )} &, &
{\Huge (}\includegraphics[scale=0.5]{fig_tree1231.eps},
\includegraphics[scale=0.5]{fig_appen1231.eps} {\Huge )}
\\
\cline{1-1} \cline{3-3}
\end{tabular}

\vspace{10pt}

\framebox{{\Huge (}\includegraphics[scale=0.5]{fig_appen1230.eps},
\includegraphics[scale=0.5]{fig_tree2131.eps} {\Huge )} $\approx$
{\Huge (}\includegraphics[scale=0.5]{fig_tree1231.eps},
\includegraphics[scale=0.5]{fig_appen1230.eps} {\Huge )}
}

\vspace{10pt}

\framebox{{\Huge (}\includegraphics[scale=0.5]{fig_tree1230.eps},
\includegraphics[scale=0.5]{fig_appen1231.eps} {\Huge )} $\approx$
{\Huge (}\includegraphics[scale=0.5]{fig_appen1231.eps},
\includegraphics[scale=0.5]{fig_tree1230.eps} {\Huge )}
}

\chapter{Divisors}

Let 
$X = \{ X_N , N \in \cN ; 
\pi^M_N : X_M \rightarrow X_N , M \geq N \}$
be a generalized scheme.

\section{Meromorphic functions}
\label{sec8.1}
The sheaf of
\emph{meromorphic functions} of level
$N \in \cN$,
$\cK_N \in \cG\cR / X_N$,
is the sheaf of generalized rings over 
$X_N$ associated to the presheaf whose sections over an open set 
$U \subseteq X_N$
are the elements of the quotient ring
$S^{-1}_U \cO_{X_N}(U)$,
where $S_U$ is the multiplicative set in $\cO_{X_N}(U)_{[1]}$ 
given by
\be
S_U = \left\{ a \in \cO_{X_N}(U)_{[1]} , \ \mbox{for all} \
     M \geq N , \ \mbox{all open sets} \
     V \subseteq (\pi_N^M)^{-1}(U),
\right.
\ee
\[
\left.
 \ \mbox{all} \
      b_1, b_2 \in \cO_{X_M}(V), \quad
 (\pi_N^M)^{\sharp} a \circ b_1 =
      (\pi_N^M)^{\sharp} a \circ b_2 \Rightarrow b_1 = b_2\right\}
\]
\vspace{10pt}

We have the canonical injections 
$\cO_{X_N} \hookrightarrow \cK_N$,
and these are compatible:
\be
\mbox{for} \ M \geq N , \mbox{we have a commutative diagram in} \
 \cG\cR/ X_N
\ee
\[\begin{array}{rcl}
\cO_{X_N} & \hookrightarrow & \cK_N \\
(\pi_N^M)^{\sharp} \downarrow & & \downarrow\\
(\pi_N^M)_*\cO_{X_M} & \hookrightarrow & (\pi_N^M)_* \cK_M
\end{array}
\]
We have an exact sequence of sheaves of abelian groups on 
$X_N$
given by the units,
\be
* \rightarrow \cO^*_{X_N} \rightarrow
\cK^*_N \rightarrow \cK_N^*/ \cO^*_{X_N} \rightarrow *
\ee

\section{Cartier divisors}
\label{sec8.2}
A \emph{Cartier divisor} on $X$,
of level $N$, 
is a global section of the sheaf $\cK_N^*/\cO_{X_N}^*$.
The abelian group of Cartier divisors is denoted
\[
Div (X)_N = \cK^*_N/\cO_{X_N}^* (X_N)
\]
For an invertible meromorphic function $f \in \cK^*_N(X_N)$,
let $\underline{div(f)}$
denote the image of $f$ in $Div(X)_N$.

Thus a divisor of level $N$, $D \in Div(X)_N$,
is represented by an open covering $\{U_{\alpha} \}$ of $X_N$,
and local equations $f_{\alpha} \in \cK^*_N (U_{\alpha})$,
such that
\be
f_{\alpha_1}/f_{\alpha_2} \in \cO^*_{X_N} (U_{\alpha_1} \cap U_{\alpha_2}) \ \mbox{for all} \ \alpha_1 , \alpha_2
\ee
Two such collections $\{U_{\alpha}, f_{\alpha}\}$
and $\{V_{\beta}, g_{\beta} \}$
represent the same divisor if and only if there exists 
a common refinement
$\{W_{\gamma}\}$,
and elements $u_{\gamma} \in \cO^*_{X_N}(W_{\gamma})$,
such that for all $\alpha, \beta, \gamma$,
\be
W_{\gamma} \subseteq U_{\alpha} \cap V_{\beta}
\Rightarrow
f_{\alpha } = u_{\gamma} \circ g_{\beta} \ \mbox{on} \ W_{\gamma}
\ee

\section{Associated invertible module}
\label{sec8.3}
For a divisor $D \in Div (X)_N$,
represented by $\{ U_{\alpha}, f_{\alpha}\}$,
we define the subsheaf $\cO_{X_N}(D) \subseteq \cK_N$
by $\cO_{X_N}(D) |_{U_{\alpha}} = 
f^{-1}_{\alpha} \circ \cO_{X_N}|_{U_{\alpha}}$.

The sheaf $\cO_{X_N}(D)$
is an $\cO_{X_N}$-submodule of $\cK_N$, 
in the sense that the operations of multiplication 
and contraction of $\cK_N$ satisfy,
\be
\cO_{X_N} \circ \cO_{X_N}(D) \circ \cO_{X_N} , \quad
(\cO_{X_N}, \cO_{X_N}(D)) , \quad
(\cO_{X_N}(D), \cO_{X_N}) \subseteq \cO_{X_N}(D)
\ee
Moreover, $\cO_X(D)$
is an  \emph{invertible} $\cO_{X_N}$-sheaf,
in the sense that it is locally isomorphic to $\cO_{X_N}$.

\section{Effective divisors}
\label{sec8.4}
A divisor $D \in Div(X)_N$ is said to be \emph{effective} 
if any one of the following equivalent conditions holds:
\begin{itemize}
\item[(8.4.1)] If $\{U_{\alpha}, f_{\alpha} \}$ represents $D$, 
       then $f_{\alpha} \in \cO_{X_N}(U_{\alpha})$.
\item[(8.4.2)] $\cO_{X_N} \subseteq \cO_{X_N}(D) \subseteq \cK_N$
\item[(8.4.3)] $\cO_{X_N}(-D) \subseteq \cO_{X_N}$
\end{itemize}

Thus for effective divisor $D$, the sheaf 
$\cO_{X_N}(-D)$ is an $\cO_{X_N}$-ideal,
which is homogeneous (and locally principal).
The quotient sheaf $\cO_D = \cO_{X_N}/ E(\cO_{X_N}(-D))$
is supported in the closed set $supp(D) \subseteq X_N$,
and $(supp(D)), \cO_D)$ is a closed subscheme of 
$(X_N, \cO_{X_N})$.

Note that for $M \geq N$,
we have via (8.1.2) a homomorphism of abelian groups
\[
(\pi_N^M)^{\sharp}: Div (X)_N \rightarrow Div(X)_M
\]

\section{Divisors}
\label{sec8.5}
We let $B(X)$
denote the set of all $D=\{D_N\}$,
$D_N \in Div(X)_N$,
which are \emph{monotone}:
\be
\mbox{for} \ M \geq N,
(\pi_N^M)^{\sharp}D_N = D_M \circ d_N^M ,
\ \mbox{with} \ d_N^M \in Div(X)_M \ \mbox{effective};
\ee
and are \emph{bounded}:
\be
\mbox{there exists} \ N_0 \in \cN, \ \mbox{and} \ 
d_{N_0}\in Div(X)_{N_0} \ \mbox{such that for} \
M \geq N_0,
\ee
\[D_M \circ (\pi^M_{N_0})^{\sharp}d_{N_0} \ \mbox{is effective}. 
\]

For $D = \{D_N\}$ , $D'=\{D'_N\}$
 in $B(X)$,
 we write $D \geq D'$ if and only if

\noindent (8.5.3) for all $N$, there exists $\tau(N) \geq N$,
such that for all $M \geq \tau(N)$ we have 
$(\pi_N^M)^{\sharp}D_N = D'_M \circ d_N^M$,
with $d_N^M\in Div(X)_M$
effective.

We write $D \approx D'$ if $D \geq D'$
and $D \leq D'$.
The relation $\approx$ is an equivalence relation on 
$B(X)$.
The collection of $\approx$-equivalence classes is defined 
to be the set of divisors of $X$.

\[
Div(X) = B(X)/\approx
\]
The structures of abelian groups on $Div(X)_N$, 
induce a structure of an abelian group on $Div(X)$.
The relation $\geq$ on $B(X)$ induces a partial-order on 
$Div(X)$,
making it an ordered abelian group.

\section{Example: $Div(\overline{spec \mathbb Z})$}
\label{sec8.6}
Let $\cX = \overline{spec \mathbb Z} = 
\{\cX_N\}_{N \geq 2 \ \mbox{square free}}$
cf. (6.4),
and fix $N = \prod\limits_{i=1}^{k} p_i$.

The open sets of $\cX_N$ are 
$U_M = spec \mathbb Z[\frac{1}{M}]$, any $M$,
and $U_M \coprod \{\eta\}$ for $M|N$.
A cofinal system of coverings of $\cX_N$ is given by
\be
V_{\eta} = U_N \coprod \{\eta\} , \quad
V_{p_i} = U_{N \cdot M_i} \coprod \{p_i\}, i = 1, \ldots, k
\ee
The sheaf $K_N^*$ is the constant sheaf $\mathbb Q^*$. 
Thus $D \in Div(X)_N$ is represented by 
$\{f_{\eta}, f_{p_1}, \ldots , f_{p_k}\} \in (\mathbb Q^*)^{k+1}$, 
with
\be
f_{p_i}/f_{p_j} \in \cO^*_{\cX_N}(V_{p_i} \cap V_{p_j})=
\mathbb Z[\frac{1}{NM_iM_j}]^*
\ee
\be
f_{\eta}/f_{p_i} \in \cO^*_{\cX_N}(V_{\eta}\cap V_{p_i})=
\mathbb Z[\frac{1}{NM_i}]^*
\ee
and with
\be
f_{p_i}\in \mathbb Q^* \ \mbox{modulo} \
\cO^*_{\cX_N}(V_{p_i}) = \mathbb Z[\frac{1}{M_i \cdot N/p_i}]^*
\ee
\be
f_{\eta} \in \mathbb Q^* \ \mbox{modulo} \
\cO^*_{\cX_N}(V_{\eta}) = \{ \pm 1\}
\ee

It follows that we have identification,
\be
Div(\cX)_N \cong (\bigoplus\limits_{p|N}\mathbb Z \cdot [p])
\oplus \log(\mathbb Q^+) \cdot [\eta]
\ee
\[D = \{ f_{\eta}, f_{p_1}, \ldots , f_{p_k}\}/\approx \mapsto 
\sum\limits_p \nu_p(D) \cdot [p]
\]
\[\nu_{\eta}(D) = -\log|f_{\eta}|, \quad
\nu_p(D) = \nu_p(f_p) \ \mbox{for} \
p|N,\]
\[\hspace{3cm}
(\nu_p = p\mbox{-adic valuation})
\]
The homomorphism 
$(\pi_N^M)^{\sharp}: Div(\cX)_M \rightarrow Div(\cX)_N$
is given for $M = N \cdot \prod\limits_{j=1}^l q_j$ by
\be
(\pi_N^M)^{\sharp}: (\bigoplus\limits_{p|N}\mathbb Z[p])
\oplus \log(\mathbb Q^+) \cdot [\eta] 
\rightarrow (\bigoplus\limits_{p|M}\mathbb Z[p])
\oplus \log(\mathbb Q^+) \cdot [\eta] 
\ee
\[\sum\limits_{p|N}\nu_p(D) \cdot [p]
 - \log|f_{\eta}|\cdot[\eta] \mapsto
\sum\limits_{p|N}\nu_p(D) \cdot [p]
+\sum\limits_{j=1}^l\nu_{q_j}(f_{\eta}) \cdot [q_j] -
\log|f_{\eta}|\cdot [\eta]
\]

The sequence $D_N \in Div(\cX)_N$ 
is monotone if and only if for all $p$, 
including $p = \eta$,
the filter $\nu_p(D_N)$
is decreasing. We obtain an identification,
\be
Div(\cX) = \left(B(\cX)/\hspace{-.25cm} \ \approx\right) \cong (\bigoplus\limits_{p}\mathbb Z[p])
\oplus \log(\mathbb R^+) \cdot [\eta] 
\ee
\[D = \{D_N\}/\approx \quad \mapsto \sum\limits_p \nu_p(D)\cdot [p] +
\nu_{\eta}(D) \cdot[\eta]
\]
\[\mbox{with} \ \nu_p(D) = \lim\limits_N \nu_p(D_N)
\]
For $f \in \mathbb Q^*$ we have
\be
div(f) = \sum\limits_p \nu_p (f) \cdot [p] + 
\nu_{\eta}(f)[\eta]
\ee
\[
\nu_{\eta}(f) = - \log|f|_{\eta} , \quad
\nu_p(f) = p\mbox{-adic valuation.}\]

We obtain the exact sequence (the ''completed-cohomology-sequence''
 associated to (8.1.3),
\[\begin{array}{cl}
* \rightarrow \{\pm 1\} \rightarrow \mathbb Q^* \xrightarrow{div}&
Div(\cX) \rightarrow \mathbb R^+ \rightarrow *\\
& D \mapsto e^{\nu_{\eta}(D)}\cdot \prod\limits_{p}p^{\nu_p(D)}
\end{array}
\]

\section{Conjectures}
\label{sec8.7}
Let $\cX = \{\cX_{N,M}\} = \overline{spec \mathbb Z} \prod\limits_{\mathbb F[\pm 1]} \overline{spec \mathbb Z}$,
and let $\cD = Div(\cX)$.

There should exist an intersection pairing
\be
( \ , \ ) : \cD \times \cD \rightarrow \mathbb R = \log (\mathbb R^+)
\ee
with the usual properties:
\begin{itemize}
\item[(8.7.2)] Bilinear: $(D_1 \circ D_2 , C) = (D_1, C) +(D_2 , C)$
\item[(8.7.3)] Symmetric: $(D,C) = (C,D)$
\item[(8.7.4)] Linear-equivalence invariant: 
$(D \circ div(f), C ) = (D,C)$ for
$f \in \cK^*_{N,M}(\cX_{N,M})$.
\end{itemize} 
\setcounter{equation}{4}

Moreover, there should exist a collection of functions 
${\mbox{\Large \it{f}}} \subseteq \cW(\mathbb R^+)$,
\be
\begin{array}{cl}
\cW(\mathbb R^+) = \{ &
f: \mathbb R^+ \rightarrow \mathbb R, f(x) \ \mbox{and} \ f'(x) \
\mbox{are continuous except for}\\
& \mbox{ finitely many points} \
  x_1, \ldots , x_N , \ \mbox{where }\\
& f(x_i \pm 0) \
\mbox{and} \ f'(x_i\pm 0) \ \mbox{exist, and} \ 
f(x) \ \mbox{compactly} \\
& \mbox{ supported
 ( or} \ |f(x)| \ \mbox{and} \ |f'(x)| \ \mbox{are} \\
& \leq c \cdot \max (x, x^{-1})^{-(\frac{1}{2}+ \varepsilon)},
\varepsilon > 0 ) \}
\end{array}
\ee
and a mapping
\be
Fr: {\mbox{\Large \it{f}}} \rightarrow Div (\cX)
\ee
such that for $f_1, f_2 \in {\mbox{\Large \it{f}}}$, 
the intersection pairing is given by
\be
(Fr(f_1), Fr(f_2) ) = \hat{f}_1(0) \hat{f}_2(1) +
     \hat{f}_1(1)\hat{f}_2(0) -
\sum\limits_{\xi(\alpha) =0} \hat{f}_1(\alpha) \cdot 
\hat{f}_2(1-\alpha)
\ee
with $\hat{f}(\alpha)= \int\limits_0^{\infty} f(x) x^{\alpha}
\frac{dx}{x}$,
the Mellin transform,
and the sum is over the non-trivial zeros of Riemann's zeta 
$\xi(\alpha)$.

Contemplating the geometric analog, 
one might further conjecture that the mapping $Fr$ is 
$\mathbb Z$-linear, and satisfy
\be
Fr(f_1 * f_2) = Fr(f_1) * Fr(f_2)
\ee
with ordinary convolution of functions
$f_1 * f_2 (y) = 
\int\limits_0^{\infty} \frac{dx}{x} f_1(x) f_2(\frac{y}{x})$
and (non-commutative) composition of divisors
\[ D_1 * D_2 = pr_{1,3}((D_1 \times \overline{spec \mathbb Z})
\cdot(\overline{spec \mathbb Z} \times D_2))
\] and that moreover
\be
Fr(f^{\natural}) = Fr(f)^{\natural}
\ee
with
\be
f^{\natural}(x) = x^{-1} \cdot f(x^{-1}),
\ \mbox{and} \
D^{\natural} = J^{\sharp}(D)
\ee
 \[ \mbox{where} \
J: \cX \xrightarrow{\sim} \cX \ \mbox{is interchanging factors}
\]
Thus there could be a family of functions 
$f_i \in \cW(\mathbb R^+)$,
and divisors $D_i \in Div(\cX)$, $i \in I$,
such that ${\mbox{\Large \it{f}}}$ is the $\mathbb Z$-spann
\be
{\mbox{\Large \it{f}}} = \left\{ \sum\limits_{i,j} m_{i,j} f_{i_1} * \cdots f_{i_k} *
f_{j_1}^{\natural}* \cdots * f_{j_l}^{\natural}, \quad 
m_{i,j} \in \mathbb Z \right\} 
\ee
and 
\be
Fr\left(\sum\limits_{i,j} m_{i,j} f_{i_1}* \cdots f_{i_k}*
f_{j_1}^{\natural}* \cdots * f_{j_l}^{\natural}
\right) =
\sum\limits_{i,j} m_{i,j} D_{i_1}* \cdots * D_{i_k} * 
D_{j_1}^{\natural}* \cdots *D_{j_l}^{\natural}
\ee
If ${\mbox{\Large \it{f}}}$ is rich enough to localize the zeros of Riemann's zeta 
(i.e. is such that if there is a zero $\alpha_0$ of zeta,
$\xi(\alpha_0) =0$,
with $\alpha_0 \neq 1- \widebar{\alpha_0}$, 
then there is $f \in {\mbox{\Large \it{f}}}$,
with $\sum\limits_{\xi(\alpha)=0} \hat{f}(\alpha) \cdot
\overline{\hat{f}(1-\widebar{\alpha})}$ negative),
then the Riemann Hypothesis would follow from the Castelnuovo-Severi inequality,
itself a consequence of the Riemann-Roch theorem for $\cX$,
cf. \cite{H89}.

\end{document}